\documentclass[12pt]{article}
\usepackage{graphicx, fullpage}
\usepackage{amsmath, amssymb,verbatim,amsthm}
\usepackage{natbib}
\setcitestyle{authoryear,round}
\usepackage{xcolor}
\usepackage{setspace}
\usepackage{caption}
\usepackage{dsfont}
\usepackage{subcaption}
\usepackage{hyperref}
\usepackage{listings}

\usepackage{amsfonts} 
\DeclareMathSymbol{\shortminus}{\mathbin}{AMSa}{"39}

\newtheorem{theorem}{Theorem}[section]

\newtheorem{assumption}{Assumption}[section]

\newtheorem{lemma}{Lemma}[section]
\newtheorem{example}{Example}[section]
\newtheorem{corollary}{Corollary}[section]

\newtheorem{remark}{Remark}[section]

\newcommand{\cov}{\mathop{\rm {\mathbb C}ov}\nolimits}%
\newcommand{\var}{\mathop{\rm {\mathbb V}ar}\nolimits}%
\newcommand{\Ex}{\mathrm{E}}

\newcommand{\R}{\mathds{R}}
\newcommand{\C}{\mathds{C}}
\newcommand{\N}{\mathds{N}}
\newcommand{\Z}{\mathds{Z}}

\newcommand{\eps}{\varepsilon}

\newcommand{\BS}[1]{\boldsymbol{#1}}

\newcommand{\ind}{\mathds{1}}

\newcommand{\gu}[1]{v(#1)}
\renewcommand{\tilde}{\widetilde}

\title{Inverse covariance operators of multivariate nonstationary time series}

\author{{Jonas Krampe and Suhasini Subba Rao}}
\date{\today}

\begin{document}

\maketitle

\begin{abstract}
\vspace{0.5em}
For multivariate stationary time series many important properties, such as partial correlation, graphical models and autoregressive representations are encoded in the inverse of its spectral density matrix. This is not true for nonstationary time series, where the pertinent information lies in the inverse infinite dimensional covariance matrix operator associated with the multivariate time series. This necessitates the study of the covariance of a multivariate nonstationary time series and its relationship to its inverse. 
We show that if the rows/columns of the infinite dimensional covariance matrix decay at a certain rate then the rate 
(up to a factor) transfers to the rows/columns of the inverse covariance matrix. This is used to obtain a nonstationary autoregressive representation of the time series and a Baxter-type bound between the parameters of the autoregressive infinite representation and the corresponding finite autoregressive projection. The aforementioned results lay the foundation for the subsequent   analysis of locally stationary time series. In particular, we show that smoothness properties on the covariance matrix transfer to (i) the inverse covariance (ii) the parameters of the vector autoregressive representation and (iii) the partial covariances. All results are set up in such a way that the constants involved depend only on the eigenvalue of the covariance matrix and can be applied in the high-dimensional settings with non-diverging eigenvalues.


\end{abstract}

\section{Introduction}

Several important properties in multivariate analysis are encrypted
within the inverse covariance of the underlying random vector. For
example, the
partial correlation, regression parameters and the network
corresponding to the (Gaussian) graphical model. For multivariate
time series the covariance is now an infinite dimensional matrix.
Nevertheless, analogous to classical multivariate analysis many
interesting properties in time series are encoded in the inverse infinite dimensional
variance matrix. They include (i) the partial covariance
between different components of time series after conditioning on the
other time series (ii) time series graphical models which takes into
account the conditional relationships over the entire time series and
(iii) vector autoregressive representations which yield information on
Granger causality. 
For stationary time
series, however, it is rare to directly deduce these relationships from the
inverse covariance, as these quantities have an equivalent  representation in terms of the finite dimensional inverse spectral density matrix corresponding to the autocovariance of the time series.
For example, the partial covariance can be expressed in terms of the
partial spectral coherence (which is a function of the inverse
spectral density matrix; see, \cite{b:pri-81},
Chapter 9.2). The stationary time series graphical model can be
deduced from the zero and non-zeroes of the inverse spectral density
matrix  (see, \cite{p:dah-00b}) and the vector autoregressive
regressive representation can be deduced from the causal factorisation
of the inverse spectral density matrix (see \cite{wiener1958prediction}). 
However, once one moves away from stationarity, a rigorous
understanding of the above properties can only be achieved by directly
studying the inverse of the infinite dimensional covariance matrix
(and its relationship to the corresponding covariance). This is the main objective of this paper, which we make precise below.

Let $\{X_{t} = (X_{t}^{(1)},\ldots,X_{t}^{(p)})^{\top};t\in
\mathbb{Z}\}$ denote a $p$-dimensional multivariate time series  with
$p\times p$-dimensional covariance matrix
$C_{t,\tau}=\cov[X_{t},X_{\tau}]$ for all $t,\tau\in \mathbb{Z}$.  Using $\{C_{t,\tau}\}_{t,\tau}$ we 
 define the linear operator or, equivalently, infinite dimensional  matrix ${\BS C} = (C_{t,\tau};\tau,t\in
\mathbb{Z})$. Under suitable conditions on ${\BS
  C}$, the inverse  ${\BS D} = {\BS  C}^{-1} = (D_{t,\tau};t,\tau\in
\mathbb{Z})$ exists. \cite{p:bas-sub-21}, Section 2, show
that a graphical model for nonstationary time series  can be
defined from the structure of $\BS D$ (based on zero, Toeplitz
and non-Toeplitz submatrices in $\BS D$). This general framework does not impose any conditions  on the nonstationary structure of the time
series. However, in order to learn the network from data
\cite{p:bas-sub-21} focus
on locally stationary time series; by now a widely accepted and used class of
nonstationary time series. Specifically,  smoothness
conditions are placed on the inverse covariance $\BS D$, and the subsequent
analysis is done under these conditions. 
However, most
locally stationary conditions are stated in terms of the covariance
rather than the inverse covariance. This leads to the question 
``do smoothness conditions on $\BS C$ transfer to smoothness on $\BS
D$?'' and provided the initial motivation for this paper. It naturally  
lead to further questions on the "transfer" of smoothness on
$\BS C$ to (a) vector autoregressive representations and (b) the partial covariance. 
Therefore, our aim is to develop a suite of tools that answer such questions. 
To the best of our knowledge there exists very few results in
this area. One notable exception is the recent work of \cite{p:din-zho-21}, but the aims and results in their work are different to those of this paper. 
\cite{p:din-zho-21} specifically focus on the
univariate nonstationary time series $(X_1,\dots,X_n)$ (with $n\rightarrow \infty$). They show that there exists an autoregressive representation of increasing order over the time points, whose coefficients decay at a certain rate. The results are used to test for correlation stationarity.  
In contrast, we work within 
the multivariate time series framework, and allow for both 
low and high dimensional time series. The latter case is important because often
to make meaningful conditional statements about components in the time series (in terms of Granger causality and conditional covariance)  the number of time series included in the analysis may need to be extremely large.
We summarise the main results below.

In order to reconcile $\BS C$ and its inverse $\BS D$,
in Section \ref{sec:inverserate} we show if 
$\|C_{t,\tau}\|_2\leq
K|t-\tau|^{-\kappa}$ for $t\neq \tau$ and some $\kappa>1$ ($\|\cdot\|_2$
denotes the induced $\ell_2$/spectral norm), then
$\|D_{t,\tau}\|_2\leq \mathcal{K}(1+\log
  |t-\tau|)^{\kappa}(|t-\tau|)^{-\kappa+1}$. 
This leads to a nonstationary VAR$(\infty)$ representation of the
time series $\{X_{t}\}_{t}$ where the corresponding VAR parameters decay 
at the same rate. We use this result to obtain a
Baxter-type bound between the parameters of 
autoregressive infinite representation and the corresponding finite autoregressive projection. It is noteworthy that  
the constant $\mathcal{K}$ depends only on the eigenvalues of $\BS C$, but not on the dimension $p$. Hence, if the eigenvalues of 
$\BS C$ do not grow with dimension $p$, these results hold for arbitrary dimension.

The results in Section \ref{sec:inverserate} are instrumental to
proving the results in Section \ref{sec:LS}, where we focus on 
locally stationary time series. In terms of second order structure, a
time series is called second
order locally stationary if its covariance structure can locally be
approximated by a smooth function $\BS C(u)$. We show in Section \ref{sec:LS0}
that $\BS C(u)$ is an autocovariance of a stationary time series. 
 In Section \ref{sec:LS1} we
show that locally stationary conditions based on the covariance
structure imply that its inverse covariance can locally be approximated by a smooth function $\BS D(u)$, which is the inverse autocovariance 
of a stationary time series i.e. 
$\BS D(u) = \BS C(u)^{-1}$.
We use this result to show that the parameters of the vector autoregressive representation of the time series can be approximated by a smooth function.
Finally, in Section \ref{sec:LS2}, we show that the smoothness conditions on the nonstationary covariance transfer to smoothness conditions on the partial covariances. We use this result to justify using an estimator of the local spectral density function to 
estimate the local partial spectral coherence  
(as was done in \cite{p:par-14}) and the local partial correlation. 
The proof of the results can be found in the Appendix. 
\section{Rate of decay of the inverse covariance}\label{sec:inverserate}

\subsection{Notation and assumptions}

In order to derive the results in this paper we need to define
the space on which the operator ${\BS C}$ is acting. This requires the
following notation. 

Let  $\R$ denote the real numbers, $\Z$ all (positive and negative) integers and $\N$ strictly positive integers. 
For $u,v\in \mathbb{R}^{p}$ let $\langle u,v\rangle = u^{\top}v$ and 
$\|v\|_2$ denote the Euclidean distance. 
We use $\ell_2$ and $\ell_{2,p}$ to denote the sequence spaces 
$\ell_{2} = \{u=(\ldots,u_{-1},u_{0},u_{1},\ldots)
; u_{j}\in \mathbb{R}\textrm{ and }\sum_{j\in \mathbb{Z}}u_{j}^{2}<\infty\}$
 and 
$\ell_{2,p}=\{v=(\ldots,v_{-1},v_{0},v_{1},\ldots) ;v_{j}\in \mathbb{R}^{p}
\textrm{ and }\sum_{j\in \mathbb{Z}}\|v_{j}\|_{2}^{2}<\infty\}$. 
On the spaces $\ell_{2}$ and $\ell_{2,p}$ we define the two inner products
$\langle u,v \rangle = \sum_{j\in \mathbb{Z}}u_{j}v_{j}$ (for  $u,v\in \ell_{2}$) and 
$\langle x,y \rangle = \sum_{j\in \mathbb{Z}}\langle x_{j},y_{j}
\rangle$ 
(for $x=(\ldots,x_{-1},x_{0},x_{1},\ldots),y=(\ldots,y_{-1},y_{0},y_{1},\ldots)\in \ell_{2,p}$).
For $x\in \ell_{2,p}$, let $\|x\|_{2} = \langle x,x\rangle$. Furthermore, for  $x \in \ell_{2,p}$ and $s\in \Z, a \in {1,\dots,p}$, we use 
$x_s^{(a)}$ to denote the $s$th element of the $a$th (column) space. Suppose $\{A_{s_1,s_2}\}_{s_1,s_2}$ are $p_1\times p_2$-dimensional matrices, using
this we define the infinite
dimensional matrix $\BS A=(A_{s_1,s_2};s_1,s_2\in \mathbb{Z})$. Under
suitable conditions on $\BS A$, $\BS A$ is a
linear operator $\BS A : \ell_{2,p_1}\to \ell_{2,p_2}$ in the sense that if $\BS Ax = y$, then $y = (\ldots,y_{-1},y_{0},y_{1},\ldots)$ where for all
$t\in \Z$, $y_{t}\in \mathbb{R}^{p_2}$ and
 $y_{t} = \sum_{\tau\in \mathbb{Z}}A_{t,\tau}x_{\tau}$. 
 Furthermore, we define $\|\BS A\|_2=\sup_{\|x\|_2=1} \|\BS A x \|_2$.
 All operators are written in bold uppercase letters.

\begin{assumption} \label{ass}
Let $\gu\cdot =
\max(1,|\cdot|)$. 
\begin{itemize}
\item[(i)] The covariance operator is positive definite with $\lambda_{\sup} =
  \sup_{v \in \ell_{2,p}, \|v\|_{2}=1}\langle v, \BS Cv \rangle<\infty$ and 
$0<\lambda_{\inf} =  \inf_{v \in \ell_{2,p}, \|v\|_{2}=1}\langle v, \BS Cv \rangle.$
    \item[(ii)] There exists some $\kappa>1$ such that for all $t \not =\tau$ we have for the $p\times p$-dimensional sub-matrices
$$\|C_{t,\tau}\|_2\leq K\gu{t-\tau}^{-\kappa},$$
where $K<\infty$ is some positive constant. 
\end{itemize}
\end{assumption}
Since $\BS C$ is positive definite, the inverse covariance operator
exists with $\BS D=\BS C^{-1} = (D_{t,\tau};t,\tau\in \mathbb{Z})$. We
mention that the condition $\lambda_{\sup}<\infty$ is implied by Assumption
\ref{ass}(ii). 

The results in this paper allow for both low and high
dimensional multivariate time series and the assumptions used are
specifically designed to allow for this. For high dimensional time
series, the condition that the largest eigenvalue is bounded excludes 
time series with dynamic factors but allows for high dimensional sparse
time series.\footnote{By dynamic factors we refer to the common component described in the representation given in \cite{forni2000generalized}. The common component contains (if any) the diverging eigenvalues of the process. The bounded eigenvalues define in this decomposition the so-called idiosyncratic component. Hence, if eigenvalues of $\BS C$ diverge with $p$, the results can be applied to the idiosyncratic component of this decomposition.}
Popular examples include  high dimensional sparse time series regression and vector autoregressive (VAR) models which have recently received
considerable attention; see, for example,  \cite{basu2015}\footnote{Note that the finite sample error bounds derived in \cite{basu2015} for the Lasso express the dependence of the processes also in terms of $\lambda_{\inf}$ and $\lambda_{\sup}$.}, \cite{krampe2021bootstrap}, \cite{krampe2021sparsity}, 
(in the context of stationary VAR models) 
and \cite{ding2017sparse} (for time-varying VAR models). 
The condition that $\lambda_{\inf}>0$
omits co-linearity, where one component in the time series can be perfectly explained by other
components. 
 Assumption \ref{ass}(ii) quantifies the pairwise dependencies between
the components (over time) and is stated in terms of the 
(induced) $\ell_2$-norm $\|\cdot\|_2$ of the $p\times p$ matrices.
However, no conditions are placed on the $\ell_{1}$-norm, which can grow
with dimension $p$ (as sparsity usually does in the sparse regression context). All results in this paper are derived in terms of
the $\|\cdot\|_{2}$-norm. Thus we show that if the pairwise
interactions are controlled in the $\ell_2$ sense as $p$ grows, then
the conditional interactions are also controlled in the $\ell_2$-sense.

Throughout this paper we use $\mathcal{K}$ to denote a generic
constant that only depends on $\lambda_{\inf}, \lambda_{\sup}, K,
\kappa$ and whose value may change from line to line.
We define $\gu\cdot =
\max(1,|\cdot|)$ and $\zeta(j)$ as follows; 
for $|j|\leq 1$ let $\zeta(j)=1$ and for $|j|>1$ let $\zeta(j) =
\log|j|/|j|$.

\subsection{The inverse covariance}\label{sec:inverse}

In the following theorem we obtain a bound on the rate of decay of the
matrices $D_{t,\tau}$ that make up the inverse covariance $\BS
D = \BS C^{-1}$. 
$\BS C$ is a bi-infinite matrix in the sense that the
entries $C_{t,\tau}$ span  $t,\tau\in \mathbb{Z}$. 
We will also
consider the one-sided infinite dimensional matrix 
$
\BS C(-\infty,T) = (C_{t,\tau};t,\tau \leq T).
$
As will be clear later in the paper, the inverse of $\BS C(-\infty,T)$ contains (up to a factor) the 
AR prediction coefficients and the following result will be used to obtain a bound on its rate of 
decay.
\begin{theorem}\label{thm.1}
Under Assumption~\ref{ass}, 
for all $t,\tau\in \Z$ we have
\begin{eqnarray}
\label{eq:Dttau}
\|D_{t,\tau}\|_2\leq \mathcal{K}\zeta(t-\tau)^{\kappa-1},
\end{eqnarray}
where $\mathcal{K}$ is a constant depending on
$K,\kappa,\lambda_{\inf}$, and $\lambda_{\sup}$ only and $\zeta(j) = \gu{\log[\gu{j}]}/\gu{j}$.
For $t,\tau\leq T$
\begin{eqnarray}
\label{eq:DttauInv}
\|[\BS C(-\infty;T)]^{-1}]_{t,\tau}\|_2\leq \mathcal{K}\zeta(t-\tau)^{\kappa-1}.
\end{eqnarray}
\end{theorem}
\begin{proof}
  The key ingredient in the proof is Lemma \ref{lem.demko} (in
  Appendix \ref{sec:proof1}), which bounds the entries of the inverse
  of a banded matrix operator (and is a generalisation of Proposition
  1 in \cite{demko1984decay}). The details of the proof are in Appendix \ref{sec:proof1}.
\end{proof}  

The above result shows that if the pairwise interaction between the
components is bounded with a certain rate in the $\ell_{2}$-sense then the conditional
interactions are also bounded with a certain rate in the
$\ell_{2}$-sense, see  Remark~\ref{rmk.dimension} for a discussion on the role of the dimension $p$. 

\begin{remark}
In the case entries in $\BS C$  decay geometrically or are banded,
then the entries of $\BS D$ decay at a geometric rate.
\end{remark}

\begin{remark}[An alternative representation of the covariance ${\BS
    C}$ and its inverse]\label{remark:alternativegroup}
We recall that we defined ${\BS C}$ as ${\BS C} = (C_{t,\tau};t,\tau\in \mathbb{Z})$, where $C_{t,\tau}$ are
$p\times p$-dimensional matrices. An alternative method for defining
${\BS C}$ is to group the covariances according to component i.e. $\widetilde{{\BS C}} = (C^{(a,b)};1\leq
a,b\leq p)$ where $[C^{(a,b)}]_{t,\tau} = C_{t,\tau}^{(a,b)}=
\cov[X_{t}^{(a)},X_{\tau}^{(b)}]$. $\widetilde{{\BS C}}$ is simply a
permutation of ${\BS C}$, thus $\widetilde{{\BS D}} = \widetilde{{\BS
    C}}^{-1}$ is a permutation of ${\BS D}$. In certain applications, such as nonstationary 
    graphical models or condition covariance between two components of a time series, 
the representations $\widetilde{\BS C}$ and $\widetilde{\BS D}$ may be
more useful in the analysis 
than $\BS C$ and $\BS D$ (see, for example, \cite{p:bas-sub-21}).
\end{remark}

We now compare Theorem \ref{thm.1} with the classical result for
stationary time series. 
For this, suppose $\BS C = (C_{t-\tau};t,\tau\in
\mathbb{Z})$ is a block Toeplitz operator from $\ell_{2,p}$ to $\ell_{2,p}$, where
$\BS C$ satisfies the rate and positive definiteness
conditions in Assumption~\ref{ass}. Then $\BS D=\BS C^{-1} = (D_{t-\tau};t,\tau\in
\mathbb{Z})$ exists and is also a block Toeplitz operator. For block Toeplitz operators \cite{cheng1993baxter,meyer2015vector} work with a global condition on the sequence $(C_s)_{s\in \Z}$ instead of the individual one used in this paper. They showed that if the global condition 
$\sum_{s \in \Z} (1+|s|^\kappa) \|C_s\|_2<\infty$ holds, then  $\sum_{s \in \Z} (1+|s|^\kappa) \|D_s\|_2<\infty$. 
The global condition implies for all $t,\tau\in \mathbb{Z}$ that $\|C_{t-\tau}\|_2 \leq K \gu{t-\tau}^{-\kappa}$ and $\|D_{t-\tau}\|_2 \leq K \gu{t-\tau}^{-\kappa}$. 
Conversely, the individual condition that yields this global condition is $\|C_{t-\tau}\|_2 \leq K \gu{t-\tau}^{-\kappa-1-\eps}$. In summary, even for block Toeplitz matrices, at the individual level if $\|C_{t-\tau}\|_2 \leq K \gu{t-\tau}^{-\kappa-\eps}$ then the above arguments yield
\begin{eqnarray}
\label{eq:BSttau}
\|D_{t-\tau}\|_2 \leq K \gu{t-\tau}^{-\kappa+1},
\end{eqnarray}
which is (without the log-factor) same as the rate derived in Theorem \ref{thm.1}. 
To the best of our knowledge, it is an open question if this rate at the individual level for the inverse can be improved for stationary as well general nonstationary time series.

\subsection{Vector Autoregressive representation and Baxter's
  inequality} \label{sec:2.baxter}

It is well known that for stationary time series the
entries of $\BS C(-\infty,T)$ are closely related to vector
autoregressive (VAR$(\infty)$)
parameters of the underlying time series. The same is true for
nonstationary time series. Precisely, under Assumption \ref{ass} and
by using the projection theorem 
the bottom row of $\BS C(-\infty,T)^{-1}$ contains the VAR$(\infty)$
coefficients in the linear projection of $X_{T}$ onto the space
spanned by $\overline{\textrm{sp}}(X_{T-1},X_{T-2},\ldots)$ i.e.,
\begin{eqnarray}
\label{eq:ARinfty}
X_{T} = \sum_{j=1}^{\infty}\Phi_{T,j}X_{T-j} +\varepsilon_{T}, \textrm{ where
  }\Phi_{T,j} =  -([\BS
  C(-\infty,T)^{-1}]_{T,T})^{-1}[\BS
  C(-\infty,T)^{-1}]_{T,T-j},
\end{eqnarray}
where $\varepsilon_{T}$ is uncorrelated with $\{X_{T-j}\}_{j=1}^{\infty}$.
Substituting the bound in Theorem \ref{thm.1} 
into (\ref{eq:ARinfty}) gives 
\begin{eqnarray}
\label{eq:ARrate}
\|\Phi_{T,j}\|_2 \leq \mathcal{K}\zeta(t-\tau)^{\kappa-1}.
\end{eqnarray}
In practice, it is often not possible to estimate the infinite number of AR parameters from a finite data set. Therefore one often estimates the parameters of the projection of $X_{T}$ onto the finite past $\overline{\textrm{sp}}(X_{T-1},\ldots,X_{T-d})$ i.e.,
\begin{eqnarray}
\label{eq:ARfinite}
X_{T} = \sum_{j=1}^{d}\Phi_{T,d,j}X_{T-j} +\varepsilon_{T,d}.
\end{eqnarray}
The above is analogous to the best fitting
VAR$(d)$ parameters for stationary time series.  In stationary time
series the difference between the finite past projection and the
corresponding infinite past projection is called the Baxter
inequality;
see Section~6 in \cite{hannan1988statistical}, \cite{cheng1993baxter}, and \cite{meyer2015vector}.
In the same spirit, we now obtain a Baxter-type inequality for nonstationary multivariate time series,   
between the VAR$(\infty)$ coefficients $\{\Phi_{T,j}\}_{j}$ and the finite prediction coefficients 
$\{\Phi_{T,d,j}\}_{j}$. 

The coefficients 
$\{\Phi_{T,d,j}\}_{j}$ are embedded in the bottom row of the finite
dimensional matrix $\BS C(T-d,T)^{-1}$ where
$\BS C(T-d,T) = (C_{t,\tau};T-d+1\leq t,\tau\leq T)$. Thus the
coefficients $\{\Phi_{T,j}\}_{j}$ and $\{\Phi_{T,d,j}\}_{j}$ are
connected through $\BS C(T-d,T)$ and $\BS C(-\infty,T)$ and their
inverses. Due to this connection we use Theorem \ref{thm.1} and the
block operator inverse identity (see equation (\ref{eq:inverseblock1}) in Appendix \ref{appendix:back}) to prove the result below. 

\begin{theorem}[Baxter type
  inequality]\label{lemma:onesided}
Suppose Assumption \ref{ass} holds with $\kappa>3/2$. Let $\{\Phi_{T,j}\}_{j}$ and
$\{\Phi_{T,d,j}\}$ be defined as in (\ref{eq:ARinfty}) and
(\ref{eq:ARfinite}) respectively. Then for $d\in \N, j=1,\dots,d$ we have 
\begin{eqnarray}
\label{eq:baxter1}
\sup_{T}\|\Phi_{T,d,j}-\Phi_{T,j}\|_{2} \leq \mathcal{K}\zeta(d)^{\kappa-3/2}\zeta(d-j)^{\kappa-3/2}.
\end{eqnarray}
Furthermore, if Assumption \ref{ass} holds with $\kappa>5/2$ we have 
\begin{eqnarray}
\label{eq:baxter2}
\sup_{T}\sum_{j=1}^{d}\|\Phi_{T,d,j}-\Phi_{T,j}\|_{2} \leq \mathcal{K}\zeta(d)^{\kappa-3/2}.
\end{eqnarray}
\end{theorem}
\begin{proof}
  In Appendix \ref{sec:proof1}.
\end{proof}
  
Inequality (\ref{eq:ARrate}) and Theorem~\ref{lemma:onesided} are
related to Theorem 2.4 in \cite{p:din-zho-21}, who obtain
autoregressive approximations for nonstationary univariate time
series. However, it is important to note that there are some
differences in the autoregressive representations derived in both
papers. The autoregressive representation derived in
\citep{p:din-zho-21} is based on the finite vector $(X_1,\dots,X_n)$
and their aim is to build an autoregressive representation of
increasing order over the time points of the data vector, i.e., $X_i$
is represented as an AR$(i-1)$ model. In contrast, we derive an
autoregressive representation of a time series $\{X_t; t \in \Z\}$
where each time point has a  $VAR(\infty)$ representation. 
In the stationary context, building an autoregressive representation
of an increasing order relates to the Cholesky decomposition of
$\var(X_1,\dots,X_n)^{-1}$ where the $i$th model is given by the $i$th
line. The AR$(\infty)$ model using the entire time series can be considered
as a limit of this, see Section~2 in \cite{krampe2021estimating} for
further discussion. With this  difference in mind, we now compare the
rates in Section \ref{lemma:onesided} with the results in Theorem~2.4
in \cite{p:din-zho-21}. Their decay rate for the autoregressive
coefficients matches with that derived in \eqref{eq:ARrate}. In terms
of Baxter's inequality, they show $\max_{T>b} \max_{1\leq j \leq b} |
\Phi_{T,T-1,j}-\Phi_{T,b,j}| \leq C (\log b)^{\kappa-1} b
^{-\kappa+3}$. Using Theorem~\ref{lemma:onesided} we compare the
coefficients of the two finite AR models (order $T-1$ and order $b$),
and obtain tighter bounds for their result. To be precise 
\begin{align*}
    \max_{T} \| \Phi_{T,T-1,j}-\Phi_{T,b,j}\|_2 &\leq 
    \max_{T}  (\| \Phi_{T,T-1,j}-\Phi_{T,j}\|_2+\| \Phi_{T,j}-\Phi_{T,b,j}\|_2) \\
    &\leq  
        \mathcal{K}\left(\zeta(T-j)^{\kappa-3/2}+\zeta(b-j)^{\kappa-3/2}\right) (\log(b)/b)^{\kappa-3/2}.
\end{align*}
The above leads to the bound 
$\max_{T>b} \max_{1\leq j \leq b} |
\Phi_{T,T-1,j}-\Phi_{T,b,j}|=O(b^{-\kappa+3/2}\log^{\kappa-3/2} b)$ instead of
$O(b^{-\kappa+3}\log^{\kappa-1}b)$ (given in \cite{p:din-zho-21}). 

We now compare Theorem \ref{lemma:onesided} 
to the stationary set-up. 
\cite{meyer2015vector} showed that under the following
global condition on the vector autoregressive parameters
$\sum_{s \in \Z} (1+|s|^\kappa) \|\Phi_s\|_2<\infty$, 
that 
\begin{eqnarray}
\label{eq:Phidsum}
\sum_{j=1}^d (1+j)^\kappa \|\Phi_{d,j}-\Phi_{j}\|_{2} \leq 
\mathcal{K} \sum_{j=d+1}^\infty (1+j)^\kappa \|\Phi_{j}\|_{2},
\end{eqnarray}
noting that we have dropped $T$ as it is not necessary under stationarity.
(\ref{eq:Phidsum}) implies $\sum_{j=1}^d \|\Phi_{d,j}-\Phi_{j}\|_{2}\leq \mathcal{K} d^{-\kappa}$.
Based on the discussion at the end of Section \ref{sec:inverse}, 
at the individual level this means if $\|C_{s}\|_{2}\leq K\gu{s}^{-\kappa-\eps}$,
then $\sum_{j=1}^d \|\Phi_{d,j}-\Phi_{j}\|_{2}\leq \mathcal{K} d^{-\kappa+1}$, whereas Theorem \ref{lemma:onesided} gives 
$\sum_{j=1}^d \|\Phi_{d,j}-\Phi_{j}\|_{2}\leq 
\mathcal{K} d^{-\kappa+3/2}$. Thus stationarity of the time series yields a 
better approximation bound between the finite and infinite AR parameters than the bound in Theorem \ref{lemma:onesided}. 

\section{Locally stationary time series}\label{sec:LS}

The first rigorous treatment of locally stationary time series was
given in  \citep{p:dah-97, p:dah-00a}. This was done by
representing $\{X_{t,T}\}_{t=1}^{T}$ in terms of a Cram\'er representation
$
X_{t,T} = \int_{0}^{2\pi}A_{t,T}(\omega)dZ(\omega),
$
where $\{Z(\omega);\omega\in [0,2\pi]\}$\ is an orthogonal increment process
and the time-varying  transfer function $A_{t,T}(\omega)$ can  locally be approximated  
by the Lipschitz smooth function $A(\omega;\cdot)$ i.e. 
$\|A_{t,T}(\omega)-A(\omega;u)\|_{2}\leq K(|t/N-u|+1/N)$.   This definition immediately leads to certain smoothness properties on the covariance structure of the time series.
More recently, several authors have extended this
definition to nonlinear time series 
cf. \citep{p:dah-sub-06, p:sub-06, p:zhou-09, p:vog-12, p:tru-19, p:dah-19, p:kar-21}.  
In this section, we return, in some sense, to the original formulation of local
stationarity and focus on the locally stationary second order structure. However, unlike \citep{p:dah-97, p:dah-00a}, we work  within the time domain and not the frequency domain.
We start by introducing the locally stationary setting, i.e., we impose certain smoothness conditions on the nonstationary time series. In Section \ref{sec:LS0}
we  obtain bounds on the eigenvalues of the underlying covariance. 
Using Theorem \ref{thm.1}, 
in Section \ref{sec:LS1} we show that smoothness conditions placed on the covariance structure transfer over to the inverse covariance and the parameters in the nonstationary
AR$(\infty)$ representation. Finally, in Section \ref{sec:LS2}, we apply these results to show that the smoothness conditions also transfer to the partial covariances. 

\subsection{Assumptions}

We start by defining an infinite array, where for each $N\in \N$
we associate a (non)stationary multivariate time series $\{X_{t,N};t\in \mathbb{Z}\}$ and 
covariance $C_{t,\tau}^{(N)} = \cov[X_{t,N},X_{\tau,N}]$ (for all
$t,\tau\in \mathbb{Z}$). For each $N$ we 
define the infinite dimensional covariance matrix 
${\BS C}^{(N)} = (C_{t,\tau}^{(N)};t,\tau\in
\mathbb{Z})$.  In the assumptions below we explicitly connect the 
sequence of infinite dimensional covariance matrices 
$\{{\BS C}^{(N)}\}_{N\in \N}$ through
$N$, which plays the role of a
smoothing parameter. We mention that it is standard practice in the
locally stationary literature to define $X_{t,N}$ on a triangular
array i.e. $\{X_{t,N}\}_{t=1}^{N}$. However, to avoid confusion, we do not link $N$ to sample size. It is also worth pointing out that we use $N\in \N$ to simplify the exposition, we could, without loss of generality,  allow $N$ to be a non-integer and define it on
$N\in [\alpha,\infty)$ (for some $\alpha >0$).

\begin{assumption}\label{assum:LS}
\begin{itemize}
\item[(i)] Eigenvalue condition: There exists some $N_0\geq 1$ where 
\begin{eqnarray*}
0 < \lambda_{\inf} \leq \inf_{N\geq N_0}\lambda_{\inf}({\bf C}^{(N)} )
  \leq \sup_{N \geq N_0}\lambda_{\sup}({\bf C}^{(N)} ) \leq \lambda_{\sup} <\infty.
\end{eqnarray*}
\item[(ii)] Covariance decay condition: For all $N$, $t$ and $\tau$
$
\|\BS C_{t,\tau}^{(N)}\|_2 \leq \frac{K}{\gu{t-\tau}^{\kappa}}.
$
\item[(iii)] Smoothness condition: There
exists a Lipschitz continuous matrix function $\{C_{r}(\cdot), r \in
Z\}$ where (a) $C_{r}(u) = C_{-r}(u)^{\top}$, (b) for all $u,v\in \R, r \in \Z$ 
$\sup_{u}\|C_{r}(u)\|_2\leq K/\gu{r}^{\kappa} $, and (c) $\|C_r(u)-C_r(v)\|_2\leq \frac{K|u-v|}{\gu r^\kappa}$,
such that for all $N$
\begin{eqnarray}
\label{eq:LScovariance2}
\|C_{t,\tau}^{(N)} - C_{t-\tau}(t/N)\|_2 \leq 
\frac{K}{\gu{t-\tau}^{\kappa-1}}
\min\left(\frac{1}{N},\frac{2}{\gu{t-\tau}}\right).
\end{eqnarray}
We assume that $\kappa >3$. 
\end{itemize}
\end{assumption}
Note that the above assumptions imply that 
\begin{eqnarray*}
\|C_{t,\tau}^{(N)} - C_{t-\tau}(u)\|_2 \leq
\frac{K}{\gu{t-\tau}^{\kappa-1}}\min\left[
\left(|u-\frac{t}{N}|+\frac{1}{N}\right), 
\frac{2}{\gu{t-\tau}}\right].
\end{eqnarray*}
Furthermore, the sequence $\{C_{r}(\cdot), r \in
Z\}$ defines the infinite dimensional matrix operator $\BS
C(\cdot)=(C_{t-\tau}(\cdot);t,\tau\in \mathbb{Z})$ (from $\ell_{2,p}$ to
$\ell_{2,p}$), where $\BS C(\cdot)$ is block Toeplitz.

Assumption~\ref{assum:LS}(i) and (ii) can be viewed as Assumption~\ref{ass} within the framework of an infinite array. 
Assumption~\ref{assum:LS}(iii) places smoothness conditions on the covariance i.e., the (potentially) non-Toeplitz-operator $\BS C^{(N)}$ can locally be approximated by a block Toeplitz-operator $\BS C(\cdot)$, where the approximation error is determined by the smoothing parameter $N$. The use of 
$\min$ in Assumption~\ref{assum:LS}(iii) is not standard within the locally stationary literature. This arises because the time series $\{X_{t,N}\}_{t}$ is defined on $t\in \mathbb{Z}$ and not $t=1,\ldots,N$ (the typical locally stationary set-up). If 
$|t-\tau|<2N$, then Assumption~\ref{assum:LS}(iii) implies that
$\|C_{t,\tau}^{(N)} - C_{t-\tau}(t/N)\|_2 \leq 
\frac{K}{N\gu{t-\tau}^{\kappa-1}}$ (the classical locally stationary condition). On the other hand, 
if $|t-\tau|\geq 2N$, then the smoothing parameter $N$ does not improve on the individual terms $C_{t,\tau}^{(N)} $ and 
$C_{t-\tau}(t/N)$
(which are extremely small) and we have  
$\|C_{t,\tau}^{(N)} - C_{t-\tau}(t/N)\|_2 \leq 
\frac{2K}{\gu{t-\tau}^{\kappa}}$.
 To distinguish these two cases all the relevant results will be stated with $\min$. 
 
\begin{remark}[The role of dimension $p$] \label{rmk.dimension}
In Assumption \ref{assum:LS} we have not included the dimension $p$ as an additional variable. This is to reduce cumbersome notation. However, it is possible to state Assumption \ref{assum:LS} in terms of uniform bounds over a three dimensional array where the eigenvalues are uniformly bounded over both $N$ and $p$
(and $\BS C^{(N)}$ and $\BS C(u)$ are indexed with $p$ too). 
If these assumptions hold, then the results in this section hold for high dimensional $p$ too. 
\end{remark}

Assumption~\ref{assum:LS} is satisfied by a wide range of locally
stationary time series. In 
Example~\ref{exam:tvMA1} (below) and \ref{exam:tvMA2} we define the time-varying Vector Moving Average (tv-VMA) model and show that this model satisfies Assumption~\ref{assum:LS}. 

\begin{example}[The time-varying vector MA$(\infty)$(tv-VMA) process]\label{exam:tvMA1}
Consider the tv-VMA$(\infty)$
\begin{eqnarray*}
X_{t,N} = \sum_{j=0}^{\infty}\Psi_{t,j}^{(N)}\varepsilon_{t-j} 
 = \sum_{j=1}^{\infty}\Psi_{t,j}^{(N)}\varepsilon_{t-j} + \Psi_{t,0}\varepsilon_{t},
\qquad t\in \mathbb{Z},
\end{eqnarray*}
where $\{\varepsilon_{t}\}_{t}$ are uncorrelated random variables with zero
mean and variance $I_{p}$. In order for the process to be well defined
certain summability or decay conditions need to be imposed on the coefficients $\{\Psi_{t,j}\}$. 
We assume that $\sup_{N\in \N}\sup_{t\in
  \mathbb{Z}}\|\Psi_{t,j}^{(N)}\|_2\leq K\gu{j}^{-\kappa}$. 
With this, we have 
\begin{align*}
C^{(N)}_{t,\tau} &= \cov(\sum_{j=0}^{\infty}\Psi_{t,j}^{(N)}\varepsilon_{t-j},
\sum_{j=0}^{\infty}\Psi_{\tau,j}^{(N)}\varepsilon_{\tau-j})=\sum_{j\in
  \mathbb{Z}}\Psi_{t,j}^{(N)}(\Psi_{\tau,j+\tau-t}^{(N)})^{\top},
\end{align*}
where we set  $\Psi_{t,j}^{(N)}=0$ for $j<0$. Using the above decay
condition on $\Psi_{t,j}^{(N)}$  and Lemma~\ref{lem.SSR.F1} we have
$\|C^{(N)}_{t,\tau}\|_2\leq K\gu{t-\tau}^{\kappa}$; thus Assumption
\ref{assum:LS}(ii) holds.  We now introduce the locally stationary
approximation to  $\{X_{t,N}\}$. Analogous to  \cite{p:dah-97} and \cite{p:dah-pol-06}
(for the case $p=1$), we assume there exists a
Lipschitz continuous matrix function $\Psi_{j}(\cdot)$
where $\sup_{u\in  \mathbb{R}}\|\Psi_{j}(u)\|_2\leq K\gu{j}^{-\kappa}$, $\sup_{u\in  \mathbb{R}}\|\Psi_{j}(u)-\Psi_{j}(v)\|_2\leq K|u-v|\gu{j}^{-\kappa}$, and 
$\|\Psi_{t,j}^{(N)} - \Psi_{j}(t/N)\|\leq K\gu{j}^{-\kappa}/N$.
Using this, we define the stationary process $\{X_{t}(u)\}_{t}$ where 
$X_{t}(u) = \sum_{j=0}^{\infty}\Psi_{t,j}(u)\varepsilon_{t-j}$ which has 
autocovariance $C_{r}(u) = \sum_{j\in
  \mathbb{Z}}\Psi_{j}(u)\Psi_{j+r}(u)^{\top}$ (where we set $\Psi_{j}(u)=0$ for $j<0$). 
Note $\sup_{u}\|C_{r}(u)\|_2\leq K/\gu{r}^{\kappa} $ (this follows from Lemma~\ref{lem.SSR.F1}). Furthermore, under these conditions we have 
\begin{eqnarray*}
\|C^{(N)}_{t,\tau} - C_{t-\tau}(t/N)\|_{2}&\leq& \sum_{j\in \mathbb{Z}}
\|\Psi_{t,j}^{(N)} - \Psi_{j}(t/N) \|_{2}\|\Psi_{\tau,j+\tau-t}^{(N)}\|_{2} \\
&& + \sum_{j\in \mathbb{Z}}
\|\Psi_{j}(t/N) \|_{2}
\big(\|\Psi_{j+\tau-t}(t/N)-\Psi_{j+\tau-t}(\tau/N)\|_{2}\\
&&+\|\Psi_{j+\tau-t}(\tau/N)-\Psi_{\tau,j+\tau-t}^{(N)}\|_{2}
\big)
\\
&\leq & 
\frac{\mathcal{K}}{N}\sum_{j\in \mathbb{Z}}\left(\frac{1}{\gu{j}^\kappa\gu{j+t-\tau}^{\kappa-1}} + \frac{|t-\tau|}{\gu{j}^\kappa\gu{j+t-\tau}^\kappa} \right)\\
&\leq &
\frac{\mathcal{K}}{N\gu{t-\tau}^{\kappa -1}}.
\end{eqnarray*}
Thus Assumption \ref{assum:LS}(iii) holds. We observe that this example illustrates why the rate drops from 
$\kappa$ to $\kappa-1$ in $\|C^{(N)}_{t,\tau}- C_{t-\tau}(t/N)\|_2$; there is an additional  "cost" due to the inclusion of the term $|t-\tau|$. 

In Example \ref{exam:tvMA2} (in Section \ref{sec:LS0})
we show that Assumption \ref{assum:LS}(i) is also satisfied (for sufficiently large $N$).
\end{example}

\subsection{Properties of the locally stationary covariance}\label{sec:LS0}

In this subsection we show that positive definiteness of $\BS C^{(N)}$
transfers to $\BS C(\cdot)$ under the stated smoothness condition. Conversely, we show that also the other direction holds i.e., for a sufficiently large 
$N_0$ positive definiteness of $\BS C(u)$ implies that $\BS C^{(N)}$
is also positive definite (for $N>N_0$).

\begin{theorem}[Positive definiteness of $\BS C(u)$] \label{lem.Cu.pos.def}
Suppose Assumption~\ref{assum:LS} holds. Then, for all $u \in \R$
$\{C_{r}(u)\}_{r}$ is a positive definite sequence where 
$\lambda_{\inf} \leq \lambda_{\inf}({\BS C}(u)) \leq
\lambda_{\sup}({\BS C}(u)) \leq \lambda_{\sup}$.   
\end{theorem}
\begin{proof}
In Appendix \ref{sec:proof2}.
\end{proof}

Under the above theorem, $\{C_{r}(u)\}_{r}$ is a
positive definite sequence. Consequently by Kolmogorov's extension theorem there exists a stationary multivariate time series $\{X_{t}(u)\}_{t\in
  \mathbb{Z}}$ which has $\{C_{r}(u)\}_{r\in \mathbb{Z}}$ as its
autocovariance function. This justifies calling $\{X_{t,N}\}_{t\in
  \mathbb{Z}}$ a ``locally''  second order stationary time series.
A further implication of Lemma~\ref{lem.Cu.pos.def} is that 
the inverse of $\BS C(u)$ exists, which we denote by $\BS D(u)=\BS
C(u)^{-1} = \{D_{t-\tau}(u);t,\tau\in \mathbb{Z}\}$. Like $\BS C(u)$, $\BS D(u)$ is also block Toeplitz and by Theorem~\ref{thm.1} the $p\times
p$-dimension matrix $D_{t-\tau}(u)$ has the bound
\begin{eqnarray}
\label{eq:BSttau2}
\sup_{u}\|D_{t-\tau}(u)\|_2 \leq \mathcal{K} \zeta(t-\tau)^{-\kappa+1}.
\end{eqnarray}

For a given nonstationary time series model, 
Assumption \ref{assum:LS}(i) is difficult to directly verify. However, we now show that given a positive definite sequence  $\{C_{r}(u)\}_{r}$ which satisfies Assumption \ref{assum:LS}(ii,iii), 
then Assumption \ref{assum:LS}(i) holds. For the
univariate case, a similar result is given in \citep[Proposition 2.9]{p:din-zho-21}. 

\begin{theorem}\label{lemma:spectral}
Suppose $\{X_{t,N}\}_{t\in \mathbb{Z}}$ is a locally stationary time
series whose covariance $C^{(N)} = (C_{t,\tau}^{(N)};t,\tau\in
\mathbb{Z})$ satisfies Assumption \ref{assum:LS}(ii,iii). Let 
$f(\omega;u) = \sum_{r\in \mathbb{Z}}C_{r}(u)\exp(ir\omega)$ be the local spectral density. 
If
\begin{eqnarray}
\label{eq:fomegaspectral}
0<\gamma_{\inf} \leq \inf_{u}\inf_{\omega}\lambda_{\min}(f(\omega;u))
\leq \sup_{u}\sup_{\omega}\lambda_{\max}(f(\omega;u))\leq
\gamma_{\sup}<\infty,
\end{eqnarray}
 then there exists a $N_0$, $\lambda_{\inf}$ and
$\lambda_{\sup}$ where for all $N\geq N_0$ we have 
\begin{eqnarray*}
0 < \lambda_{\inf} \leq \lambda_{\inf}({\bf C}^{(N)} ) \leq \lambda_{\sup}({\bf C}^{(N)} ) \leq \lambda_{\sup} <\infty.
\end{eqnarray*}
\end{theorem}
\begin{proof}
In Appendix \ref{sec:proof2}.
\end{proof}

Equipped with the above results, we return to Example \ref{exam:tvMA1}.

\begin{example}[Example ~\ref{exam:tvMA1}, continued]\label{exam:tvMA2}
We define the local spectral density as
\begin{eqnarray*}
f(\omega;u)  = [\sum_{j=0}^{\infty}\Psi_{j}(t/N)\exp(-ij\omega)] [\sum_{j=0}^{\infty}\Psi_{j}(t/N)\exp(ij\omega)]^{\top}.
\end{eqnarray*}
Under the conditions of Example~\ref{exam:tvMA1} we have 
$\sup_{u}\sup_{\omega}\lambda_{\max}(f(\omega;u))\leq \sum_{j\in \Z} K \gu{j}^{-\kappa}=:\gamma_{\sup}<\infty$. 
Furthermore, if we have a non-vanishing filter in the sense 
\begin{eqnarray*}
\inf_{u \in \R, z \in \C, |z|=1} \lambda_{\min}(\sum_{j=0}^\infty \Psi_j(u) z^j) \geq \gamma_{\inf}^{1/2}>0,
\end{eqnarray*}
 then $\inf_{u}\inf_{\omega}\lambda_{\min}(f(\omega;u))\geq \gamma_{\inf}$. Thus the conditions in Theorem~\ref{lemma:spectral} 
are satisfied, and for a sufficiently large $N_0$, there exists $0 < \lambda_{\inf} \leq \lambda_{\sup} <\infty $ such that for all $N\geq N_0$ we have 
\begin{eqnarray*}
0 < \lambda_{\inf} \leq \lambda_{\inf}({\bf C}^{(N)} ) \textrm{ and }
\lambda_{\sup}({\bf C}^{(N)} ) \leq \lambda_{\sup} <\infty.
\end{eqnarray*} 
\end{example}

In summary, the results in this section tell us the following. If an array of nonstationary time series satisfy Assumption 
\ref{assum:LS}, then there exists a stationary time series $\{X_{t}(u)\}$ whose covariance is $\{C_{r}(u)\}$. Conversely, if we define a nonstationary
time series $\{X_{t,N}\}_{t}$ with covariance $C^{(N)}$ and an accompanying stationary time series 
$\{X_{t}(u)\}_{t}$ whose covariances satisfy (\ref{eq:fomegaspectral})
and Assumption \ref{assum:LS}(ii,iii), 
then the positive definite condition in Assumption \ref{assum:LS}(i) holds. 
One important application of this result is 
given in Example \ref{exam:tvMA1}. However, the same result holds for
more general models, including the models which satisfy 
the physical dependence conditions considered in \cite{p:zhou-09}, \cite{p:dah-19},  
\cite{p:kar-21}, \cite{p:zha-21} and \cite{p:din-zho-21}). 
In the following theorem we
make this precise.

\begin{theorem}\label{theorem:physical}
Suppose that $\{X_{t,N}\}_{t}$ is a zero 
mean multivariate time series of dimension $p$ with the 
causal representation $X_{t,N} = G_{t,N}(\mathcal{F}_{t})$ where $\mathcal{F}_{t} =
(\varepsilon_{t},\varepsilon_{t-1},\ldots)$
and $\{\varepsilon_{t}\}$ are independent, identically distributed (iid)
random vectors of dimension $p$. Associated with $\{X_{t,N}\}$ we define the multivariate stationary
time series $\{X_{t}(u)\}_{t}$ where $X_{t}(u) = G(u,\mathcal{F}_t)$. 
Using $X_{t,N}$ and $X_{t}(u)$ we define the error process
\begin{eqnarray*}
e_{t,N} = X_{t,n} - X_{t}(t/N) = E_{t,N}(\mathcal{F}_{t}).
\end{eqnarray*}
and difference process ${X}_{t}^{v_1,v_2} = (X_{t}(v_1) - X_{t}(v_2))$.
Suppose $\{\widetilde{\varepsilon}_{t}\}_{t}$ are iid random vectors
that are independent of $\{\varepsilon_{t}\}$ but with the same
distribution and define
$\mathcal{F}_{t|\{t-j\}} =
(\varepsilon_{t},\varepsilon_{t-1},\ldots,
\varepsilon_{t-j+1},\widetilde{\varepsilon}_{t-j},\linebreak \varepsilon_{t-j-1},\ldots).$
Then
we define the coupled processes as $X_{t,N|\{t-j\}} = G_{t,n}(\mathcal{F}_{t|\{t-j\}})$,
$X_{t|\{t-j\}}(u) = G(u, \mathcal{F}_{t|\{t-j\}})$,
$X_{t|\{t-j\}}^{v_1,v_2} = G(v_1, \mathcal{F}_{t|\{t-j\}})-G(v_2,
\mathcal{F}_{t|\{t-j\}})$ and  $e_{t,n|\{t-j\}} =
E_{t,N}(\mathcal{F}_{t|\{t-j\}})$. Suppose the following hold:
\begin{itemize}
\item[(A)] Spectral-norm physical dependence 
\begin{eqnarray*}
&&\sup_{N}\sup_{t}\|\var(X_{t,N} - X_{t,N|\{t-j\}})\|_{2} \leq K\delta_{j} \\
&&\sup_{u} \|\var(X_{t|\{t-j\}}(u) - X_{t|\{t-j\}}(u))\|_{2}\leq
                                                                   K\delta_{j}
  \\
&&\sup_{t}\|\var({X}_{t}^{v_1,v_2}-X_{t|\{t-j\}}^{v_1,v_2})\|_{2} \leq K|v_1-v_2|\delta_{j} \\
&& \sup_{N}\sup_{t} \|\var(e_{t,N} - e_{t,N|\{t-j\}}) \|_{2} \leq
  KN^{-1}\delta_{j},
\end{eqnarray*}
where $\delta_{j} = \gu{j}^{-\kappa}$, $\kappa>3$ and $K$ is a finite
constant. 
\item[(B)] Let $C_{r}(u) = \cov[X_{0}(u),X_{r}(u)]$ and 
$f(\omega;u) = \sum_{r\in \mathbb{Z}}C_{r}(u)\exp(ir\omega)$. 
Then we assume the spectral density matrices satisfy 
$$0<\inf_{u,\omega}\lambda_{\inf}f(\omega;u)\leq 
\sup_{u,\omega}\lambda_{\sup}f(\omega;u)<\infty.$$
\end{itemize}
Under the above conditions, Assumption \ref{assum:LS}(ii,iii) is
satisfied (with the same $\kappa$ as that given in the conditions)
and for a sufficiently
large  $N_0$, Assumption \ref{assum:LS}(i) is satisfied.
\end{theorem}
\begin{proof}
In Appendix \ref{sec:proof2}.
\end{proof}

We observe that in the theorem above the physical dependence condition (A) is described in terms of a spectral-norm of a variance. 
\begin{remark}
It is worth mentioning that Condition (A) in Theorem \ref{theorem:physical} is equivalent to 
$$\|\var(X_{t,N} - X_{t,N|\{t-j\}})\|_{2}=(\max_{\|x\|_2=1}\Ex[ |x^\top(X_{t,N} - X_{t,N|\{t-j\}})|^2])^{1/2}.$$ Using the latter representation a generalisation to a bound on the $q$th moment:
\\*$(\max_{\|x\|_2=1} \Ex[ |x^\top(X_{t,N} - X_{t,N|\{t-j\}})|^q])^{1/q}$ is possible, thus generalising physical dependence in terms of any norm.  Recently, \cite{p:yi-22} used  such a generalisation. 
\end{remark}
\begin{example}[Locally stationary stochastic recurrence
  equations]\label{exam:sre}
  We now show that the nonstationary stochastic
  recurrence models studied in \cite{p:sub-06} and \cite{p:dah-19}
  satisfy the conditions in  Theorem \ref{theorem:physical}.
  
  Let us suppose that $\{X_{t,N}\}$ has the representation
  \begin{eqnarray*}
X_{t,N} = A(t/N,\varepsilon_{t})X_{t-1,N} + b(t/N,\varepsilon_{t})
  \end{eqnarray*}
  where $\{\varepsilon_{t}\}$ are iid random vectors.
The above model includes time-varying random coefficient vector
autogressive models, time-varying vector GARCH models and Bilinear
models (if the $\varepsilon_{t}$ in $A(t/N,\varepsilon_{t})$ were
changed to $\varepsilon_{t-1}$) as special cases. 

Based on the  above model we define the stationary time series model
 \begin{eqnarray*}
X_{t}(u) = A(u,\varepsilon_{t})X_{t-1}(u) + b(u,\varepsilon_{t}).
  \end{eqnarray*}
 Suppose $\sup_{u}\|\Ex[A(u,\varepsilon_{t})
A(u,\varepsilon_{t})^{\top}]\|_{2}<\rho<1$, $\sup_{u}\|\Ex[b(u,\varepsilon_{t})
b(u,\varepsilon_{t})^{\top}]\|_{2}<\infty$ and for all $v_1$ and $v_2$ 
$\|\Ex[(A(v_1,\varepsilon_{t})-A(v_2,\varepsilon_{t}))(A(v_1,\varepsilon_{t})-A(v_2,\varepsilon_{t}))^{\top}]\|_{2}\leq
K|v_1-v_2|$ and $\|\Ex[(b(v_1,\varepsilon_{t})-b(v_2,\varepsilon_{t}))(b(v_1,\varepsilon_{t})-b(v_2,\varepsilon_{t}))^{\top}]\|_{2}\leq
K|v_1-v_2|$. Under these conditions it can be shown that $X_{t,N}$ and
$X_{t}(u)$ almost surely have the causal solution
\begin{eqnarray*}
X_{t,N} &=& g_{t,N}(\mathcal{F}_{t}) = 
  \sum_{s=0}^{\infty}\prod_{i=0}^{s-1}A((t-i)/N,\varepsilon_{t-i})
  b((t-s)/N,\varepsilon_{t-s}) \\
X_{t}(u) &=&  g(u,\mathcal{F}_{t}) =\sum_{s=0}^{\infty}\prod_{i=0}^{s-1}A(u,\varepsilon_{t-i})
  b(u,\varepsilon_{t-s}). 
\end{eqnarray*}
Further in Lemma \ref{lemma:physicaldependencemarkov} we show that 
Condition (A) in Theorem \ref{theorem:physical}  holds with
$$\delta_{j} = K(\sum_{s=j}^{\infty}s^{3/2}\rho^{(s-1)/2})^{2}.$$
\end{example}

\subsection{Locally stationary approximations of the inverse covariance}\label{sec:LS1}

In this section we show that properties on 
the covariance operator $\BS C^{(N)}$ transfer to the inverse
covariance operator 
$\BS D^{(N)} = (\BS C^{(N)})^{-1}$. 
Specifically, in the following theorem we show that  the relationship between 
 $\BS C^{(N)}$ and $\BS C(u)$
 in Assumption \ref{assum:LS}(ii,iii) carry over to $\BS D^{(N)}$ and $\BS D(u)= \BS C(u)^{-1}$ up to a (small)  loss in rate. 
This result is used to show "approximate" smoothness of the time-varying VAR coefficients in representation (\ref{eq:ARinfty}). 

\begin{theorem}\label{thm.D.smooth}
Suppose Assumption~\ref{assum:LS} holds.  Then for all  $t,\tau \in \Z$,
$D_{t-\tau}(u)$ is Lipschitz, in the sense that for all $u,v\in \R$
\begin{eqnarray}
\|D_{t-\tau}(u)-D_{t-\tau}(v)\|_2\leq \mathcal{K} |u-v| \zeta(\tau-t)^{\kappa-1}. \label{lem.D.cont}
\end{eqnarray}
Furthermore, we have for all $t,\tau \in \Z$
\begin{eqnarray}
 \left\|\left[\BS D^{(N)}-\BS D(t/N) \right]_{t,\tau}\right\|_2
\leq \mathcal{K} \zeta(t-\tau)^{\kappa-2}\min(1/N,2\zeta(t-\tau)) \label{lem.Dn.smooth},
\end{eqnarray}
where $\mathcal{K}$ is a finite constant that is independent of $u,v,t,\tau$.
\end{theorem}
\begin{proof}
In Appendix \ref{sec:proof3}.
\end{proof}

An important consequence of Theorem~\ref{thm.D.smooth} is that when working with $\BS C$ and $\BS D$ it is enough to put smoothness conditions on one of them as the smoothness transfers to the other. In particular,  
conditions can be stated in terms of the covariance of the original time series.
 Furthermore, we note that differentiability conditions also transfer from $C_{r}(u)$ to $D_{r}(u)$. E.g., if one starts with the condition that for all $r$ $\sup_{u} \| \frac{d C_r(u)}{du}\|_2 \leq K\zeta(r)^{\kappa-1}$, then using the same arguments as those used in the proof of Theorem~\ref{thm.D.smooth} (outlined after the proof of Theorem~\ref{thm.D.smooth} in Appendix \ref{sec:proof3}) we have
\begin{align} 
\label{eq.claim.diff}
    \|\frac{d D_r(u)}{du}\|_2\leq \mathcal{K} \zeta(r)^{\kappa-1}.
\end{align}
Smoothness and differentiability conditions on 
$\BS D^{(N)}$ and $\BS D(u)$ are used in 
\cite{p:bas-sub-21} (stated in Assumption 4.2) to obtain certain rates of decay on the Fourier transform of $\BS D^{(N)}$.
Theorem \ref{thm.D.smooth} and (\ref{eq.claim.diff}) 
show that these conditions can be equivalently stated in terms of smoothness and differentiability conditions on covariance 
$\BS C^{(N)}$ and $\BS C(u)$.
It is worth noting that the loss in the rate of decay for the inverse in Section~\ref{sec:inverserate} is also present in Theorem \ref{thm.D.smooth}.

We now state a result that is analogous to Theorem
 \ref{thm.D.smooth}, but for one-sided matrices. This result will be
 useful in proving Theorem \ref{lemma:smoothVAR}  (below) on
 smoothness properties of  time-varying VAR representations. 
\begin{theorem} \label{cor.Dn.smooth}
Suppose Assumption \ref{assum:LS} holds and let
$\BS C^{(N)}(-\infty,T) = (C_{t,\tau}^{(N)};t,\tau \leq T)$
  and $\BS C(-\infty,T;u) = (C_{t,\tau}(u);t,\tau\leq T)$
Then for all $t,\tau \leq
T$ we have
\begin{eqnarray*}
 &&\|\left[\BS C^{(N)}(-\infty,T)^{-1}-\BS C(-\infty,T; T/N)^{-1}
     \right]_{t,\tau}\|_2 \\ 
&\leq& \mathcal{K}\zeta(t-\tau)^{\kappa-2}\min(1/N,2\zeta(t-\tau))
\end{eqnarray*}
\end{theorem}
\begin{proof}
In Appendix \ref{sec:proof3}
\end{proof}

We now apply 
Theorem~\ref{thm.D.smooth} to the popular 
time-varying VAR model. 
Let us suppose that $\{X_{t,N}\}$ has the tv-VAR$(d)$ representation 
\begin{eqnarray}
\label{eq:tvVAR}
X_{t,N} = \sum_{j=1}^{d}\Phi_j(t/N)X_{t-j,N} + \Sigma(t/N)^{1/2}\varepsilon_{t}, \qquad
  t\in \mathbb{Z},
\end{eqnarray}
where $\{\varepsilon_{t}\}_{t}$ are uncorrelated random vectors with  variance $I_p$. 
In contrast to the tv-VAR representation given in \eqref{eq:ARfinite}, the tv-VAR model is defined with Lipschitz conditions
on the matrices $\Phi_{j}(\cdot)$ and $\Sigma(\cdot)$. The tv-VAR$(d)$ model with smooth AR coefficients as defined in (\ref{eq:tvVAR}) is attractive because its coefficients are straightforward to interpret and
has been used in econometrics and 
in neuroscience (see, for example, \cite{ding2017sparse,Shojaie_tvVAR,yan2021time}).
Let $\BS C^{(N)}$ denote the covariance corresponding to $\{X_{t,N}\}$. 
Obtaining a rate of decay for the covariance  by directly analyzing
$\BS C^{(N)}$ is unwieldy (see \cite{kunsch1995note} for the univariate proof). However, we show below that starting with the inverse $\BS D^{(N)} = 
(\BS C^{(N)})^{-1}$ 
(which is a banded matrix, since $X_{t,N}$ has a tvVAR$(p)$ representation) we can use 
Theorem~\ref{thm.1} and \ref{thm.D.smooth}, to transfer the  information on the rate of decay of the inverse covariance operator to the covariance operator itself. 

\begin{corollary}[Application of Theorem \ref{thm.D.smooth} to tvVAR models] \label{cor.var}
Suppose that the multivariate time series $\{X_{t,N}\}_{t}$ has the
time-varying VAR$(d)$ representation in (\ref{eq:tvVAR}), where we assume there exists a 
$\delta>0$ and $\gamma$
where 
\begin{eqnarray}
\label{eq:roots}
\inf_{u \in \R, z \in \C, |z|\leq 1+\delta} \lambda_{\min}(I_{p}-\sum_{j=1}^{d} \Phi_j(u) z^j) \geq \gamma_{}>0,
\end{eqnarray}
and the
matrices $\Phi_{j}(\cdot)$ are Lipschitz continuous in the sense that
$\|\Phi_{j}(u)-\Phi_{j}(v)\|_{2} \leq K|u-v|$. We further assume that $\Sigma(\cdot)$
is Lipschitz continuous in the sense that $\|\Sigma(u)-\Sigma(v)\|_{2} \leq K|u-v|$ and for all $u\in \R$ $\Sigma(u)$ is positive definite (with eigenvalues that are bounded from above and away from zero uniformly over all $u$). Let $\BS C^{(N)}$ denote the covariance operator of $\{X_{t,N}\}_{t}$ and 
$
C_{r}(u) = \int_{0}^{2\pi}f(\omega;u)\exp(-ir\omega)d\omega,
$
where $f(\omega;u)=[I_p-\sum_{j=1}^{d}\Phi_{j}(u)\exp(-ij\omega)]^{-1}\Sigma(u)
([I_p-\sum_{j=1}^{d}\Phi_{j}(u)\exp(ij\omega)]^{-1})^{\top}$. Then, there exists an $N_0$ and $0<\rho<1$ such that for all $N>N_0$ we have 
$
\|C_{t,\tau}^{(N)}\|_{2} \leq \mathcal{K}\rho^{|t-\tau|},
$
$\|C_{r}(u)-C_{r}(v)\|_2\leq \mathcal{K}|u-v|\rho^{|t-\tau|}$, and  
$
\|C_{t,\tau}^{(N)}-C_{t-\tau}(t/N)\|_{2} \leq \mathcal{K}\rho^{|t-\tau|}/N.$
\end{corollary}
\begin{proof}
In Appendix \ref{sec:proof3}.
\end{proof}

\begin{remark}[Differentiability of the tv-VAR covariance]
As mentioned after Theorem~\ref{thm.D.smooth} in equation (\ref{eq.claim.diff}),
smoothness conditions in terms of differentiability transfer between
$\BS C(\cdot)$ and $\BS D(\cdot)$. For the tv-VAR model
this implies that smoothness conditions formulated in terms of
differentiability of the transition matrices $\Phi_j(\cdot)$ 
transfer to $\BS D(\cdot)$, then from equation (\ref{eq.claim.diff})
they transfer to $\BS C(\cdot)$. 
\cite{ding2017sparse}, Lemma 3.1 
also prove that differentiability of $\Phi_{1}(\cdot)$ implies differentiability of the covariance for tv-VAR$(1)$ models. They show this result by directly connecting the covariance to $\Phi_{1}(\cdot)$ through 
the tv-VAR$(1)$ model. However, their proof requires the additional condition that $\|\Phi_1 \|_1=\max_{\|x\|_1=1}\|\Phi_1 x\|_1<1$, which places quite strict conditions on the VAR parameters.
\end{remark}

We apply Corollary \ref{cor.var} to the time-varying ARCH process. 
\begin{example}[The time-varying ARCH$(p)$ process]
The time-varying ARCH model is defined as follows. Let 
\begin{eqnarray*}
X_{t,N} = \sigma_{t,N}Z_{t} \qquad \sigma_{t,N}^{2} = a_{0}(t/N)+\sum_{j=1}^{p}a_{j}(t/N)X_{t-j}^{2},
\end{eqnarray*}
where $\{Z_{t}\}$ are iid random variables with mean zero and variance
one and the coefficients $a_{j}(\cdot)$ are Lipschitz continuous and such that
$\inf_{u}a_{0}(u)>0$, for $1\leq j \leq p$ $a_{j}(\cdot)\geq 0$
and are such that
$\sup_{u}(\Ex[Z_{0}^{4}])^{1/2}\sum_{j=1}^{p}a_{j}(u)<1$. Under these
conditions it can be shown that $\sup_{t,N}\Ex[X_{t,N}^{4}]<\infty$
and that $X_{t,N}^{2}$ has the tvAR$(p)$ representation 
\begin{eqnarray*}
X_{t,N}^{2} = a_{0}(t/N)+\sum_{j=1}^{p}a_{j}(t/N)X_{t-j,N}^{2}+\varepsilon_{t,N}
\end{eqnarray*}
where $\varepsilon_{t,N} = \sigma^{2}_{t,N}(Z_{t}^{2}-1)$. The
condition that $\sup_{u}\sum_{j=1}^{p}a_{j}(u)<1$ (which is implied by the condition
$\sup_{u}(\Ex[Z_{0}^{4}])^{1/2}\sum_{j=1}^{p}a_{j}(u)<1$) implies that
(\ref{eq:roots}) holds. Thus by applying Corollary \ref{cor.var},
Assumption \ref{assum:LS} holds for the time-varying ARCH process.
\end{example}

We have shown in (\ref{eq:ARinfty}) that under certain 
conditions all nonstationary time series have an AR$(\infty)$ representation.
But there is no guarantee that the AR parameters are smooth.
Below, we show below that under the locally stationary conditions in Assumption~\ref{assum:LS} a smooth \emph{approximation} is possible.

We recall from (\ref{eq:ARinfty}) that $\{X_{T,N}\}_{t}$ has the
representation 
\begin{eqnarray}
\label{eq:VARinfty1}
X_{T,N}^{} = \sum_{j=1}^{\infty}\Phi_{T,j}^{(N)}X_{T-j,N}^{} +\varepsilon_{T,N}, 
\end{eqnarray}
where $\{\varepsilon_{T,N}\}_{t\in \mathbb{Z}}$ are uncorrelated random vectors with $\Sigma_{T,N}=\var[\varepsilon_{T,N}]$. We have shown in Section \ref{sec:LS0} that
under Assumption \ref{assum:LS} there exists a stationary time series $\{X_{t}(u)\}_{t}$ with autocovariance $\{C_{r}(u)\}_{r}$. Using the arguments leading to \eqref{eq:ARinfty}, it can be shown that $\{X_{t}(u)\}_{t}$ has the VAR$(\infty)$ representation
\begin{eqnarray}
\label{eq:VARinfty2}
X_{t}(u) = \sum_{j=1}^{\infty}\Phi_{j}(u)X_{t-j}(u) +\varepsilon_{t}(u),
\end{eqnarray}
where $\varepsilon_{t}(u)$ are uncorrelated random vectors with variance $\Sigma(u) = \var[\varepsilon_{t}(u)]$. 
In the following theorem we show that 
$\{\Phi_{T,j}^{(N)}\}$ can be approximated by the stationary VAR coefficients $\{\Phi_{j}(u)\}$.

\begin{theorem}\label{lemma:smoothVAR}
Suppose the array of time series $\{X_{T,N}\}_{t}$ satisfy
Assumption \ref{assum:LS} and let $\{\Phi_{t,j}^{(N)}\}_{j}$ be defined as in (\ref{eq:VARinfty1}) with $\Sigma_{T}^{(N)} = \var[\varepsilon_{T,N}]$. Additionally, let $\{X_{t}(u)\}_{t}$ be the locally stationary
approximation defined in (\ref{eq:VARinfty2}).  
\begin{itemize}
\item[(i)] Then for all $T\in \mathbb{Z}$ and $j\geq 1$ we have
\begin{eqnarray*}
\label{eq:smooth1}
\|\Sigma_{T}^{(N)} -  \Sigma(T/N) \|_2 &\leq& \frac{\mathcal{K}}{N} \nonumber\\
  \textrm{ and }
\| \Phi_{T,j}^{(N)} - \Phi_{j}(T/N)\|_2 &\leq& \mathcal{K} \zeta(j)^{\kappa-2} \min(2\zeta(j),1/N)
\end{eqnarray*}
\item[(ii)] For all $u_1,u_2\in \mathbb{R}$ and $j\geq 1$
\begin{eqnarray*}
\label{eq:smooth2}
\|\Sigma(u_1) -  \Sigma(u_2) \|_2 &\leq& \mathcal{K}|u_1-u_2| \\
  \textrm{ and }
\| \Phi_{j}(u_1) - \Phi_{j}(u_2)\|_2 &\leq& \mathcal{K}\zeta(j)^{\kappa-1}|u_1-u_2|.
\end{eqnarray*}
\end{itemize}
\end{theorem}
\begin{proof}
In Appendix \ref{sec:proof3}
\end{proof}

\begin{remark}[Approximation and estimation by finite order tv-VAR]
From above theorem, if a process is locally stationary then
it can be approximated by an time-varying VAR$(\infty)$ time series
with slowly varying parameters. Consequently by using \cite{meyer2015vector} this
implies that the locally stationary time series can be approximated
with a finite order time-varying VAR$(d)$ with slowly varying
parameters. More precisely, let $\{\Phi_{d,j}(u)\}_{j=1}^{d}$ denote
the finite order (stationary) VAR$(d)$ parameters associated with 
the vector autocovariance $\{C_{r}(u)\}_{r}$. Then by 
\cite{meyer2015vector} we have 
$\sup_{u}\sum_{j=1}^{d}\|\Phi_{d,j}(u)-\Phi_{j}(u)\|_{2}\leq
\mathcal{K}d^{-\kappa}$. This result together with Theorems \ref{lemma:onesided},
\ref{lemma:smoothVAR} and the triangle inequality gives
\begin{eqnarray*}
\sum_{j=1}^{d}\|\Phi_{T,d,j}-\Phi_{d,j}(T/N)\|_{2} \leq \mathcal{K}\left(N^{-1}+\zeta(d)^{\kappa-3/2}\right).
\end{eqnarray*}
A potential application of Theorem \ref{lemma:smoothVAR} is that it
could be used (i) in forecasting and (ii) to
 develop a bootstrap procedure for nonstationary time
series by transferring the widely used stationary AR-sieve to the locally
stationary setup. Both procedures would require estimators of the the
finite order time-varying VAR parameters $\{\Phi_{d,j}(T/N)\}$. One
could estimate this using local kernel methods or the sieve
estimation method described in \cite{p:din-zho-19}.
\end{remark}

\begin{remark}[Innovations and Kolomogorov's formula]
An immediate implication of the above
 result is that the time varying innovation variance 
 $\Sigma_{t}^{(N)}$ can be approximated by Kolomogorov's formula 
\begin{eqnarray*}
\textrm{det}[\Sigma_{t}^{(N)} ]  = \int_{-\pi}^{\pi}\log
  \textrm{det}[f(t/N;\omega) ]d\omega + O(1/N)
\end{eqnarray*}
where $f(u;\omega) = \sum_{r\in \mathbb{Z}}C_{r}(u)\exp(ir\omega)$.
A similar result was obtained in
\cite{p:liu-tan-omb-21}, Proposition~1 for a specific
class of locally stationary time series. 
\end{remark}

\subsection{The partial covariance of a locally stationary time series}\label{sec:LS2}

The partial covariance is commonly used in the analysis of time series as a measure of
linear dependence between two time series after accounting for all the
other components in the time series. For stationary time series, the analysis is 
typically conducted through the partial spectral coherence which 
is the standardized Fourier transform of the partial covariance,  and
is, conveniently, a function of the spectral density matrix function (cf. \cite{b:pri-81,b:bri-01,p:dah-00b,p:kra-pap-21}).
For nonstationary time series the time-varying partial spectral coherence can be 
defined as a function of the localized inverse spectral density, as was done in \cite{p:par-14}.
However, as far as we are aware, there are 
no results that connect this definition 
(of the time-varying partial spectral coherence) 
to the actual partial covariance of the underlying nonstationary time series. 

We use the results on inverse covariances (developed in Section \ref{sec:LS1}) 
to show that the partial covariance of a locally stationary time series 
(as defined in Assumption \ref{assum:LS})
can be approximated by a smooth function, which, in turn, is the 
partial covariance of the locally stationary approximation $\{X_{t}(u)\}_{t}$.
We show below that this result can be used to justify using the time-varying partial spectral coherence as an approximation of the Fourier transform of the localized 
partial covariance. 

We start by defining the partial covariance for nonstationary time series. For this, let $\mathcal{H}^{(N)} = \overline{\textrm{sp}}(X_{t,N}^{(c)};t\in \mathbb{Z},1\leq
c\leq p)$ denote the space spanned by the entire multivariate time series. Furthermore, let $\mathcal{S} \subseteq \{1,\dots,p\}=:V$ be a set of indices referring to components of the time series and $\mathcal{H}^{(N)} - (X^{(c)};c\in \mathcal{S}) =
\overline{\textrm{sp}}[X_{s,N}^{(c)};s\in \mathbb{Z},c\in
\mathcal{S}^{\prime}]$ be the space spanned by the entire time series  of the components in $\mathcal{S}^\prime$ only, where $\mathcal{S}^{\prime}$ denotes the complement of $\mathcal{S}$.
Let $P_{\mathcal{M}}(Y)$ denote the
orthogonal projection of $Y\in \mathcal{H}^{(N)}$ onto the subspace
$\mathcal{M}$. For any $\mathcal{S} \subseteq V$, we define the residual of
$X_{t,N}^{(a)}$ after projecting on 
$ \mathcal{H}^{(N)} -
  (X^{(c)};c\in \mathcal{S})$  as 
\begin{eqnarray}
\label{eq:XtNcond}
X_{t,N}^{(a)|\shortminus\mathcal{S}} &:=& X_{t,N}^{(a)} - P_{\mathcal{H}^{(N)} -
  (X^{(c)};c\in \mathcal{S})}(X_{t,N}^{(a)}), t\in \mathbb{Z}. 
\end{eqnarray}
In the definitions below we focus on the two sets
$\mathcal{S}=\{a,b\}$ and $\mathcal{S}=\{a\}, a,b \in V, a \not=b.$
Using the above, we define the  partial covariance
\begin{eqnarray}
\label{eq:partialTS1}
\Delta_{t,\tau,N}^{\shortminus\{a,b\}} = \left(
\begin{array}{cc}
\rho_{t,\tau,N}^{(a,a)|\shortminus\{a,b\}} & \rho_{t,\tau,N}^{(a,b)|\shortminus\{a,b\}} \\
\rho_{t,\tau,N}^{(b,a)|\shortminus\{a,b\}} &\rho_{t,\tau,N}^{(b,b)|\shortminus\{a,b\}}\\
\end{array}
\right) :=
\cov\left[ \left(
\begin{array}{c}
X_{t,N}^{(a)|\shortminus\{a,b\}} \\
X_{t,N}^{(b)|\shortminus\{a,b\}} \\
\end{array}
 \right),
\left(
\begin{array}{c}
X_{\tau,N}^{(a)|\shortminus\{a,b\}} \\
X_{\tau,N}^{(b)|\shortminus\{a,b\}} \\
\end{array}
 \right)
\right]
\end{eqnarray}
and self partial covariance
\begin{eqnarray}
\label{eq:partialTS2}
\rho_{t,\tau,N}^{(a,a)|\shortminus\{a\}} = \cov[X_{t,N}^{(a)|\shortminus\{a\}}, X_{\tau,N}^{(a)|\shortminus\{a\}}].
\end{eqnarray}
As will become clear in the proof of the following theorem
$\Delta_{t,\tau,N}^{\shortminus\{a,b\}}$ and
$\rho_{t,\tau,N}^{(a,a)|\shortminus\{a\}}$ can be expressed in terms
of the matrix operator ${\BS C}^{(N)}$ and its inverse. Under Assumption \ref{assum:LS}
and by Theorem \ref{lem.Cu.pos.def} there exists a stationary time
series $\{X_{t}(u)\}_{t}$ which has covariance ${\BS C}(u)$, that
locally approximates ${\BS C}^{(N)}$. Using ${\BS C}(u)$ we will define the
partial covariances corresponding to the stationary time series
$\{X_{t}(u)\}_{t}$. In the theorem below we show that the partial
covariances of $\{X_{t}(u)=(X_{t}^{(1)}(u),\ldots,X_{t}^{(p)}(u))^{\top}\}_{t}$ locally approximates the partial
covariance of $\{X_{t,N} = (X_{t,N}^{(1)},\ldots,X_{t,N}^{(p)})^{\top}\}_{t}$. To do this, analogous to (\ref{eq:XtNcond}),
(\ref{eq:partialTS1}) and (\ref{eq:partialTS2}) we define 
\begin{eqnarray}
\label{eq:Xucond}
X_{t}^{(a)|\shortminus\mathcal{S}}(u) &:=& X_{t,N}^{(a)}(u) - P_{\mathcal{H}_u -
  (X^{(c)}_u;c\in \mathcal{S})}(X_{t}^{(a)}(u))  \textrm{ for }t\in \mathbb{Z},
\end{eqnarray}
\begin{eqnarray}
\label{eq:partialTS1u}
\Delta_{t-\tau}^{\shortminus\{a,b\}}(u) = \left(
\begin{array}{cc}
\rho_{u,t-\tau}^{(a,a)|\shortminus\{a,b\}} & \rho_{u,t-\tau}^{(a,b)|\shortminus\{a,b\}} \\
\rho_{u,t-\tau}^{(b,a)|\shortminus\{a,b\}} &\rho_{u,t-\tau}^{(b,b)|\shortminus\{a,b\}}\\
\end{array}
\right) :=
\cov\left[ \left(
\begin{array}{c}
X_{t}^{(a)|\shortminus\{a,b\}}(u) \\
X_{t}^{(b)|\shortminus\{a,b\}}(u) \\
\end{array}
 \right),
\left(
\begin{array}{c}
X_{\tau}^{(a)|\shortminus\{a,b\}}(u) \\
X_{\tau}^{(b)|\shortminus\{a,b\}}(u) \\
\end{array}
 \right)
\right]
\end{eqnarray}
and self partial covariance
\begin{eqnarray}
\label{eq:partialTS2u}
\rho_{t-\tau}^{(a,a)|\shortminus\{a\}}(u) = \cov[X_{t}^{(a)|\shortminus\{a\}}(u), X_{\tau}^{(a)|\shortminus\{a\}}(u)].
\end{eqnarray}
We note that a key ingredient in the proof of the theorem below is that
the partial covariance can be expressed as 
\begin{eqnarray*}
\var\left[X_{t,N}^{(e)|\shortminus\{a,b\}}; t\in \mathbb{Z}, e\in\{a,b\}\right] = 
\BS C_{\mathcal{S},\mathcal{S}}  - 
\BS C_{\mathcal{S},\mathcal{S}^{\prime}}
\BS C_{\mathcal{S}^{\prime},\mathcal{S}^{\prime}}^{-1}
\BS C_{\mathcal{S},\mathcal{S}^{\prime}}^{\top}, 
\end{eqnarray*}
where $\mathcal{S} = \{a,b\}$,  
$\BS C_{\mathcal{S},\mathcal{S}} = 
(\BS C^{(e,f)};e,f\in \mathcal{S})$
(similarly for $\BS C_{\mathcal{S},\mathcal{S}^{\prime}}$ and $\BS C_{\mathcal{S}^{\prime},\mathcal{S}^{\prime}}$) and 
$\BS C^{(e,f)} = (\cov[X_{t,N}^{(e)},X_{\tau,N}^{(f)}];t,\tau \in \mathbb{Z})$. The presence of $\BS C_{\mathcal{S}^{\prime},\mathcal{S}^{\prime}}^{-1}$ in the above expression explains why the results in the previous sections (in particular Theorem \ref{thm.D.smooth}) are necessary for proving the result. 

\begin{theorem}\label{thm.part}
Suppose Assumption~\ref{assum:LS} holds and let further 
$\Delta_{t,\tau,N}^{\shortminus\{a,b\}}$,
$\rho_{t,\tau,N}^{(a,a)|\shortminus\{a\}}$,
$\Delta_{t-\tau}^{\shortminus\{a,b\}}(u)$, and $\rho_{t-\tau}^{(a,a)|\shortminus\{a\}}(u)$
be defined as in (\ref{eq:partialTS1}), (\ref{eq:partialTS2}),
(\ref{eq:partialTS1u}) and (\ref{eq:partialTS2u}). Then 
 for all $a,b \in \{1,\dots,p\}$
\begin{eqnarray}
\label{eq:partial1}
\|\Delta_{t,\tau,N}^{\shortminus\{a,b\}} - \Delta_{t-\tau}^{\shortminus\{a,b\}}(t/N)\|_2 
&\leq& \mathcal{K} \zeta(t-\tau)^{\kappa-2}\min(1/N,\zeta(t-\tau)) \\
\label{eq:partial2}
\|\Delta_{t-\tau}^{\shortminus\{a,b\}}(u)-\Delta_{t-\tau}^{\shortminus\{a,b\}}(v)\|_2&\leq&
\mathcal{K} |u-v| \zeta(t-\tau)^{\kappa-1}
 \\
\label{eq:partial3}
\|\rho_{t,\tau,N}^{(a,a)|\shortminus\{a\}} - \rho_{t-\tau}^{(a,a)|\shortminus\{a\}}(t/N)\|_2 
&\leq&  \mathcal{K} \zeta(t-\tau)^{\kappa-2}\min(1/N,\zeta(t-\tau))  \\
\textrm{ and } 
\label{eq:partial4}
\|\rho_{t-\tau}^{(a,a)|\shortminus\{a\}}(u)-\rho_{t-\tau}^{(a,a)|\shortminus\{a\}}(v)\|_2
&\leq& \mathcal{K} |u-v| \zeta(t-\tau)^{\kappa-1},
\end{eqnarray}
where $\mathcal{K}$ is a positive generic constant.
\end{theorem}
\begin{proof}
In Appendix \ref{sec:proof3}.
\end{proof}
The above result provides the tools to prove the following. 
Let $\{X_{t,N}\}_{t}$ be an array of nonstationary time series that satisfy Assumption \ref{assum:LS} and 
$\{C_{r}(u)\}_{r}$ the corresponding stationary approximation covariance.
Let $f(\omega;u) = 
\sum_{r\in \mathbb{Z}}C_{r}(u)e^{ir\omega}$ and 
$\Gamma(\omega;u) = f(\omega;u)^{-1}$. 
Using the stationary partial spectral coherence  (see \cite{b:pri-81}, Section 9.3 and \cite{p:dah-00b}), the
 localized (complex) partial 
spectral coherence is defined as
\begin{eqnarray*}
g_{a,b}(\omega;u) = -\frac{\Gamma^{(a,b)}(\omega;t/N)}
 {(\Gamma^{(a,a)}(\omega;t/N)\Gamma^{(b,b)}(\omega;t/N))^{1/2}},
\end{eqnarray*}
where $\Gamma^{(a,b)}(\omega)$
denotes the $(a,b)$ entry of the matrix $\Gamma(\omega;u)$.
Under Assumption \ref{assum:LS} (for 
$\kappa >3$) and 
by using Theorem \ref{thm.part} it can be shown that
\begin{eqnarray*}
\frac{\sum_{r\in \mathbb{Z}}\rho_{t,t+r,N}^{(a,b)|\shortminus\{a,b\}}\exp(ir\omega)}{\sqrt{\sum_{r\in \mathbb{Z}}\rho_{t,t+r,N}^{(a,a)|\shortminus\{a,b\}}\exp(ir\omega)\sum_{r\in \mathbb{Z}}\rho_{t,t+r,N}^{(b,b)|\shortminus\{a,b\}}\exp(ir\omega)}}  = g_{a,b}\left(\omega;t/N\right)  + O(N^{-1}).
\end{eqnarray*}
In other words, the estimated  local partial spectral coherence
(based on an estimator of the local spectral density function) is
an estimator of the Fourier transform of the partial covariances of the nonstationary time series localised about time point $t$. 
This justifies using local spectral density estimation approaches for 
estimating the partial covariance.

\subsection*{Acknowledgements}
The authors would like to thank two anonymous referees for their
comments and suggestions, which greatly improved all aspects of the
paper. 
JK was supported by the Research Center (SFB) 884 ``Political Economy
of Reforms''(Project B6), 
funded by the German Research Foundation (DFG), and acknowledges the partial support of DFG 
(travel grant 493207657) and National Institute of Health (grants R01GM135926 and R21NS120227). SSR research was partially supported by the National
Science Foundation (grant DMS-1812128 and DMS-2210726).

\bibliographystyle{apalike}
\bibliography{bib_tsreg}

\appendix

\section{Supplementary material} \label{section.proofs}

\subsection{Summary of results in supplementary material}

In order to navigate the appendix we summarize below the contents and
main results in the appendix. 
\begin{itemize}
\item Appendix A gives all the background lemmas. 

In Appendix \ref{appendix:back} we state all the block operator
identities that are required in this paper. In Appendix
\ref{appendix:matrix} we state and derive several matrix norm
inequalities. This includes a Cauchy-Schwarz type bound for the
spectral norm  of cross covariance 
matrices (Lemma \ref{lemma:CauchySchwarzMV}). 
\item In Appendix B we prove the results in Section 2. A fundamental
  result required in the proof is  Lemma \ref{lem.demko} which gives a bound on the
entries for the inverse of block banded matrices.
\item In Appendix C we prove the results for Section 3. 

The proofs for Theorems \ref{lem.Cu.pos.def} and \ref{lemma:spectral}
are given in Appendix \ref{sec:proof2}. In Appendix \ref{sec:physical}
we consider models which satisfy the physical dependence conditions
first proposed in \cite{p:wu-05}. 
The proof of Theorems \ref{thm.D.smooth}, \ref{cor.Dn.smooth} and
\ref{thm.part} follow a similar set of arguments and are given in
Appendix \ref{sec:proof3}.
\end{itemize}

\subsection{Notation and background}\label{appendix:back}

Before proceeding with the proofs, we need to introduce some
notation. We define below unit vectors of appropriate dimension to
select sub-matrices or elements from the operator  $\BS A : \ell_{2,p}
\to \ell_{2,p}$. That is, $A_{s_1,s_2}=(e_{s_1} \otimes I_p)^\top \BS
A (e_{s_2} \otimes I_p)$, where $\otimes$ is the Kronecker product and 
$I_p$ denotes the identity operator in $\R^p$. Furthermore,
$A_{s_1,s_2}^{(a,b)}=(e_{s_1} \otimes e_a)^\top 
\BS A (e_{s_2} \otimes e_b)$ and we introduce the short notation for
this unit vector as $e_{(a,s)}=(e_s \otimes e_a)$. 

In the proofs below we will often consider sub-matrices, where one
column or row has been 
removed. To set-up the matrix notation for this, 
let $I$ denote the identity operator in $\ell_2$ and $I_{-k}$ the identity operator after removing the $k$th row, i.e., for $u \in \ell_2$, $I_{-k} u=(\dots,u_{-1},u_{0},u_{1},\dots,u_{k-1},u_{k+1},\dots)$. The same notation is used for 
operators in $\R^p$ and similar spaces. This results in the following operations applied to an operator $\BS A : \ell_{2,p}
\to \ell_{2,p}$:
\begin{itemize}
\item $(e_{s_1} \otimes e_a)^\top 
\BS A (e_{s_2} \otimes e_b)=A_{s_1,s_2}^{(a,b)}$. 
    \item $(I_{-k} \otimes I_p) \BS A$ removes from the infinite dimensional matrix $p$ rows (of infinite length) so that $A_{k,i}^{(a,b)}$ is removed for all $i \in \Z, a,b \in \{1,\dots,p\}$. 
    \item $\BS A (I_{-k} \otimes I_p)^\top$ removes from the infinite dimensional matrix $p$ columns (of infinite length) so that $A_{i,k}^{(a,b)}$ is removed for all $i \in \Z, a,b \in \{1,\dots,p\}$.
    \item $(I_{-k} \otimes I_p) \BS A (I_{-k} \otimes I_p)^\top=:\tilde {\BS A}$ is an infinite dimensional matrix where $A_{i,k}^{(a,b)}$ and $A_{k,i}^{(a,b)}$ are removed for all $i \in \Z, a,b \in \{1,\dots,p\}$. 
    \item $(I_{-k} \otimes I_p)^\top  \tilde{\BS A} (I_{-k} \otimes
      I_p)=\BS B$ is an infinite dimensional matrix where $p$ zero
      columns and rows (of infinite length) are added so that
      $B_{i,k}^{(a,b)}=0, B_{k,i}^{(a,b)}=0$ for all $i \in \Z, a,b
      \in \{1,\dots,p\}$. Additionally, for all $s_{1}\neq k$ and
      $s_{2}\neq k$ we have
    $((I_{-k} \otimes I_p)^\top  \tilde{\BS A} (I_{-k} \otimes I_p))_{s_1,s_2}=(\BS A)_{s_1,s_2}.
    $
    \item $(I_{-k} \otimes I_{-a}) \BS A$ removes from the infinite
      dimensional matrix $p-1$ rows (of infinite length) so that
      $A_{k,i}^{(c,b)}$ is removed for all $i \in \Z, a,b \in
      \{1,\dots,p\}, c\not=a$.
      \\Similarly, for the other operations used above.\\
    We denote $(I_{-k} \otimes I_{-a})=:I_{-(a,k)}$. 
    \item We have that $(I_{-k} \otimes I_p)^\top (I_{-k} \otimes
      I_p)$ is the identity operator
      on the reduced space and 
$(I_{-k} \otimes I_p)(I_{-k} \otimes I_p)^\top+(e_{k} \otimes
I_p)(e_{k} \otimes I_p)^\top=\BS I=(I \otimes I_p)$, where $\BS I$ is
the identify on the full space. Furthermore, $(e_{k} \otimes I_p)^\top
(I_{-k} \otimes I_p)=0$.
\item For $x \in \ell_q, q \in [1,\infty]$ we define $\|x\|_q=(\sum_{l \in Z} x_l^q)^{1/q}$ and $\|x\|_\infty=\max_{l \in \Z} |x_l|$.
For an operator $\BS B : \ell_2 \to \ell_2$, we also define the $\ell_q$-induced norms, that is for $q \in [1,\infty]$ we set $\|\BS B\|_q=:\sup_{\|x\|_q=1,x\in \ell_2} \|\BS B x\|_q,$ where $\|\BS B\|_\infty=\sup_{s_1 \in \Z} \sum_{s_2\in \Z} |B_{s_1,s_2}|$.
\end{itemize}


An important tool in the proofs is the inversion 
and manipulation of infinite dimensional (block) matrices. Under
certain conditions on both the matrices and the spaces we can treat
these in much the same way as finite dimensional matrices. An identity
that we will make frequent use of is the analogous version of the
block inversion identity but for infinite dimensional
operators. Suppose that $\BS U:(S_{1},S_{2})\rightarrow (S_{1},S_{2})$
where $S_{1}$ and $S_{2}$ are two Hilbert spaces and 
\begin{eqnarray*}
  \BS U =
  \left(
  \begin{array}{cc}
   \BS A & \BS B \\
  \BS  C &\BS D \\
   \end{array}
  \right).
\end{eqnarray*}
If the eigenvalues of $\BS U$ are bounded away from zero and from
infinite, 
then using equation (1.7.4) in \cite{b:tre-08}, page 43 (setting
$\lambda=0$) for the inversion of block operator matrices we have
\begin{eqnarray}
\label{eq:inverseblock1}
 \BS U^{-1} =
  \left(
  \begin{array}{cc}
    \widetilde{\BS A} & \widetilde{\BS B} \\
    \widetilde{\BS C} & \widetilde{\BS D} \\
   \end{array}
  \right) &=&
   \left(
  \begin{array}{cc}
    \widetilde{\BS A} \quad& - \widetilde{\BS A}\BS B\BS D^{-1}\\
 -\BS D^{-1}\BS C\widetilde{\BS A}\quad  & \BS D^{-1}+\BS D^{-1}\BS
                                      C\widetilde{\BS A}\BS B\BS D^{-1} \\
   \end{array}
  \right) 
\end{eqnarray}
where from Definition 1.6.1 in  \cite{b:tre-08}, page 35 we have
\begin{eqnarray}
\label{eq:schur1}  
  \widetilde{\BS A}  = (\BS A-\BS B\BS D^{-1}\BS C)^{-1} \textrm{ and
  }
   \widetilde{\BS D}  = (\BS D-\BS C\BS A^{-1}\BS B)^{-1}.
\end{eqnarray}
An immediately consequence of the above is that the difference in the
block diagonal entries is
\begin{eqnarray}
\label{eq:schur2}  
  \BS A- \widetilde{\BS A}^{-1} = \BS B\BS D^{-1}\BS C = \BS B(
\tilde{\BS D}-\tilde{\BS C}\tilde{\BS A}^{-1}\tilde{\BS B})\BS C.
\end{eqnarray}
We will make frequent use of (\ref{eq:inverseblock1}) and
(\ref{eq:schur2}) in the proofs.


\subsection{Some background results}\label{appendix:matrix}

\begin{lemma} \label{lemma:operatornorm}
Suppose that $\{A_{\ell}\}_{\ell=1}^{\infty}$ is a sequence of
$p\times p$ dimensional matrices and $\sum_{\ell=1}^{\infty}\|A_{\ell}\|_{2}^{2}<\infty$.
Define the sequence space $\ell_{2,p}^{+} = \{w=(v_{1},v_{2},\ldots) :
v_j \in \R^p\}$ and the linear operator $\BS A = (A_{\ell};\ell \geq
0)$, where $\BS A:\ell_{2,p}^{+}\rightarrow \ell_{2,p}^{+}$.
Then 
\begin{eqnarray*}
\|\BS A\|_{2} \leq (\sum_{\ell=1}^{\infty} \| A_{\ell}\|_2^2)^{1/2}
\end{eqnarray*}
\end{lemma}
\begin{proof}
Let $x=(x_{1},x_{2},\ldots)$ where $x_{l}\in \mathbb{R}^{p}$. By
definition of the $\|\cdot\|_2$ operator norm we have
\begin{align*}
    \|\BS A\|_2 & = \sup_{\|x\|_2=1, x \in \ell_{2,p,1}} x^\top \BS
                  A^\top \BS A x = \sup_{\|x\|_2=1, x \in
                  \ell_{2,p,1}} (\sum_{l_1,l_2=1}^\infty x_{l_1}^\top A_{l_1}^\top A_{l_2} x_{l_2})^{1/2} \\
    &\leq \sup_{\|x\|_2=1, x \in \ell_{2,p,1}} \sum_{l=1}^\infty
      \|x_l\|_2 \|A_{l}\|_2 \\
    &\leq \sup_{\|x\|_2=1, x \in \ell_{2,p,T}} (\sum_{l=1}^\infty
      \|x_l\|_2^2)^{1/2} (\sum_{l=1}^\infty \|A_{l}\|_2^2)^{1/2}
      \textrm{ (by the Cauchy-Schwarz inequality)}\\
    &= (\sum_{l=1}^\infty \| A_{l}\|_2^2)^{1/2},
\end{align*}
thus proving the result. 
\end{proof}

We use the following result in the proof of Lemma \ref{lemma:propB} and Theorem \ref{lem.Cu.pos.def}.
\begin{lemma}\label{lem.l2.infty}
Let $\BS B$ be a symmetric linear operator from $\ell_{2,p}$ to $\ell_{2,p}$ with $\|\BS B\|_2<\infty$. 
Then,
$$
\|\BS B\|_2 \leq \max_{s_1} \sum_{s_2\in \Z} \|B_{s_1,s_2}\|_2
$$
\end{lemma}
\begin{proof}
To prove the result we define the following operator based on
$\BS B$. Let 
 $\BS {\widetilde B}=(\|B_{s_1,s_2}\|_2)_{s_1,s_2}$ be an operator
 from $\ell_2$ to $\ell_2$. Since
 $\BS B$ is symmetric, we have 
\begin{align*}
\|\BS B\|_2&=\sup_{\|x\|_2=1} x^\top \BS B x= \sup_{\|x\|_2=1} \sum_{s_1,s_2\in \Z} x_{s_1}^\top B_{s_1,s_2} x_{s_2} \leq \sup_{\|x\|_2=1} \sum_{s_1,s_2\in \Z} \|x_{s_1}\|_2 \|B_{s_1,s_2}\|_2 \|x_{s_2}\|_2\\
&=\|\BS {\widetilde B}\|_2\leq \|\BS {\widetilde
                                                                                                                                                                                                                  B}\|_\infty=\max_{s_1} \sum_{s_2\in \Z} \|B_{s_1,s_2}\|_2.
\end{align*}
This proves the result. 
\end{proof}

The following lemma is a generalisation of the Cauchy-Schwarz
inequality to the spectral norm of matrices.
\begin{lemma}\label{lemma:CauchySchwarzMV}
    Let $X$ and $Y$ be finite dimensional random vectors (not necessarily of
    the same dimension). 
Then, we have 
\begin{eqnarray*}
    \|\cov(Y,X)\|_2^2 \leq \|\var(X)\|_2 \|\var(Y)\|_2.
\end{eqnarray*}
A generalisation of the above result is to the case that
$A$ and $B$ denote two conformable random matrices. Then 
\begin{eqnarray*}
\|\Ex(AB)\|_{2}^{2} \leq \|\Ex(AA^{\top})\|_{2}\|\Ex(BB^{\top})\|_{2} 
\end{eqnarray*}
\end{lemma}
\begin{proof}
To prove the result we start by first assuming that
$\var(X)$ is strictly positive definite and later relax this condition
to the case that $\var(X)$ is non-negative definite.
 Let 
    \begin{eqnarray*}
    \var((X^\top,Y^\top)^\top)) =\begin{pmatrix}
    \Sigma_{1,1} & \Sigma_{1,2} \\
    \Sigma_{2,1} & \Sigma_{2,2} 
    \end{pmatrix}.
    \end{eqnarray*}
    Since $\Sigma$ is a positive semi-definite matrix, $\Sigma_{2,2}-\Sigma_{2,1} \Sigma_{1,1}^{-1} \Sigma_{1,2}$ is a
positive semi-definite matrix.
Hence, we have (see, for example, \cite[p.76]{lutkepohl1996handbook})
$$
\|\Sigma_{2,2}\|_2 \geq \|\Sigma_{2,1} \Sigma_{1,1}^{-1} \Sigma_{1,2}\|_2=\|\Sigma_{2,1} \Sigma_{1,1}^{-1/2}\|_2^2.
$$
Thus,
\begin{eqnarray*}
\|\cov(Y,X)\|_2^2= \| \Sigma_{2,1} \Sigma_{1,1}^{-1/2}
\Sigma_{1,1}^{1/2}\|_2^2 \leq
\| \Sigma_{2,1} \Sigma_{1,1}^{-1/2}\|_{2}^{2}
\|\Sigma_{1,1}^{1/2}\|_2^2\leq 
\|\Sigma_{2,2}\|_2 \|\Sigma_{1,1}\|_2.
\end{eqnarray*}

We now generalise the proof to the case that $\Sigma_{11} = \var(X)$
is non-negative definite. For this note that we have the eigenvalue decomposition $\var(X)=B\Lambda B^*$. In the case that $\var(X)$ is only positive semi-definite but not positive definite, we have for some $r<p$ that $\Lambda=\operatorname{diag}(\lambda_1,\dots,\lambda_r,0,\dots,0)$, where $\lambda_j>0$ are the ordered positive eigenvalues. Let $R=\{1,\dots,r\}$. We then define $Z=I_{p,R}^\top B^* X$. Then note that $\var(Z)=I_{p,R}^\top \Lambda I_{p,R}$, i.e, positive definite, $\|\var(Z)\|_2=\|\var(X)\|_2$, and $|B\|_2=1$. Furthermore, we have $X=B I_{p,R} Z$. This implies with the previous result for positive definite variances
\begin{align*}
    \|\cov(Y,X)\|_2 &= \|\cov(Y,Z) I_{p,R}^\top B^* \|_2\leq \|\var(Y)\|_2 \|\var(Z)\|_2 \|I_{p,R}^\top B\|_2 \\
    &= \|\var(X)\|_2 \|\var(Y)\|_2.
\end{align*}

For the generalisation to matrices, suppose that $A$ and $B$ are random matrices, where 
\begin{eqnarray*}
\Ex
\begin{pmatrix}
B\\
A \\
\end{pmatrix}
\begin{pmatrix}
B^{\top} & A^{\top} 
\end{pmatrix}
=\Ex\begin{pmatrix}
BB^\top  &BA^\top\\
A B^\top & AA^\top\\
\end{pmatrix}.
\end{eqnarray*}
Let $\Sigma_{1,1}=\Ex BB^\top$, $\Sigma_{1,2}=\Ex B A^\top$, $\Sigma_{2,1}=\Ex A B^\top$, and $\Sigma_{2,2}=\Ex A A^\top$. Then, we can follow the previous arguments.
\end{proof}

\begin{remark}[Generalisation of Lemma \ref{lemma:CauchySchwarzMV} to Infinite dimensional operators]
    Suppose that the eigenvalues of the symmetric positive semi-definite operator
  $\BS \Sigma$ are bounded, and 
  \begin{eqnarray*}
   \BS \Sigma =
   \left(
  \begin{array}{cc}
   \BS A & \BS B \\
  \BS  B^{*} &\BS D \\
   \end{array}
  \right).   
\end{eqnarray*}
By using the same
  arguments as those in Lemma \ref{lemma:CauchySchwarzMV} we have
  \begin{eqnarray*}
\|\BS B\|_{2}^{2} \leq \|\BS A\|_{2}\|\BS D\|_{2}.
  \end{eqnarray*}
\end{remark}

An application of the Lemma \ref{lemma:CauchySchwarzMV} is in
obtaining a bound for 
the spectral norm of the variance of infinite sums. Suppose the random
matrix $Y$ has the representation 
\begin{eqnarray*}
Y = \sum_{j=0}^{\infty}U_{j},
\end{eqnarray*}
where $\{U_{j}\}$ are random matrices. Then 
\begin{eqnarray*}
\|\Ex[YY^{\top}]\|_{2} \leq \sum_{j_1,j_2=0}^{\infty}\|\Ex[U_{j_1}U_{j_2}^{\top}]\|_{2}.
\end{eqnarray*}
By applying Lemma \ref{lemma:CauchySchwarzMV} to bound
$\|\Ex[U_{j_1}U_{j_2}^{\top}]\|_{2}$ we have
\begin{eqnarray}
\label{eq:ExY}
\|\Ex[YY^{\top}]\|_{2} \leq \left(\sum_{j=0}^{\infty}\|\Ex[U_{j}U_{j}^{\top}]\|_{2}^{1/2}\right)^{2}.
\end{eqnarray}
The above bound will be used to prove the results in Example
\ref{exam:sre}. 

The following lemma is used in the proofs of Theorems
\ref{thm.D.smooth}, 
\ref{lemma:smoothVAR} and \ref{thm.part}.
\begin{lemma} \label{lem.SSR.F1}
Let $\gu{j} = \max(1,|j|)$ and $\zeta(j) =
\gu{\log[\gu{j}]}/\gu{j}$. For all 
$y\in \R$ and $p\geq 2$ we have 
\begin{eqnarray}
\label{eq:F1.1}
\sum_{j \in \Z} \gu{j}^{-p} \gu{j+y}^{-p}\leq (\pi^2+3) \gu{y-1}^{-p}
\end{eqnarray}
and 
\begin{eqnarray}
\label{eq:F1.1b}
\sum_{j \in \Z} \zeta(j)^{p} \zeta(j+y)^{p}\leq 20 \zeta(y-1)^p
\end{eqnarray}
Further, suppose that $p,q,r\geq 2$ then
\begin{eqnarray}
\label{eq:F1.2}
\sum_{j \in \Z} \gu{j}^{-q} \gu{j+y}^{-p}\leq (\pi^2+3) \gu{y-1}^{-\min(p,q)},
\end{eqnarray}
\begin{eqnarray}
\label{eq:F1.2b}
\sum_{j \in \Z} \zeta(j)^{p} \zeta(j+y)^{q}\leq 20 \zeta(y-1)^{\min(p,q)},
\end{eqnarray}
\begin{eqnarray}
\label{eq:F1.3}
\sum_{s_1,s_2\in \Z} \gu{s_1+t}^{-p} \gu{s_1+s_2}^{-q}\gu{s_2+\tau}^{-r} \leq (\pi^2+3)^{2} \gu{t-\tau-2}^{-\min(p,q,r)},
\end{eqnarray}
and
\begin{eqnarray}
\label{eq:F1.3b}
\sum_{s_1,s_2\in \Z} \zeta(s_1+t)^{p} \zeta(s_1+s_2)^{q}\zeta(s_2+\tau)^{-r} \leq 400 \zeta(t-\tau-2)^{\min(p,q,r)}
\end{eqnarray}
\end{lemma}
\begin{proof}
First note that $\sum_{k=1}^\infty k^{-2}=\pi^2/6$. 
The strategy is to split the sum in several parts and for each part we
pull one of the factors out of say,
of $\gu{j}^{-p} \gu{j+y}^{-p}$, leverage on the pulled factor and show that the remaining sum is finite. 

We first prove (\ref{eq:F1.1}). 
Without loss of generality, let $y>0$. We have
\begin{align*}
    \sum_{j \in \Z} \gu{j}^{-p} \gu{j+y}^{-p} =I_1+I_2+I_3,
\end{align*}
where
\begin{align*}
    I_1 &=\sum_{j=0}^\infty \gu{j}^{-p} \gu{j+y}^{-p} \leq (\pi^2/6+1)
          \gu y^{-p}, \\
    I_2&=2\gu{y-1}^{-p}+
    \sum_{j=-y+2}^{-y/2} \gu{j}^{-p} \gu{j+y}^{-p}+\sum_{j=-y/2+1}^{-2} \gu{j}^{-p} \gu{j+y}^{-p}\\
    &\leq 2\gu{y-1}^{-p}+2\gu{y/2}^{-p} 2^{-p+2}
    \leq 2\gu{y-1}^{-p}+\gu{y}^{-p} 8 (\pi^2/6-1)\\
    &\leq 2\gu{y-1}^{-p}+\gu{y}^{-p}  (2/3\pi^2) \\
    I_3&=\sum_{j=-\infty}^{-y} \gu{j}^{-p} \gu{j+y}^{-p} \leq (\pi^2/6+1) \gu y^{-p}.
\end{align*}
The bounds for  $I_{1}, I_2$ and $I_{3}$ prove (\ref{eq:F1.1}).

 
To proof \eqref{eq:F1.1b}, note first that
$$\sum_{k=1}^\infty \zeta(k)^{2}=1+\sum_{k=2}^\infty \zeta(k)^{2}\leq
1+\int_1^\infty (\log(x)/x)^2 dx=1+2.$$
We will also use that $\zeta(\cdot)$ is monotonic decreasing after $\zeta(3)$, and $\zeta(1)=1$, $\zeta(2)=\zeta(4)<\zeta(3)$.
We start by follow the arguments as in the proof of (\ref{eq:F1.1}) by
splitting the sum into three parts we have $\sum_{j \in \Z} \zeta(j)^{p}
\zeta(j+y)^{p}  = I_1+I_2+I_3$ (where $I_{1},I_{2}$ and $I_{3}$ are
the same as those in the proof of  (\ref{eq:F1.1}) but with
$\zeta(\cdot)$ replacing $\gu{\cdot}^{-1}$). Without loss of generality we
prove the result for  $y\geq 3$. For $y\geq 3$ and using the
monotonicity property of $\zeta(\cdot)$ we have  
\begin{eqnarray*}
  I_1=\sum_{j=0}^\infty \zeta(j)^p \zeta(j+y)^p\leq 3 \zeta(y)^p
\end{eqnarray*}
and by the same argument
$$I_3=\sum_{j=-\infty}^{-y}\zeta(j)^p \zeta(j+y)^p\leq 3\zeta(y)^p.$$

Next we bound $I_{2}$. For this we use that
$$\zeta(y/2)^p \sum_{j=2}^\infty \zeta(j)^p\leq \zeta(y) 2^p \sum_{j=2}^\infty \zeta(j)^p \leq  \zeta(y) (\sum_{j=1}^\infty \zeta(j)^p+\sum_{j=1}^\infty \zeta(j+1/2)^p) \leq 6 \zeta(y).$$

 This gives 
\begin{align*}
I_2&=\sum_{j=-y}^{-1}\zeta(j)^p \zeta(j+y)^p=
2\zeta(y-1)^p+\sum_{j=-y+2}^{-y/2}\zeta(j)^p \zeta(j+y)^p+\sum_{j=-y/2+1}^{-2}\zeta(j)^p \zeta(j+y)^p\\
&\leq 2\zeta(y-1)^p+12 \zeta(y)\\
&\leq 14 \zeta(y-1)^p.
\end{align*}
Thus we have bounds for the terms $I_{1},I_{2}$ and $I_{3}$ in $\sum_{j \in \Z} \zeta(j)^{p}
\zeta(j+y)^{p}$, which proves \eqref{eq:F1.1b}.

The proof of (\ref{eq:F1.2}) uses that $\gu{j}^{-p}>\gu{j}^{-q}$, then 
the result immediately follows from (\ref{eq:F1.1}).

To prove  (\ref{eq:F1.3}), let us suppose wlog that $p\leq q \leq r$,
then by using (\ref{eq:F1.2}) we have 
\begin{eqnarray*}
\sum_{s_1,s_2\in \Z} \gu{s_1+t}^{-r}
  \gu{s_1+s_2}^{-p}\gu{s_2+\tau}^{-q} 
 &=& \sum_{s_1\in \Z} \gu{s_1+t}^{-r}\sum_{s_{2}\in \mathbb{Z}}
  \gu{s_1+s_2}^{-p}\gu{s_2+\tau}^{-q} \\
&\leq&  (\pi^2+3) \sum_{s_1\in \Z} \gu{s_1+t}^{-r}\gu{s_{1}-\tau-1}^{-p} \\
&\leq&  (\pi^2+3)^{2} \sum_{s_1\in \Z} \gu{t-\tau-2}^{-p}
\end{eqnarray*}
where the last two lines follow from (\ref{eq:F1.2}). This proves the
result. 
\eqref{eq:F1.2b} and \eqref{eq:F1.3b} follow analogously. 
\end{proof}

\section{Proof of results in Section~\ref{sec:inverserate}}\label{sec:proof1}


The proof of Theorem \ref{thm.1} is based on decomposing 
$\BS C^{-1}$ in terms of the inverse of a banded block matrix and its remainder, and balancing these two terms. 
 An important result on the inverse of banded matrices is given in 
\cite{demko1984decay}, Theorem 2.4. Specifically, they consider positive definite infinite dimensional matrices of the form $\BS A:\ell_{2}\rightarrow \ell_{2}$ where 
$\BS A = (A_{t,\tau};t,\tau\in
\mathbb{Z})$ ($A_{t,\tau}\in \mathbb{R}$). They show that if 
$\BS A$ has bandwidth $M$ (in the sense $A_{t,\tau}=0$ if $|t-\tau|>M$) and $\BS A^{-1} =\BS B = (B_{t,\tau};t,\tau\in
\mathbb{Z})$, then 
\begin{eqnarray}
\label{eq:demkorate}
|B_{t,\tau}| \leq \frac{(1+\sqrt{r})^2}{b} \rho^{\lfloor |t-\tau|/M\rfloor+1},
\end{eqnarray}
where $\rho=(\sqrt{r}-1)/(\sqrt{r}+1)$, $r=b/a$, $b =
  \sup_{v \in \ell_{2}, \|v\|_{2}=1}\langle v, \BS Av \rangle$, and
$a =\inf_{v \in \ell_{2}, \|v\|_{2}=1}\langle v, \BS Av
  \rangle$. An interesting application of this results is given in
 \cite{p:din-zho-21}, who use it to obtain a
  rate of decay for the parameters in an autoregressive approximation.
  As our results are in the multivariate (possibly high dimensional)
  setting we require a bound on the block entries of a banded matrix
  (and not just the individual entries). Thus in the following lemma
  we obtain a generalisation of (\ref{eq:demkorate}) for block matrices.
  
\begin{lemma}
\label{lem.demko}
Let $\BS A$ be a positive definite linear operator on $\ell_{2,p}$ 
where $\BS A = (A_{t,\tau};t,\tau\in \mathbb{Z})$ and 
$A_{t,\tau}$ is a $p\times p$ dimensional matrix. We suppose that 
$\BS A$ is block-banded with bandwidth $M$ and block-size $p$ in the
sense that for all $s_1,s_2$ with $|s_1-s_2|>M$, $A_{s_1,s_2}=0$.
Let $b =  \sup_{v \in \ell_{2,p}, \|v\|_{2}=1}\langle v, \BS Av \rangle,$ and 
$a =  \inf_{v \in \ell_{2,p}, \|v\|_{2}=1}\langle v, \BS Av
  \rangle$. Furthermore, $r=b/a$, $\rho=(\sqrt{r}-1)/(\sqrt{r}+1)$. 
Let $\BS B = \BS A^{-1} = (B_{t,\tau};t,\tau\in \mathbb{Z})$ (where $B_{t,\tau}$ is a $p\times p$ dimensional matrix).
Then, the following bound holds for all $p\times p$ sub-matrices and $t\not = \tau$
$$
\|B_{t,\tau}\|_2 \leq \frac{(1+\sqrt{r})^2}{b} \rho^{\lfloor |t-\tau|/M\rfloor+1}
$$
where $\lfloor x \rfloor$ denotes the largest integer less than or
equal to $x$.

Let $\BS{\tilde A}=(I_{-k} \otimes I_p)^\top \BS A (I_{-k} \otimes
I_p)$ be a sub-matrix without the $kth$ $p$-dimensional
row and column, where $k \in \Z$. Then, for $\tilde {\BS B}=\tilde {\BS A}^{-1}$ with $\tilde B_{t,\tau}=(((I_{-k} \otimes I_p) \tilde {\BS A}^{-1} (I_{-k} \otimes I_p)^\top )_{t,\tau};t,\tau\in \mathbb{Z})$ the following bound holds for all $p\times p$ sub-matrices and $t\not = \tau$
$$
\|\tilde B_{t,\tau}\|_2 \leq \frac{(1+\sqrt{r})^2}{b} \rho^{\lfloor |t-\tau|/M\rfloor+1}.
$$
\end{lemma}
\begin{proof}
The proof is based on the proof of Proposition 2.2 in
\cite{demko1984decay}, with a modification to allow for block
matrices. We use the notation from Proposition 2.2 in
\cite{demko1984decay}. More precisely, 
let $\Pi_n$ denote the space of polynomials up to order
$n$. A
key ingredient in the proof is the following classical result from spectral
theory. Suppose $A$ is a positive definite operator, then
\begin{eqnarray*}
\|\BS A^{-1}-p(\BS A)\|_2\leq\max_{x \in [a,b]}|1/x-p(x)|,
\end{eqnarray*}
where $p$ is a real polynomial and recall $b =
  \sup_{v \in \ell_{2,p}, \|v\|_{2}=1}\langle v, \BS Av \rangle$ and 
$a =
  \inf_{v \in \ell_{2,p}, \|v\|_{2}=1}\langle v, \BS Av \rangle$. Set
  $r=b/a$, $\rho=(\sqrt{r}-1)/(\sqrt{r}+1)$.  For any complex valued
  function $f$ on $K$, define the norm 
  $\|f\|_K=\sup\{|f(z)| : z \in
  K\}$ (thus $\|1/x-p(x)\|_{[a,b]}=\max_{x \in
    [a,b]}|1/x-p(x)|$). Proposition 2.1, 
\cite{demko1984decay} show that 
  \begin{align} \label{eq.lem.1}
  \inf\{\|1/x-p(x)\|_{[a,b]} : p \in \Pi_n \}=\frac{(1+\sqrt{r})^2}{b} \rho^{n+1}.      
  \end{align}
Using this result we define the polynomial 
\begin{eqnarray}
\label{eq:pstar}
p_n^{*}=\arg_{p \in \Pi_n}\inf
\{\|1/x-p(x)\|_{[a,b]} : p \in \Pi_n \}.
\end{eqnarray}
We note for any polynomial $p_n$ of order $n$ and $M$
block-banded matrix $\BS A$ with block size $p$, 
if $|t-\tau|\geq nM$  then  $p_n(\BS A)_{t,\tau}
\equiv 0$ where 
$p_n(\BS A))_{t,\tau}$ denotes the 
$(t,\tau)$ $p\times p$ dimension block matrix in $p_n(\BS A)$.

For a given $t$ and $\tau$, set $n=\lfloor |t-\tau|/M\rfloor$. Let $p_{n}^{*}$ be defined as in (\ref{eq:pstar}).
Then by definition of $n$ 
we have $p_n^{*}(\BS A)_{t,\tau}=0$.
Since $B_{t,\tau} = (\BS A^{-1})_{t,\tau}$ this 
gives 
\begin{eqnarray*}
\|B_{t,\tau}\|_2&=\|(\BS A^{-1}-p_n^{*}(\BS A))_{t,\tau}\|_2\leq \|\BS A^{-1}-p_n^{*}(\BS A)\|_2=
\|1/x-p_n^*(x)\|_{[a,b]}\\
&=\frac{(1+\sqrt{r})^2}{b} \rho^{\lfloor |t-\tau|/M\rfloor+1},
\end{eqnarray*}
where the last part follows from (\ref{eq.lem.1}). 
This completes the proof of the first assertion.  

For the second assertion, our strategy is to extend $\BS{\tilde A}$ such that is an operator from $\ell_{2,p}$ to $\ell_{2,p}$ and possesses the same banded-scheme as $\BS A$. Then, we apply the results we derived in the first assertion to this extended $\BS{\tilde A}$, hence we obtain an inverse with the desired properties. Lastly, we show that when shrinking the inverse of the extended $\BS{\tilde A}$ to the space of $\BS{\tilde A}$, we obtain an inverse of $\BS{\tilde A}$. 
This idea can be formalised as follows. The extended  $\BS{\tilde A}$ is obtained by 
$(I_{-k} \otimes I_p) \BS{\tilde A} (I_{-k} \otimes I_p)^\top+c
(e_{k} \otimes I_p)(e_{k} \otimes I_p)^\top=:\BS E$. We have that $\BS E$
is a block-banded matrix. Additionally, if we set $c=\|(e_{k} \otimes I_p)^\top \BS A (e_{k} \otimes I_p)\|_2$ the largest and smallest eigenvalues of $E$ can be bounded by those of $A$, that is its largest and smallest eigenvalues are bounded above from $b$ and below from $a$, respectively. Hence, the previous assertion applies to $\BS E$. 
We now show $(I_{-k} \otimes I_p)^\top \BS D^{-1} (I_{-k} \otimes
I_p)=(\BS{\tilde A})^{-1}$ which gives the assertion. For this, we
show $(I_{-k} \otimes I_p)^\top \BS D^{-1} (I_{-k} \otimes I_p) \tilde
{\BS A} =(I_{-k} \otimes I_p)^\top (I_{-k} \otimes I_p)$ and use the
uniqueness of the inverse operator.  The calculation is
\begin{align*}
    (I_{-k} \otimes I_p)^\top \BS D^{-1} (I_{-k} \otimes I_p) \tilde {\BS A}=&
    (I_{-k} \otimes I_p)^\top\! ((I_{-k} \otimes I_p) \BS{\tilde A} (I_{-k} \otimes I_p)^\top\!\!\!+c (e_{k} \otimes I_p)(e_{k} \otimes I_p)^\top)^{-1} \\
    &\times 
    (I_{-k} \otimes I_p)(\tilde {\BS A}(I_{-k} \otimes I_p)^\top+c (e_{k} \otimes I_p)(e_{k} \otimes I_p)^\top\\
    &-c (e_{k} \otimes I_p)(e_{k} \otimes I_p)^\top) (I_{-k} \otimes I_p) \\
    =&(I_{-k} \otimes I_p)^\top (I_{-k} \otimes I_p)+0.
\end{align*}
Thus, $(I_{-k} \otimes I_p)^\top \BS D^{-1} (I_{-k} \otimes I_p)$ is an inverse of $\tilde {\BS A}$ and the second assertion follows.
\end{proof}

We now apply the above result to a specific banded matrix (required 
in the proof of Theorem \ref{thm.1}).  Define the integer set 
$t^{c} = \{\tau\in \mathbb{Z}, \tau\neq t\}$, and 
$\BS C_{t^{c},t^{c}}=(I_{-t} \otimes I_p)^\top \BS C (I_{-t}\otimes I_p)$
(this is operator $\BS C$ but with the t$th$ block row and column removed). 
We define $\BS B_{M}$ as the $M$th banded version of $\BS C_{t^{c},t^{c}}$ as
follows. For all $p\times p$ sub-matrices and $s_1,s_2\in \Z$ let
\begin{eqnarray}
\label{eq:BSdefine}
((I_{-t} \otimes I_p) \BS B_M (I_{-t} \otimes I_p)^\top)_{s_1,s_2}=\ind(|s_1-s_2|\leq M)
  ((I_{-t} \otimes I_p) \BS C_{t^{c},t^{c}} (I_{-t} \otimes I_p)^\top)_{s_1,s_2}, 
\end{eqnarray}
where $\ind$ denotes the indicator function.

The following lemma is used in the proof of Theorems~\ref{thm.1} and \ref{lemma:spectral}.
\begin{lemma}\label{lemma:propB}[Properties of ${\BS B}_{M}$]
Suppose Assumption \ref{ass} is satisfied and let ${\BS B}_{M}$ be a
(symmetric) banded matrix 
defined as in (\ref{eq:BSdefine}).
Define the space of vectors $$\ell_{2,p}^{-t} =
\{ v = (\ldots,v_{t-1},v_{t+1},v_{t+2},\ldots); v_{j}\in
\mathbb{R}^{p}, \sum_{j\neq t}\|v_{j}\|_2^{2}<\infty\}$$ and the 
eigenvalues 
\begin{eqnarray*}
a_{M} =  \inf_{v \in \ell_{2,p}^{-t},
  \| v\|_{2}=1}\langle v,\BS B_M v \rangle \textrm{ and }
 b_M= \sup_{v \in \ell_{2,p}^{-t}, \|v\|_{2}=1}\langle v,
\BS B_M v \rangle.
\end{eqnarray*} 
Then  
\begin{eqnarray}
\label{eq:BM1}
\|\BS C_{t^{c},t^{c}}-\BS B_M\|_2
\leq 2 \frac{K}{(\kappa-1)} (M-1)^{-\kappa+1},
\end{eqnarray}
\begin{eqnarray}
\label{eq:Bmeigenvalues}
a_{M}\geq \lambda_{\inf}-2 \frac{K}{(\kappa-1)} (M-1)^{-\kappa+1},
 b_{M} \leq 
\lambda_{\sup}+2 \frac{K}{(\kappa-1)} (M-1)^{-\kappa+1}
\end{eqnarray}
and if $M$ is such that $\lambda_{\inf}-2 \frac{K}{(\kappa-1)}
(M-1)^{-\kappa+1}>0$, then
\begin{eqnarray}
\label{eq:Bspec}
\|\BS B_{M}^{-1}\|_2\leq \left(\lambda_{\inf}-2 \frac{K}{(\kappa-1)} (M-1)^{-\kappa+1}\right)^{-1}.
\end{eqnarray}
The same rates to the banded matrices associated with $\BS C$ or $\BS C(-\infty;T]$.
\end{lemma}
\begin{proof} We first prove (\ref{eq:BM1}). For this, we first expand $\BS C_{t^{c},t^{c}}-\BS B_M$ with zero such that it is an operator from $\ell_{2,p}$ to $\ell_{2,p}$ again. Then, we apply Lemma~\ref{lem.l2.infty} to obtain
\begin{align*}
    \|\BS C_{t^{c},t^{c}}-\BS B_M\|_2=&\leq \sup_{s_1} \sum_{s_2} \| ((I_{-k} \otimes I_p)(\BS C_{t^{c},t^{c}}-\BS B_M) (I_{-k} \otimes I_p))_{s_1,s_2}\|_2 \\
    \leq& \sum_{|s|>M} \frac{K}{|s|^{-\kappa}} \leq 2 K \sum_{s>M} \int_{s-1}^s x^{-\kappa} dx=
    2 \frac{K}{(\kappa-1)} (M-1)^{-\kappa+1}.
\end{align*}
where the first bound on last line above follows from Assumption \ref{ass}. 

To prove (\ref{eq:Bmeigenvalues}) we use that $\BS B_M=\BS
C_{t^{c},t^{c}}+(\BS C_{t^{c},t^{c}}-\BS B_M)$
and the eigenvalues of $\BS C_{t^{c},t^{c}}$ are in
$[\lambda_{\inf},\lambda_{\sup}]$.
Thus, with \eqref{eq:BM1} we have
\begin{eqnarray}
\label{eq:Bmeigenvalues2}
\lambda_{\inf}(\BS B_{M})\geq \lambda_{\inf}-2 \frac{K}{(\kappa-1)} (M-1)^{-\kappa+1} \textrm{ and }
\lambda_{\sup}(\BS B_{M})  \leq 
\lambda_{\sup}+2 \frac{K}{(\kappa-1)} (M-1)^{-\kappa+1}.
\end{eqnarray}
The proof of (\ref{eq:Bspec}) immediately follows from (\ref{eq:Bmeigenvalues}).
\end{proof}


Using the above lemma we now prove Theorem \ref{thm.1}.

\begin{proof}[Proof of Theorem~\ref{thm.1}]
For \eqref{eq:Dttau} we focus here on the case $t\not = \tau$ and
$|t-\tau|\geq 2$. 

To motivate the proof, 
we first describe a more direct but naive approach which does not give sufficiently
sharp bounds. One strategy is to directly approximate $\BS D$
with the inverse of a (block) banded matrix, say $\BS A_{M}$, and then use the 
Neuman series to bound its error. I.e. use an expansion of the form
\begin{eqnarray*}
  \BS D = (\BS A_{M}+(\BS D - \BS A_{M}))^{-1}  =
  \BS A_{M}^{-1} + \sum_{j=1}^{\infty}[\BS A_{M}^{-1} (\BS D - \BS A_{M})]^{j},
\end{eqnarray*}
which holds when $M$ is large enough such that $\|\BS A_{M}^{-1}(\BS D -
\BS A_{M})\|_2<1$. The $\ell_2$ bound of the above is 
\begin{eqnarray*}
\|\BS D_{t,\tau}\|_{2} &\leq&  \|(\BS A_{M}^{-1})_{t,\tau}\|_{2}+
  \sum_{j=1}^{\infty}\|\BS A_{M}^{-1} (\BS D - \BS A_{M})\|_{2}^{j} \\
&\leq&  \|\BS A_{M}^{-1}\|_{2}+
  \|\BS A_{M}^{-1} (\BS D - \BS A_{M})\|_{2}^{}(1-\|\BS A_{M}^{-1}
       (\BS D - \BS A_{M})\|_{2})^{-1} \\
&=& I_{1} + I_{2}.
\end{eqnarray*}
By using  Lemma \ref{lem.demko} we can show that
\begin{eqnarray*}
I_{1}\leq 
  \frac{(1+\sqrt{r_M})^2}{\lambda_{\sup,M}}\rho_M^{\lfloor|s-\tau|/M\rfloor+1},
\end{eqnarray*}
where
$r_M=\lambda_{\sup,M}/\lambda_{\inf,M}$ 
$\rho_M=(\sqrt{r_M}-1)/(\sqrt{r_M}+1)$ and $\lambda_{\inf,M}$ and
$\lambda_{\sup,M}$ are the eigenvalues of $\BS A_{M}$. It can also be
shown that
\begin{eqnarray*}
I_{2} \leq \frac{C}{(\kappa-1)}
(M-1)^{-\kappa+1}.
\end{eqnarray*}
This leads to the rate
\begin{eqnarray*}
\|\BS D_{t,\tau}\|_{2} &\leq&
                              \frac{(1+\sqrt{r_M})^2}{\lambda_{\sup,M}}\rho_M^{\lfloor|s-\tau|/M\rfloor+1}+
     \frac{C}{(\kappa-1)}
(M-1)^{-\kappa+1}.                         
\end{eqnarray*}
However, this bound does not adequately 
utilize the rate of decay of the entries of $\BS C$. Instead we take
an indirect approach, where we rewrite $\BS D$ as the inverse of a block
matrix, where the relevant entries are the inverse of one submatrix of
$\BS C$ \emph{multipled} with another submatrix of $\BS C$. The latter term
allows us to leverage on the rate of decay of the entries of $\BS C$.
We describe this approach below. 

Define the integer set $t^{c} = \{\tau\in \mathbb{Z}, \tau\neq t\}$, and denote 
$\BS C_{t,t^{c}}=(e_t \otimes I_p)^\top \BS C (I_{-t} \otimes I_p)$ and 
$\BS C_{t^{c},t^{c}}=(I_{-t} \otimes I_p)^\top \BS C (I_{-t}
\otimes I_p)$. Without loss of generality we consider a permuted
version of $\BS C$, which contains $C_{t,t}$ in the top left hand corner 
of $\BS C$, where 
\begin{eqnarray*}
\BS C = 
\left(
\begin{array}{cc}
C_{t,t} & \BS C_{t,t^{c}} \\
\BS C_{t^{c},t} & \BS C_{t^{c},t^{c}} \\
\end{array}
\right).
\end{eqnarray*}
Setting $\BS U = {\BS C}$, $\BS A = C_{t,t}$, $\BS B =
\BS C_{t,t^{c}}$, $\BS C= \BS C_{t,t^{c}}^{\top}$, $\BS D =
\BS C_{t^{c},t^{c}}$ and applying the block matrix operator  inversion formula in
(\ref{eq:inverseblock1}) we have 
\begin{eqnarray*}
\BS D=\BS C^{-1} = 
\left(
\begin{array}{cc}
D_{t,t} \quad & -D_{t,t}^{-1} \BS C_{t,t^{c}}\BS C_{t^{c},t^{c}}^{-1} \\
-\BS C_{t^{c},t^{c}}^{-1} \BS C_{t^{c},t}D_{t,t}^{-1}\quad
& (\BS C_{t^{c},t^{c}}-\BS C_{t^{c},t}C_{t,t}^{-1}\BS C_{t,t^{c}})^{-1}\\
\end{array}
\right).
\end{eqnarray*}
Using the above $D_{t,\tau}$ can be written as 
\begin{eqnarray*}
D_{t,\tau} = -D_{t,t}^{-1} \BS C_{t,t^{c}}\BS C_{t^{c},t^{c}}^{-1} (I_{-t} \otimes I_p)^\top (e_{\tau} \otimes I_p),
\end{eqnarray*}
using that $\lambda_{\sup}^{-1}\leq D_{t,t}\leq \lambda_{\inf}^{-1}$
we have
$
\|D_{t,\tau}\|_{2} 
\leq \lambda_{\sup}\|(\BS C_{t,t^{c}}\BS C_{t^{c},t^{c}}^{-1}) (I_{-t}^\top e_\tau \otimes I_p)\|_{2}.
$
Thus for the remainder of the proof, we focus on bounding the induced $\ell_2$-norm of 
$$A_{t,\tau} = (\BS C_{t,t^{c}}\BS C_{t^{c},t^{c}}^{-1})
(I_{-t}^\top e_\tau \otimes I_p).$$

An outline in the proof is to (a) replace
${\BS C}_{t^{c},t^{c}}^{-1}$ with the inverse of a (block) banded matrix (b)
use the Neuman series to obtain a
bound on the replacement error and (c) finally balance the rate of decay of the
inverse banded matrix approximation of $(\BS C_{t,t^{c}}\BS C_{t^{c},t^{c}}^{-1})
(I_{-t}^\top e_\tau \otimes I_p)$ with the spectral norm of  the
approximation error (both
of which depend on the bandwidth $M$).

Let $\BS B_{M}$ denote the $M$th banded matrix version of $\BS C_{t^{c},t^{c}}$;
the precise definition is given in (\ref{eq:BSdefine}).  
By Lemma \ref{lemma:propB},  equation (\ref{eq:BM1}) we have the bound
\begin{eqnarray}
\label{eq:CB-diff}
\|\BS C_{t^{c},t^{c}}-\BS B_M\|_{2} \leq 2 \frac{K}{(\kappa-1)}
(M-1)^{-\kappa+1}.
\end{eqnarray} 
Using $\BS B_{M}$ we write $\BS C_{t^{c},t^{c}}^{-1}$ as a Neumann series 
\begin{eqnarray*}
\BS C_{t^{c},t^{c}}^{-1}=\BS B_M^{-1}[I+\BS
B_M^{-1}(\BS C_{t^{c},t^{c}}-\BS B_M)]^{-1} =  \BS B_M^{-1} [I+\sum_{s=1}^\infty (-1)^s [\BS B_M^{-1}(\BS C_{t^{c},t^{c}}-\BS B_M)]^s],
\end{eqnarray*}
noting that the above expanion holds, 
when
$M>1+(2K/(\kappa-1))^{1/(\kappa-1)}$ thus $\|{\BS
  C}_{t^{c},t^{c}}-\BS B_M\|_2<1$.
Substituting the above into $A_{t,\tau} =
(\BS C_{t,t^{c}}\BS C_{t^{c},t^{c}}^{-1}) (I_{-t}^\top e_\tau \otimes I_p)$
gives for all $t,\tau\in \mathbb{Z}$
\begin{eqnarray*}
A_{t,\tau} &=& \BS C_{t,t^{c}} \BS B_M^{-1} [I+\sum_{s=1}^\infty (-1)^s [\BS B_M^{-1}(\BS C_{t^{c},t^{c}}-\BS B_M)]^s] (I_{-t}^\top e_\tau \otimes I_p)\\
&=& \BS C_{t,t^{c}}\BS B_M^{-1} (I_{-t}^\top e_\tau \otimes I_p) + \BS C_{t,t^{c}}\BS B_M^{-1} \sum_{s=1}^\infty (-1)^s [\BS B_M^{-1}(\BS C_{t^{c},t^{c}}-\BS B_M)]^s] (I_{-t}^\top e_\tau \otimes I_p). 
\end{eqnarray*}
By applying the triangular inequality to the above we have $\|A_{t,\tau}\|_{2}
\leq  J_{1,t,\tau} + J_{2,t,\tau}$ where 
\begin{eqnarray*}
J_{1,t,\tau} &=&\|\BS C_{t,t^{c}}\BS B_M^{-1} (I_{-t}^\top e_\tau \otimes
          I_p) \|_{2} \nonumber\\
\textrm{ and } J_{2,t,\tau}&=& \| \BS C_{t,t^{c}}\BS B_M^{-1} \sum_{s=1}^\infty (-1)^s
         [\BS B_M^{-1}(\BS C_{t^{c},t^{c}}-\BS B_M)]^s] 
(I_{-t}^\top e_\tau \otimes I_p)\|_{2}.     
\end{eqnarray*}
We now bound $J_{1,t,\tau}$ and $J_{2,t,\tau}$.  By using the 
sub-multiplicativity  of $\|\cdot\|_2$ we bound 
$J_{2,t,\tau}$ with
\begin{eqnarray}
J_{2,t,\tau} &\leq&  \| \BS C_{t,t^{c}}\BS B_M^{-1}\|_{2} \sum_{s=1}^\infty 
         (\|[\BS B_M^{-1}\|_{2}\|(\BS C_{t^{c},t^{c}}-\BS
             B_M)]\|_{2})^{s}. 
\end{eqnarray}
By using Lemma \ref{lemma:propB}, if $M$ is such that
$\lambda_{\inf}-2 \frac{K}{(\kappa-1)}
                                (M-1)^{-\kappa+1}>0$, then
\begin{eqnarray*}
\|[\BS B_M^{-1}\|_{2}\|(\BS C_{t^{c},t^{c}}-\BS
             B_M)]\|_{2} &\leq& 
\left(\lambda_{\inf}-2 \frac{K}{(\kappa-1)}
                                (M-1)^{-\kappa+1}\right)^{-1}2
                                \frac{K}{(\kappa-1)}
                                (M-1)^{-\kappa+1}.
\end{eqnarray*}
Thus for 
$$M>1+(2K/[\min(1,\lambda_{\inf}/2)(\kappa-1)])^{1/(\kappa-1)}:=M_{2},$$
we have $\|[\BS B_M^{-1}\|_{2}\|(\BS C_{t^{c},t^{c}}-\BS
             B_M)]\|_{2} <1$. Hence, we obtain the geometric sum 
\begin{align*}
J_{2,t,\tau}&\leq \|\BS C_{t,t^{c}}\BS B_M^{-1} \|_2  \sum_{s=1}^\infty \|[\BS B_M^{-1}\|_{2}\|(\BS C_{t^{c},t^{c}}-\BS
             B_M)]\|_{2} ^s\\
             &\leq 
\frac{2 K/(\kappa-1)(M-1)^{-\kappa+1}}{1-
\|[\BS B_M^{-1}\|_{2}\|(\BS C_{t^{c},t^{c}}-\BS
             B_M)]\|_{2}}
= 2 \mathcal{K}/(\kappa-1)(M-1)^{-\kappa+1}=:\widetilde{J}_{2,t,\tau}, 
\end{align*}
where the last line of the above follows from Lemma \ref{lemma:propB}.
In summary, for $M>M_{2}$ we have 
\begin{eqnarray}
J_{2,t,\tau}&\leq& 
2 \mathcal{K}/(\kappa-1) (M-1)^{-\kappa+1} =
                   \widetilde{J}_{2,t,\tau}.
                   \label{eq.J2tilde}    
\end{eqnarray}
Next we bound $J_{1,t,\tau}$. We start by expanding  $\BS
C_{t,t^{c}}\BS B_M^{-1}$, then use the sub-multiplicativity
of $\|\cdot\|_2$ to give 
\begin{eqnarray}
\label{eq:J1bound}
   J_{1,t,\tau} &\leq &  \sum_{s\in \Z, s \not=t} \| C_{t,s}\|_2 \cdot
                 \|((I_{-t} \otimes I_p)  \BS B_M^{-1} (I_{-t} \otimes
                 I_p)^{\top})_{s,\tau} \|_{2}.
\end{eqnarray}
We bound the terms inside of the sum $\sum_{s\in \Z, s \not=t} \| C_{t,s}\|_2 \cdot
                 \|((I_{-t} \otimes I_p)  \BS B_M^{-1} (I_{-t} \otimes
                 I_p)^\top)_{s,\tau} \|_{2}$. 
Under Assumption~\ref{ass} we have $\| C_{t,s}\|_2\leq
K\gu{t-s}^{-\kappa}$. To bound the second term, we use Lemma~\ref{lem.demko}  
\begin{eqnarray}
\label{eq:BMM1}
\|((I_{-t} \otimes I_p)  \BS B_M^{-1} (I_{-t}
\otimes I_p)^{\top})_{s,\tau}\|_{2}\leq
  \frac{(1+\sqrt{r_M})^2}{\lambda_{\sup,M}}\rho_M^{\lfloor|s-\tau|/M\rfloor+1},
\end{eqnarray}
where $r_M=\lambda_{\sup,M}/\lambda_{\inf,M}$ 
$\rho_M=(\sqrt{r_M}-1)/(\sqrt{r_M}+1)$ and $\lambda_{\sup,M}$
and $\lambda_{\inf,M}$ are such that 
\begin{eqnarray*}
\lambda_{\sup,M}&\leq&
\lambda_{\sup}+2 \frac{K}{\kappa-1} (M-1)^{-\kappa+1}
\textrm{ and }
\lambda_{\inf,M}\geq \lambda_{\inf}-2 \frac{K}{\kappa-1} (M-1)^{-\kappa+1}.
\end{eqnarray*}
This gives a bound for $r_{M}, \lambda_{\sup,M}$ and $\rho_{M}$ in
terms of $r,\lambda_{\sup}, \rho$ and $M$. 
To remove the dependency  of $M$ in these we choose 
$M$ such that 
\begin{eqnarray*}
M>\left(\frac{2K}{\kappa-1}
\max(2\lambda_{\inf}^{-1},\lambda_{\sup}^{-1})\right)^{1/(\kappa-1)}+1
  :=M_{1}
\end{eqnarray*}
For $M > M_{1}$ we have 
$\lambda_{\inf,M}\geq \lambda_{\inf}/2$ and $\lambda_{\sup,M}\leq 2 \lambda_{\sup}$.
This means, $r_M\leq 4 r$ and $\rho_M\leq (2
\sqrt{r}-1)/(2\sqrt{r}+1)=:\rho,$ where
$r=\lambda_{\sup}/\lambda_{\inf}$. Substituting this into
(\ref{eq:BMM1}) gives 
\begin{eqnarray}
\label{eq:BMM}
\|((I_{-t} \otimes I_p)  \BS B_M^{-1} (I_{-t}
\otimes I_p)^{\top})_{s,\tau}\|_{2}\leq
  \frac{2(1+2\sqrt{r})^2}{\lambda_{\sup}} \rho^{\lfloor|s-\tau|/M\rfloor+1}.
\end{eqnarray}
Substituting (\ref{eq:BMM}) and $\| C_{t,s}\|_2\leq
K\gu{t-s}^{-\kappa}$ into (\ref{eq:J1bound}) we have
\begin{eqnarray}
 J_{1,t,\tau} &\leq\frac{2K(1+2\sqrt{r})^2}{\lambda_{\sup}} \sum_{s \in Z, s \not = t}  |s-t|^{-\kappa} \rho^{\lfloor|s-\tau|/M\rfloor+1} \nonumber\\
 &\leq \frac{2K(1+2\sqrt{r})^2}{\lambda_{\sup}}  \sum_{s\in Z}
  \rho^{|s|/M} \frac{1}{\gu{s-t+\tau}^\kappa}= \widetilde{J}_{1,t,\tau}. \label{eq.J1tilde}
\end{eqnarray}
Thus when $M>K_{c}$ where 
\begin{eqnarray*}
K_{c}:=\max(M_{1},M_{2})=\left(\frac{2K}{\kappa-1}
\max(2\lambda_{\inf}^{-1},\lambda_{\sup}^{-1},1)\right)^{1/(\kappa-1)}+1
\end{eqnarray*}
 the bounds $J_{1,t,\tau}$ and $J_{2,t,\tau}$ in (\ref{eq.J1tilde}) and
 (\ref{eq.J2tilde}) hold and we have
 \begin{eqnarray*}
\|D_{t,\tau}\|_{2}\leq \lambda_{\sup}
 (\widetilde{J}_{1,t,\tau}+\widetilde{J}_{2,t,\tau}),
\end{eqnarray*}
 where $\widetilde{J}_{1,t,\tau}$  and
 $\widetilde{J}_{2,t,\tau}$ are defined in (\ref{eq.J1tilde}) and
 (\ref{eq.J2tilde}) respectively. 
 
 The final part in the proof  is to balance the two bounds $\widetilde{J}_{1,t,\tau}$
and $\widetilde{J}_{2,t,\tau}$. For each 
$t,\tau\in \mathbb{Z}$ we set 
$M=M_{t-\tau}:=-\frac{|t-\tau| \log(\rho)}{2(\kappa-1) \log(|t-\tau|)}$ (note
$0<\rho <1$). When $|t-\tau|$ is 
sufficiently large i.e.,
$M_{t-\tau}\geq K_{c}$ by substituting $M_{t-\tau}$ into the bounds for $J_{1,t,\tau}$ and $J_{2,t,\tau}$ it can be shown that 
\begin{eqnarray}
\label{eq:Dttaubound1}
 \|D_{t,\tau}\|_2&\leq& 2K(1+2\sqrt{r})^2(2^\kappa+2 S_\kappa) |t-\tau|^{-\kappa+1}+
 \frac{2
                        K}{\kappa-1}\left(\frac{|\log(\rho)|}{2(\kappa-1)}\frac{|t-\tau|}{\log|t-\tau|}-1\right)^{-\kappa+1}.
\end{eqnarray}
Note that the above expression, though unwieldy gives the desired decay
$\zeta(t-\tau)^{\kappa-1}$.  However, if 
$|t-\tau|$ is small, i.e., 
$M_{t-\tau} \leq K_{c}$ the bound (\ref{eq:Dttaubound1}) does not hold and 
we use an alternative bound for
$\|D_{t,\tau}\|_2$. 
It is easily seen that $\|D_{t,\tau}\|_2\leq \|\BS D\|_2\leq \lambda_{\inf}^{-1}$. 
We rewrite the above in a similar form as (\ref{eq:Dttaubound1}) (but
with different contants)
\begin{eqnarray}
\label{eq:Dttaubound2}
\|D_{t,\tau}\|_2\leq \lambda_{\inf}^{-1}
\leq (\lambda_{\inf}\min(\lambda_{\inf}/2,\lambda_{\sup}))^{-1} \frac{2K}{\kappa-1}\left(\frac{|\log(\rho)|}{2(\kappa-1)}\frac{|t-\tau|}{\log|t-\tau|}-1\right)^{-\kappa+1}.
\end{eqnarray}
Combining (\ref{eq:Dttaubound1}) and (\ref{eq:Dttaubound2})
gives the following global bound for all $t,\tau\in \mathbb{Z}$ 
 \begin{eqnarray*}
 \|D_{t,\tau}\|_2&\leq& 2K(1+2\sqrt{r})^2(2^\kappa+2 S_\kappa) \gu{t-\tau}^{-\kappa+1}\\
 && +\max(1,(\min(\lambda_{\inf}/2,\lambda_{\sup})\lambda_{\inf})^{-1})\frac{2 K}{(\kappa-1)}\left(\frac{|\log(\rho)|}{2(\kappa-1)}\frac{\gu{t-\tau}}{\gu{\log\gu{t-\tau}}}-1\right)^{-\kappa+1} \\
 &\leq& \mathcal{K}\zeta(t-\tau)^{\kappa-1}.
\end{eqnarray*}
Note that in the proof we have carefully tracked all the constants,
to demonstrate that the constants only depend on $\lambda_{\inf}, \lambda_{\sup}, K$
and $\kappa$. To reduce notation, in the remainder of the paper we
use a generic constant $\mathcal{K}$. 

To prove (\ref{eq:DttauInv}), we only need a slight modification of
the above arguments. We define the integer set $t_{T}^{c}=\{\tau \leq T, \tau\not=t\}$ and obtain
\begin{eqnarray*}
\BS C(-\infty;T) = 
\left(
\begin{array}{cc}
C_{t,t} & \BS C_{t,t_{T}^{c}} \\
\BS C_{t_{T}^{c},t} & \BS C_{t_{T}^{c},t_{T}^{c}} \\
\end{array}
\right).
\end{eqnarray*}
This leads to 
$$[\BS C(-\infty;T)^{-1}]_{t,\tau}=-[\BS C(-\infty;T)^{-1}]_{t,t}^{-1}
\BS C_{t,t_{T}^{c}}\BS C_{t_{T}^{c},t_{T}^{c}}^{-1} 
(I_{-t} \otimes I_p)^\top (e_{\tau} \otimes I_p).$$ Then, we follow the same strategy as above. Note that the sums occurring now are going from from $-\infty$ to $T$ instead before in deriving \eqref{eq:Dttau} in which they are from $-\infty$ to $\infty$.
\end{proof}

\begin{proof}[Proof of Theorem~\ref{lemma:onesided}] 
We first recall that the coefficients $\{\Phi_{T,j}\}$ and
$\Phi_{T,d,j}$ are embedded in the last rows of 
 $\BS C(-\infty,T)^{-1}$ and $\BS C(T-d,T)^{-1}$  respectively. Therefore 
we first need to connect the inverses of 
$\BS C(-\infty,T)$ and $\BS C(T-d,T)$. For this, we write ${\BS C}(-\infty,T)$ in terms of the following block matrix
\begin{eqnarray*}
{\BS C}(-\infty,T) = 
\left(
\begin{array}{cc}
\BS C(-\infty,T-d) & \BS C(-\infty,T-d,T) \\
\BS C(-\infty,T-d,T)^{\top} & \BS C(T-d,T) \\
\end{array}
\right),
\end{eqnarray*}
where $\BS C(-\infty,T-d,T) = (C_{t,\tau};t\leq T-d,T-d+1\leq \tau  \leq T)$ and
$\BS C(T-d,T) = (C_{t,\tau};T-d+1\leq t,\tau\leq T)$. Next we represent
${\BS C}(-\infty,T) ^{-1}$ as a block operator (analogous to ${\BS C}(-\infty,T) $)
\begin{eqnarray*}
{\BS C}(-\infty,T) ^{-1}=
\left(
\begin{array}{cc}
\widetilde{\BS D}(-\infty,T-d) & \widetilde{\BS D}(-\infty,T-d,T) \\
\widetilde{\BS D}(-\infty,T-d,T)^{\top} & \widetilde{\BS D}(T-d,T) \\
\end{array}
\right).
\end{eqnarray*}
Note we have used the notation $\widetilde{\BS D}$ to show that they
are not the inverse of the corresponding submatrix of $\BS C$.
To evaluate $\BS C(T-d,T)^{-1} - \widetilde{\BS D}(T-d,T)$ we
 apply the second identity in (\ref{eq:schur2})  where we set
$\BS U = {\BS C}(-\infty,T) ^{-1}$, 
$\BS A = \widetilde{\BS D}(-\infty,T-d)$, $\BS B = \widetilde{\BS
  D}(-\infty,T-d,T)$, $\BS C = \widetilde{\BS D}(-\infty,T-d,T)^{\top}$,
$\BS D = \widetilde{\BS D}(T-d,T)$ and $\widetilde{\BS D} = {\BS C}(T-d,T)$. 
This gives 
\begin{eqnarray*}
 \BS C(T-d,T)^{-1} - \widetilde{\BS D}(T-d,T) = -\widetilde{\BS
  D}(-\infty,T-d,T)
\widetilde{\BS D}(-\infty,T-d)^{-1}\widetilde{\BS D}(-\infty,T-d,T)^{\top}. 
\end{eqnarray*}
Thus block-wise for all $1\leq t,\tau\leq d$ we have
\begin{eqnarray*}
&&[\BS C(T-d,T)^{-1} - \widetilde{\BS D}(T-d,T)]_{T-t,T-\tau} \\
&=& -
[(e_{T-t} \otimes
       I_p)^\top\widetilde{\BS D}(-\infty,T-d,T)]\widetilde{\BS D}(-\infty,T)^{-1}
[(e_{T-\tau} \otimes I_p)^\top\widetilde{\BS D}(-\infty,T-d,T)]^\top.
\end{eqnarray*}
Using the above we obtain the bound
\begin{eqnarray}
\label{eq:BCS}
&&\|[\BS C(T-d,T)^{-1} - \widetilde{\BS D}(T-d,T)]_{T-t,T-\tau}\|_{2} \nonumber\\
&\leq& \lambda_{\sup}\|(e_{T-t} \otimes
       I_p)^\top\widetilde{\BS D}(-\infty,T-d,T)\|_2
\| (e_{T-\tau} \otimes I_p)^\top\widetilde{\BS D}(-\infty,T-d,T)\|_2.
\end{eqnarray}
Next we obtain a bound for the matrix rows
$(e_{T-t} \otimes
       I_p)^\top\widetilde{\BS D}(-\infty,T-d,T) = (\widetilde{\BS D}(-\infty,T-d,T)_{T-t,\ell};\ell
<T)$. 
By applying Lemma \ref{lemma:operatornorm} 
and using Theorem \ref{thm.1} we have 
\begin{align*}
    \| (e_{T-t} \otimes
       I_p)^\top\widetilde{\BS D}(-\infty,T-d,T)\|_2 
    &\leq (\sum_{\ell=-\infty}^{T-d-1} \| \widetilde{\BS D}(-\infty,T-d,T)_{T-t,\ell}\|_2^2)^{1/2}  \\
 &\leq \mathcal{K} (\sum_{\ell=-\infty}^{T-d-1} \zeta(T-t-\ell)^{2(\kappa-1)})^{1/2}\leq \mathcal{K} \zeta(d-t)^{\kappa-3/2}.
\end{align*}
Substituting the above into (\ref{eq:BCS}) 
for all $1\leq  t,\tau \leq d$ we have 
\begin{eqnarray}
\label{eq:CT2}
\|[\BS C(T-d,T)^{-1} - \widetilde{\BS D}(T-d,T)]_{T-t,T-\tau}\|_{2} 
&\leq& \mathcal{K}\zeta(d-t)^{\kappa-3/2}\zeta(d-\tau)^{\kappa-3/2}.
\end{eqnarray}

We now return to the VAR coefficients. Using the block inverse operator
identity in (\ref{eq:inverseblock1}) it can be shown that 
$1\leq j \leq d$
\begin{eqnarray*}
\Phi_{T,d,j}-\Phi_{T,j} =  -[\BS C(T-d,T)^{-1}]_{T,T}^{-1}[\BS
  C(T-d,T)^{-1}]_{T,T-j}+
[\widetilde{\BS D}(T-d,T)]_{T,T}^{-1}[\widetilde{\BS D}(T-d,T)]_{T,T-j},
\end{eqnarray*}
(the bottom rows of $\BS C(T-d,T)^{-1}$ and
$\widetilde{\BS D}(T-d,T)$ respectively). Using the above and 
(\ref{eq:CT2}) we will prove (\ref{eq:baxter1}). Setting $t=0$ and
$\tau=j$ in (\ref{eq:CT2}) gives 
\begin{eqnarray*}
\|\Phi_{T,d,j}-\Phi_{T,j}\|_{2} &\leq&  \lambda_{\sup}\| [\BS
                                       C(T-d,T)^{-1} - \widetilde{\BS
                                       D}(T-d,T)]_{T,T-j}\|_{2}\\
&& + \lambda_{\sup}\| [\BS
                                       [\BS C(T-d,T)^{-1}]_{T,T}^{-1} - \widetilde{\BS
                                       D}(T-d,T)]_{T,T}^{-1}\|_{2}\\
&\leq &\mathcal{K}\zeta(d)^{\kappa-3/2}\zeta(d-j)^{\kappa-3/2}.
\end{eqnarray*}
This proves (\ref{eq:baxter1}). Using (\ref{eq:baxter1}) we
immediately obtain (\ref{eq:baxter2}). 
\end{proof}
Note that projection methods can also be used to prove the above result (and the same bound obtained). In this case the proof would be similar to that given in the proof of Theorem 3.2 in
\cite{p:kre-17} (in the context of spatially stationary processes).

\section{Proofs of results in Section~\ref{sec:LS}}

\subsection{Proofs of results in Section~\ref{sec:LS0}}\label{sec:proof2}

The following lemma is used in the proof of Theorem \ref{lem.Cu.pos.def}.
\begin{lemma}\label{lem.Cun.approx}
Suppose Assumption \ref{assum:LS} holds and let $ G_{u,M}(\omega)$, 
$G_{u,M}^{(N)}(\omega)$ and $ G_{u}(\omega)$ be defined as in
 (\ref{eq:GuMN0}),  (\ref{eq:GuMN}) and (\ref{eq:GuMN2}) respectively. Then 
\begin{eqnarray}
\label{eq:Cun.approx1}
\sup_{\omega}\|G_{u,M}(\omega) - G_{u,M}^{(N)}(\omega)\|_2 \leq \mathcal{K}\frac{M}{N}
\end{eqnarray}
and 
\begin{eqnarray}
\label{eq:Cun.approx2}
\sup_{\omega}\|G_{u}(\omega) - G_{u,M}(\omega)\|_2 \leq  \mathcal{K}\left(\frac{1}{M}+\frac{1}{M^{\kappa-1}}\right)
\end{eqnarray}
where $\mathcal{K}$ is a constant that only depends on $K$ and
$\kappa$. 
\end{lemma}
\begin{proof}
Under Assumption \ref{assum:LS}(iii) we have 
\begin{eqnarray*}
\|G_{u,M}^{(N)}(\omega)   - G_{u,M}^{}(\omega)\|_{2} 
&\leq& \frac{1}{M}\sum_{t,\tau=T_{u,N}-M/2+1}^{T_{u,N}+M/2}
\|C_{t,\tau}^{(N)} - C_{t-\tau}(u)\|_{2} \\
&\leq& \frac{1}{M}\sum_{t,\tau=T_{u,N}-M/2+1}^{T_{u,N}+M/2}
\left(\frac{1}{N\gu{t-\tau}^{\kappa-1}} +
       \frac{|(T_{u,N}-t)|}{N\gu{t-\tau}^{\kappa}} \right) 
\leq \mathcal{K}\frac{M}{N},
\end{eqnarray*}
this proves (\ref{eq:Cun.approx1}). To prove (\ref{eq:Cun.approx2}) we
use that
\begin{eqnarray*}
G_u(\omega)&=& G_{u,M}^{}(\omega) + \frac{1}{M}\sum_{|r|\leq
  M/2}|r|C_{r}(u)\exp(ir\omega) +   \sum_{|r|>M/2} C_r(u)
  \exp(ir\omega).
\end{eqnarray*}
Under Assumption \ref{assum:LS}(iii) we have $\| C_r(u)\|_{2}\leq
K/\gu{r}^{\kappa}$ (where $\kappa >2$), thus
 \begin{eqnarray*}
\|G_u(\omega)- G_{u,M}^{}(\omega)\|_2 \leq
\frac{1}{M}\sum_{|r|\leq
  M/2}|r|\|C_{r}(u)\|_{2} +   \sum_{|r|>M/2}\|C_r(u)\|_{2} \leq \mathcal{K}\left(\frac{1}{M}+\frac{1}{M^{\kappa-1}}\right).
\end{eqnarray*}
Thus proving the result. 
\end{proof}

We are now equiped to prove Theorem~\ref{lem.Cu.pos.def}.

\begin{proof}[Proof of Theorem~\ref{lem.Cu.pos.def}]
Our aim is to show that the $\|\cdot\|_2$-norm of the matrix function
\begin{eqnarray}
\label{eq:GuMN2}
G_{u}(\omega) = \sum_{r\in \mathbb{Z}}C_{r}(u)\exp(ir\omega)
\end{eqnarray}
is bounded above and below by the $\lambda_{\sup}$ and
$\lambda_{\inf}$ respectively (for all $\omega$). Since ${\BS C}(u)$ is a (block)
Toeplitz matrix then by Toeplitz theorem (see \cite{p:toe-11} and \cite{b:bot-00}, Theorem 1.1)
this would immediately prove that the eigenvalues of $\BS C(u)$ are bounded above and below by $\lambda_{\sup}$ and $\lambda_{\inf}$ (thus proving the result). 

For a given $u\in \mathbb{R}$ and $N\in \N$ we define the integer $T_{u,N}$ as 
$T_{u,N} = \lfloor uN \rfloor$ (where $\lfloor x \rfloor$ denotes
the largest integer smaller than $x$). Let $M \in 2 \N$ and define an $M\times M$-dimensional submatrix of ${\BS C}^{(M)}$
that is centred about $T_{u,N}$
\begin{eqnarray*}
{\BS  C}^{(N)}_{u,M} &:=&(C_{T_{u,N}+s_1,T_{u,N}+s_2}^{(N)})_{s_1,s_2=-M/2+1,\dots,M/2}=:
(I_{T_{u,N}}^{(M)} \otimes I_p)^\top \BS C^{(N)} (I_{T_{u,N}}^{(M)} \otimes
I_p). 
\end{eqnarray*}
We show below that if $M$ is sufficiently small, then ${\BS
  C}^{(N)}_{u,M}$ is an approximation of the $M\times M$-dimensional
submatrix of ${\BS C}(u)$
\begin{eqnarray*}
\BS C_M(u) &:=& (C_{s_1-s_2}(u))_{s_1,s_2=-M/2+1,\dots,M/2}=:
(I_{(u)}^{(M)}\otimes I_p)^\top \BS C(u) (I_{(u)}^{(M)} \otimes I_{p}).
\end{eqnarray*}
We start by obtaining a finite approximation of 
$G_{u}(\omega)$ in terms of $\BS C_M(u)$. Let
\begin{eqnarray}
\label{eq:GuMN0}
G_{u,M}^{}(\omega) =
\frac{1}{M}\sum_{t,\tau = T_{u,N}-M/2+1}^{T_{u,N}+M/2}C_{t-\tau}(u)
\exp(i(t-\tau)\omega) = (x_\omega\otimes I_p)^* {\BS C}_{M}(u) (x_\omega\otimes I_p),
\end{eqnarray}
where $x_\omega=1/\sqrt{M}(\exp(-it\omega))_{t=T_{u,N}-M/2+1,\dots,T_{u,N}+M/2}$.
Using ${\BS C}^{(N)}_{u,M}$ for each $M\in 2\N$ and $\omega \in [0,2\pi]$ we define the quantity 
\begin{eqnarray}
\label{eq:GuMN}
G_{u,M}^{(N)}(\omega) =
\frac{1}{M}\sum_{t,\tau = T_{u,N}-M/2+1}^{T_{u,N}+M/2}C_{t,\tau}^{(N)}
\exp(i(t-\tau)\omega) = (x_\omega\otimes I_p)^* \BS C^{(N)}_{u,M} (x_\omega \otimes I_p).
\end{eqnarray}
Since $\BS C^{(M)}_{u,N}$
is a finite dimensional submatrix of $\BS C^{(N)}$, for $N>N_0$, 
the eigenvalues of $\BS C^{(N)}_{u,M}$ are bounded above and below by 
$\lambda_{\inf}$ and $\lambda_{\sup}$ respectively. Then, since $\|x_\omega\|_2=1$ we have
\begin{eqnarray}
\label{eq:eigenbound}
\lambda_{\inf}\leq  \|G_{u,M}^{(N)}(\omega)\|_{2} \leq \lambda_{\sup}
  \textrm{ for all $N$,  $M$ and $\omega$.}
\end{eqnarray}
By using Lemma \ref{lem.Cun.approx}, equation
(\ref{eq:Cun.approx1}) we have 
\begin{eqnarray}
\label{eq:eigenbound2}
\sup_{\omega}\|G_{u,M}(\omega) - G_{u,M}^{(N)}(\omega)\|_2 \leq \mathcal{K}\frac{M}{N},
\end{eqnarray}
where $\mathcal{K}$ is a generic constant that depends only on $K$ and
$\kappa$.
The above immediately implies 
$\lambda_{\inf} - \mathcal{K}M/N\leq \|G_{u,M}(\omega)\|_2 \leq \lambda_{\sup} +\mathcal{K}M/N$.
Finally we return to $G_{u}(\omega)$. Using Lemma
\ref{lem.Cun.approx}, equation (\ref{eq:Cun.approx2}) we have
$\sup_{\omega}\|G_{u}(\omega) - G_{u,M}(\omega)\|_2 \leq  \mathcal{K}/M$.
By using this and (\ref{eq:eigenbound2}) we have 
\begin{eqnarray*}
\|G_{u}(\omega)\|_{2} = \|G_{u,M}^{(N)}(\omega)\|_2 +
  O\left( \frac{M}{N}+\frac{1}{M}\right).
\end{eqnarray*}
Finally, we set $M=2\lfloor N^{1/2}\rfloor$ and substitute it into the above, this
together with (\ref{eq:eigenbound}) gives
\begin{eqnarray*}
\lambda_{\inf} - \frac{\mathcal{K}}{N^{1/2}}\leq 
\|G_{u}(\omega)\|_{2} \leq \lambda_{\sup} + \frac{\mathcal{K}}{N^{1/2}}.
\end{eqnarray*}
As this holds for all $N>N_0$ we have that for any $\eps>0$
$\lambda_{\inf}-\eps\leq \|G_{u}(\omega)\|_{2}\leq
\lambda_{\sup}+\eps$. Thus leading to the required result. 
\end{proof}

\begin{proof}[Proof of Theorem~\ref{lemma:spectral}]
We start by giving a short overview of the proof.   
To show that Assumption \ref{assum:LS}(i) holds 
(a uniform bound on the eigenvalues of $\BS C^{(N)}$) for a sufficiently large $N$, 
we first replace the infinite dimensional matrix
$\BS C^{(N)}$ with an infinite dimensional banded matrix $\BS C^{(N)}_{M}$
(where we obtain a bound for $\|\BS C^{(N)} - \BS C^{(N)}_{M}\|_{2}$). 
The central part of the proof is to obtain a bound for the eigenvalues
of $\BS C^{(N)}_{M}$ (that is uniform over a sufficiently large
$N$). The key observation is that the banded matrix
embeds an infinite number of overlapping $(M+1)\times (M+1)$-dimensional block matrices, 
where each block matrix can be approximated by an $(M+1)\times
(M+1)$-dimensional block matrix whose entries consist of a stationary autocovariance.
We will show that a lower and upper bound for the eigenvalues of each
stationary approximation block matrix is
$\gamma_{\inf}$ and $\gamma_{\sup}$ respectively.  This yields a bound for the eigenvalues of each
$(M+1)\times (M+1)$-dimensional block matrices in $\BS C^{(N)}_{M}$.
Finally, motivated by the proof of Proposition 2.9 in
\cite{p:din-zho-21}, we show that the eigenvalues of the banded matrix
$\BS C^{(N)}_{M}$ can be bounded by the eigenvalues of
``overlapping'' $(M+1)\times (M+1)$-dimensional 
block matrices. This
will prove the result.

We start by defining the infinite dimensional (block) banded matrix $\BS C_{M}^{(N)}$ where for all $t,\tau\in \mathbb{Z}$ the entries are defined by
$
[\BS C^{(N)}_{M}]_{t,\tau} = \ind(|t-\tau|\leq M)C_{t,\tau}.
$
Without loss of generality we
assume that $M=2m$ where $m \in \N$.
Using Lemma~\ref{lemma:propB} we have  $\|\BS C^{(N)} - \BS C_{M}^{(N)}\|_{2}\leq
\mathcal{K}M^{-\kappa+1}$. 
Our aim is to obtain bounds for $x^\top\BS C_{M}^{(N)}x$ where $x = (\ldots,x_{-1},x_{0},x_{1},\ldots)^{\top} \in \ell_{2,p}, x_l \in \R^p $
and $\|x\|_2=1$. To do this we define the $(M+1)p$-dimensional shifting subsequence 
$x_{s-m,s+m} = (x_{s-m},\ldots,x_{s+m})^{\top}$ and the $(M+1)p\times
(M+1)p$ dimensional (block) matrix 
\begin{eqnarray*}
\BS C^{(N)}(s-m,s+m) = (C_{t,\tau}^{(N)};s-m\leq t,\tau \leq s+m).
\end{eqnarray*}
It can be shown that for $|t-\tau|\leq M+1$ the entries of
$\BS C^{(N)}_{M}$ can be written in terms of $\BS C^{(N)}(s-m,s+m) $
\begin{eqnarray*}
  [C^{(N)}_{M}]_{t,t+r} =
  \left\{
  \begin{array}{cc}
    \frac{1}{M-|r|}\sum_{s=0}^{2m-r} C^{(N)}(t-s,t+2m-s)_{(s+1,s+1+r-2m)} &
                                                                  r\geq
                                                                  0 \\
      \frac{1}{M-|r|}\sum_{s=0}^{2m-|r|}C^{(N)}(t-2m+s,t+s)_{(s+1+r-2m,s+1)} &
                                                                  r< 0\\
\end{array}
  \right.
\end{eqnarray*}
For each $u\in \mathbb{Z}$ we define the stationary approximation matrix $\BS C(s-m,s+m;u)$
\begin{eqnarray*}
\BS C(s-m,s+m;u) = (C_{t-\tau}(u);s-m\leq t,\tau \leq s+m).
\end{eqnarray*}
Under Assumption \ref{assum:LS}(iii) and using
Lemma~\ref{lem.l2.infty} we have 
\begin{eqnarray}
\label{eq:CNsm}
&& \|\BS C^{(N)}(s-m,s+m)-\BS C(s-m,s+m;s/N)\|_2 \nonumber\\
&\leq& 
\sup_{t\in (s-m,s+m)}\sum_{\tau=s-m}^{s+m}\|C_{t,\tau}-C_{t-\tau}(s/N)\|_{2} \leq
\mathcal{K}\frac{m}{N},
\end{eqnarray}
where $\mathcal{K}$ is a generic constant that holds for all $N$ and
$s$. The condition
$$0<\gamma_{\inf} \leq \inf_{u}\inf_{\omega}\lambda_{\inf}(f(\omega;u))
\leq \sup_{u}\sup_{\omega}\lambda_{\sup}(f(\omega;u))\leq
\gamma_{\sup}<\infty$$
implies (see, among others, \cite[Proposition
2.3]{basu2015}) that for all $u\in \mathbb{R}$
\begin{eqnarray*}
&&  \gamma_{\inf} \leq \inf_{\omega}\lambda_{\inf}[f(\omega;u)] \leq 
  \lambda_{\inf}[\BS C(s-m,s+m;u)]\\
 &&\qquad   \leq \lambda_{\sup}[\BS C(s-m,s+m;u)] \leq
 \sup_{\omega}\lambda_{\sup}[f(\omega;u)] \leq \gamma_{\sup}.
  \end{eqnarray*}
Therefore by using (\ref{eq:CNsm}) and the above we have 
\begin{eqnarray}
\left(\gamma_{\inf} -
  \mathcal{K}\frac{m}{N}\right)\|x_{s-m,s+m}\|_{2}&\leq& 
x_{s-m,s+m}^{\top}\BS C^{(N)}(s-m,s+m) x_{s-m,s+m}  \nonumber\\
&\leq&
  \left(\gamma_{\sup} +                                                     \mathcal{K}\frac{m}{N}\right)\|x_{s-m,s+m}\|_{2}.
         \label{eq:CMM2}                                           
\end{eqnarray}
This gives a bound for each block. 
Next we obtain a bound between
\begin{eqnarray}
\label{eq:BSC1}
 x^{\top}\BS C^{(N)}_{M}x = \sum_{\ell\in \mathbb{Z}}\sum_{r=-M}^{M}x_{\ell}^{\top}C_{\ell,\ell+r}x_{\ell+r}
\end{eqnarray}
 with the overlapping block matrix inner-product
\begin{eqnarray*}
X_{M}^{\top}\BS O_MX_{M} := \frac{1}{M+1}\sum_{s\in
  \mathbb{Z}}x_{s-m,s+m}^{\top}\BS C^{(N)}(s-m,s+m)x_{s-m,s+m}.
\end{eqnarray*}
Note we have not formally defined $X_{M}$ or $O_{M}$ but have simply
set it to equal the above.
Basic algebra gives
\begin{eqnarray}
\label{eq:BSC2}
X_{M}^{\top}\BS O_{M}X_{M} = \sum_{\ell\in \mathbb{Z}}\sum_{r=-M}^{M}\left(\frac{M+1-|r|}{M+1}\right)x_{\ell}^{\top}C_{\ell,\ell+r}x_{\ell+r}.
\end{eqnarray}
Using (\ref{eq:BSC1}) and (\ref{eq:BSC2}) we have 
\begin{eqnarray*}
x^{\top}\BS C^{(N)}_{M}x - X_{M}^{\top}\BS O_{M}X_{M}= \frac{1}{M+1}\sum_{\ell\in \mathbb{Z}}\sum_{r=-M}^{M}|r|x_{\ell}^{\top}C_{\ell,\ell+r}x_{\ell+r}.
\end{eqnarray*}
Hence under Assumption \ref{assum:LS}(ii) we have
\begin{eqnarray}
\label{eq:CMM3}
\left\|x^{\top}\BS C^{(N)}_{M}x - X_{M}^{\top}\BS O_{M}X_{M}\right\|_{2}
  &\leq &
\frac{1}{M+1}\sum_{\ell\in
         \mathbb{Z}}\sum_{r=-m}^{m}\frac{|r|}{\gu{r}^{\kappa}}
         \|x_{\ell}\|_{2}\|x_{\ell+r}\|_2 \nonumber\\
&\leq& \frac{2}{M+1}\left(\sum_{r=1}^{\infty}\frac{1}{\gu{r}^{\kappa-1}}\right)
        \sum_{\ell\in \mathbb{Z}} \|x_{\ell}\|_{2}^{2} = \frac{2}{M+1}\left(\sum_{r=1}^{\infty}\frac{1}{\gu{r}^{\kappa-1}}\right),
\end{eqnarray}
where the last line follows because $\|x\|_2=\sum_{\ell\in \mathbb{Z}}
\|x_{\ell}\|_{2}^{2}=1$.
Finally, we obtain an upper and lower bound for $X_{M}^{\top}\BS O_{M}X_{M}$.
We use (\ref{eq:CMM2}) to give 
\begin{eqnarray*}
\frac{(\gamma_{\inf}-\mathcal{K} m/N)}{M}\sum_{s\in
  \mathbb{Z}}\|x_{s-m,s+m}\|_{2}^{2}\leq X_{M}^{\top}\BS O_{M}X_{M}\leq \frac{(\gamma_{\sup}+\mathcal{K} m/N)}{M}\sum_{s\in
  \mathbb{Z}}\|x_{s-m,s+m}\|_{2}^{2}.
\end{eqnarray*}
Using that $\sum_{s\in
  \mathbb{Z}}\|x_{s-m,s+m}\|_{2}^{2} = (M+1)\|x\|_{2}^{2}=(M+1)$ we have 
\begin{eqnarray*}
\gamma_{\inf}-\mathcal{K} m/N\leq X_{M}^{\top}\BS O_{M}X_{M}\leq \gamma_{\sup}+\mathcal{K} m/N.
\end{eqnarray*}
Hence by using (\ref{eq:CMM3}), $\|\BS C^{(N)} - \BS C_{M}^{(N)}\|_{2}\leq
\mathcal{K}M^{-\kappa+1}$ and setting  $m=\lfloor N^{1/\kappa}\rfloor$ we have 
\begin{eqnarray*}
\gamma_{\inf}-\mathcal{K}N^{-1+1/\kappa}\leq x^{\top}\BS C^{(N)}x  \leq \gamma_{\sup}+\mathcal{K}N^{-1+1/\kappa},
\end{eqnarray*}
where $\mathcal{K}$ is generic constant that does not depend on $N$ or
$M$. Thus for a sufficiently large $N$ we have the result. 
\end{proof}

\subsection{Proofs for spectral-norm physical dependence systems}\label{sec:physical}


In order to prove Theorem \ref{theorem:physical} we require the
following lemma.

\begin{lemma}\label{lemma:covariancephysicaldependence}
Suppose $\{V_{t}\}_{t}$ and $\{U_{t}\}_{t}$ are zero
mean multivariate time series of dimension $p$ that have the
causal representation $V_{t} = v_{t}(\mathcal{F}_{t})$ and
$U_{t}=u_{t}(\mathcal{F}_{t})$ where $\mathcal{F}_{t} =
(\varepsilon_{t},\varepsilon_{t-1},\ldots)$
and $\{\varepsilon_{t}\}$ are iid
random vectors. Let $\{\widetilde{\varepsilon}_{t}\}_{t}$ are iid random vectors
with the same distribution as $\{\varepsilon_{t}\}_{t}$ but
independent of them and set
$\mathcal{F}_{t|\{t-j\}} =
(\varepsilon_{t},\varepsilon_{t-1},\ldots,
\varepsilon_{t-j+1},\widetilde{\varepsilon}_{t-j},\varepsilon_{t-j-1},\ldots)$. Let $U_{t|\{t-j\}} =
u_{t}(\mathcal{F}_{t|\{t-j\}})$ and $V_{t|\{t-j\}} = v_{t}(\mathcal{F}_{t|\{t-j\}})$.
We assume that $\{V_{t}\}$ and $\{U_{t}\}$
satisfy the spectral-norm physical dependence conditions
\begin{eqnarray*}
\sup_{t} \|\var((U_{t} - U_{t|\{t-j\}})\|_{2} \leq A\delta_{j}
  \textrm{ and }\sup_{t} [\var(V_{t} - V_{t|\{t-j\}})\|_{2} \leq B\delta_{j}
\end{eqnarray*}
where $\delta_{j} = \gu{j}^{-\kappa}$ and $\kappa >1$. 
Then 
\begin{eqnarray*}
\|\cov(U_{t},V_{\tau})\|_2 \leq AB\delta_{|t-\tau|}\sum_{j=0}^{\infty}\delta_{j}.
\end{eqnarray*}
\end{lemma}
\begin{proof}
To prove the result we write both $V_{t}$ and $U_{t}$ as the sum of
martingale differences and represent $\cov(U_{t},V_{\tau})$ as the
covariance of the martingale difference. This representation together
with the physical dependence condition will prove the result. The
details are below.

With a small abuse of notation we define the sigma algebra
$\mathcal{F}_{t} =
\sigma(\varepsilon_{t},\varepsilon_{t-1},\ldots)$. 
Since $U_{t},V_{t}\in \mathcal{F}_{t}$, almost surely we
can represent $U_{t}$ and $V_{t}$ as the infinite sum
\begin{eqnarray*}
U_{t} = \sum_{j=0}^{\infty}D_{U}(t,j) \textrm{ and }
V_{t} = \sum_{j=0}^{\infty}D_{V}(t,j), 
\end{eqnarray*}
where  $D_{U}(t,j) = \Ex(U_{t}|\mathcal{F}_{t-j}) -
  \Ex(U_{t}|\mathcal{F}_{t-j-1})$ and 
$D_{V}(t,j) = \Ex(V_{t}|\mathcal{F}_{t-j}) -
  \Ex(V_{t}|\mathcal{F}_{t-j-1})$.
Without loss of generality we assume that $t-j_1<\tau-j_2$, then by
iterated expectations
\begin{eqnarray}
\label{eq:MDiff}
\Ex\left[ D_{U}(t,j_1)D_{V}(\tau,j_2)^\top\right] &=& 
\Ex\left[
  D_{U}(t,j_1)\Ex[D_{V}(\tau,j_2)^\top|\mathcal{F}_{t-j_1}]\right]\nonumber\\
&=& \Ex\left[
  D_{U}(t,j_1)\Ex[(\Ex(V_{t}|\mathcal{F}_{t-j}) -
  \Ex(V_{t}|\mathcal{F}_{t-j-1}))^\top|\mathcal{F}_{t-j_1}]\right]
  = 0, \quad 
\end{eqnarray}
where the above holds because 
 for any $t-j_1 < \tau-i$,
$\Ex[\Ex(V_{\tau}|\mathcal{F}_{\tau-i})|\mathcal{F}_{t-j_1}] =
\Ex[V_{\tau}|\mathcal{F}_{t-j_1}]$. By  a similar argument,
$\Ex\left[ D_{U}(t,j_1)D_{V}(\tau,j_2)^\top\right] = 0$
for  $t-j_1>\tau-j_2$.

We use (\ref{eq:MDiff}) to write $\cov[V_{t},U_{\tau}] $ as the
product of martingale differences. Using (\ref{eq:MDiff}) and assuming
$t<\tau$ we have 
\begin{eqnarray*}
\cov[V_{t},U_{\tau}]   = \sum_{j=\tau-t}^{\infty}\Ex[D_{V}(t,j)D_{U}(\tau,\tau-t+j)^\top].
\end{eqnarray*}
Applying Lemma~\ref{lemma:CauchySchwarzMV} to the above gives 
\begin{eqnarray}
\label{eq:covVW}
\|\cov[V_{t},U_{\tau}]\|_2   \leq
  \sum_{j=\tau-t}^{\infty}\|\var[D_{V}(t,j)]\|_2^{1/2}
\|\var[D_{U}(\tau,\tau-t+j)]\|_2^{1/2}.
\end{eqnarray}
Finally to bound the above expectations we use the physical dependence
condition and the observation
\begin{eqnarray*}
D_{V}(t,j) = \Ex[V_{t}-V_{t|\{t-j\}}|\mathcal{F}_{t-j}].
\end{eqnarray*}
Thus by using the law of total variance we have 
\begin{eqnarray*}
\|\var[D_{V}(t,j)^{}]\|_{2}^{1/2} \leq
  \|\var[(V_{t}-V_{t|\{t-j\}})^{}]\|_{2}^{1/2} \leq B\delta_{j}.
\end{eqnarray*}
By the same argument we have
$\|\var[D_{U}(\tau,\tau-t+j)^{2}]\|_2^{1/2}\leq
A\delta_{\tau-t+j}$. Substituting these bounds into (\ref{eq:covVW})
and using that $\delta_{j} = \gu{j}^{-\kappa}$ we have
 \begin{eqnarray*}
|\cov[V_{t},U_{\tau}]|   \leq AB\delta_{\tau-t}\sum_{j=0}^{\infty}\delta_{j}.
\end{eqnarray*}
A similar bound holds for the case $t>\tau$. 
Thus proving the result. 
\end{proof}

\begin{proof}[Proof of Theorem~\ref{theorem:physical}] 
We first show that condition (A) implies that Assumption
\ref{assum:LS}(ii,iii) hold.

By using Lemma \ref{lemma:covariancephysicaldependence} 
(with $V_{t} = W_{t} = X_{t,N}$)
and the spectral-norm physical dependence condition on $\{X_{t,N}\}$ it immediately follows that
Assumption \ref{assum:LS}(ii) holds. 

We now show that $C_{r}(u)=\cov[X_{0}(u),X_{r}(u)]$ satisfies
Assumption \ref{assum:LS}(iii) parts (a),(b) and (c).
Assumption \ref{assum:LS}(iii) part (a) holds by definition of
$C_{r}(u)$. Assumption \ref{assum:LS}(iii) part (b) follows from Lemma
\ref{lemma:covariancephysicaldependence} and the spectral-norm physical dependence
condition on $\{X_{t}(u)\}_{t}$. 
To show that Assumption \ref{assum:LS}(iii) part (c) holds we treat
the case $|u-v|\leq 1$ and $|u-v|>1$ separately. For $|u-v|>1$ by using (b) we
have 
\begin{eqnarray*}
\|\cov(X_{t}(u),X_{\tau}(u))  - \cov(X_{t}(v),X_{\tau}(v)) \|_2 &\leq&
  \|\cov(X_{t}(u),X_{\tau}(u))\|_2\!+\!
\|\cov(X_{t}(v),X_{\tau}(v))\|_2 \\
&\leq& \frac{2K}{\gu{t-\tau}^{\kappa}} \leq \frac{2K|u-v|}{\gu{t-\tau}^{\kappa}}, 
\end{eqnarray*}
thus (c) holds for $|u-v|>1$. For the case $|u-v|\leq 1$,
we first use condition i) which states that 
$\tilde X_t^{u,v}=X_{t}(u)-X_{t}(v)$ satisfies the physical dependence condition
$\|\var(\tilde X_t^{u,v}-\tilde X_{t|\{t-j\}}^{u,v})\|_2 \leq |u-v| K \delta_j.$
Thus, by using the expansion $X_{\tau}(v) = X_{\tau}(u) + \tilde X_{\tau}^{u,v}$ and
applying Lemma \ref{lemma:covariancephysicaldependence}
we have 
\begin{eqnarray*}
\|C_{t-\tau}(u)-C_{t-\tau}(v)\|_2 &\leq& \!\bigg[\|\cov[X_{t}(v),
 \tilde X_\tau^{u,v}]\|_2\!+\! \|\cov[
    \tilde X_t^{u,v}, X_{\tau}(v)]\|_2\! +\!\|\cov[\tilde X_\tau^{u,v},\tilde X_t^{u,v}\|_2\bigg]\\
    & \leq& \frac{K|u-v|}{\gu{t-\tau}^{\kappa}}+\frac{K|u-v|}{\gu{t-\tau}^{\kappa}}+\frac{K|u-v|^2}{\gu{t-\tau}^{\kappa}}\leq  \frac{3K|u-v|}{\gu{t-\tau}^{\kappa}},
\end{eqnarray*}
where the last inequality follows due to the condition $|u-v|\leq 1$.
This proves that Assumption \ref{assum:LS}(ii)(c) holds. 

Finally to prove that (\ref{eq:LScovariance2}) holds,  we use a similar
technique as above. We focus on the case $|t-\tau|>N$ and
$|t-\tau|\leq N$ separately. For $|t-\tau|>N$ and by using the above
for the bounds for 
$C_{t,\tau}^{(N)}$ and $C_{t-\tau}(t/N)$ it can be shown that
$\|C_{t,\tau}^{(N)} - C_{t-\tau}(t/N)\|_2\leq 2K\gu{t-\tau}^{-\kappa}$. On
the other hand for $|t-\tau|\leq N$ we use that
\begin{eqnarray*}
X_{t,N} &=& X_{t}(t/N)+e_{t,N} \\
\textrm{ and } X_{\tau,N}  &=&
X_{\tau}(\tau/N)+e_{\tau,N} =
                               X_{\tau}(t/N)+\{X_{\tau}(\tau/N)-X_{\tau}(t/N)\}+e_{\tau,N}.
\end{eqnarray*}
Substituting the above into $\cov(X_{t,N},X_{\tau,N})$ gives 
\begin{eqnarray*}
\cov(X_{t,N},X_{\tau,N}) = \cov\left[X_{t}(t/N)+e_{t,N},X_{\tau}(t/N) + \{X_{\tau}(\tau/N)-X_{\tau}(t/N)\}+e_{\tau,N}\right].
\end{eqnarray*}
Expanding out the above covariance and using that 
$\cov[X_{t}(t/N), X_{\tau}(t/N)]=C_{t-\tau}(t/N)$ we have
\begin{eqnarray*}
&&\|C_{t,\tau}^{(N)} - C_{t-\tau}(t/N)\|_2 \\ 
&\leq& \|\cov[X_{t}(t/N),  (X_{\tau}(\tau/N)-X_{\tau}(t/N))]\|_2+
\|\cov[X_{t}(t/N), e_{\tau,N}]\|_2 \\
&& + \|\cov[e_{t,N}, X_{\tau}(t/N)]\|_2+ \|\cov[e_{t,N}, (X_{\tau}(t/N)-X_{\tau}(\tau/N))]\|_2
+\|\cov[e_{t,N},e_{\tau,N}]\|_2.
\end{eqnarray*}
Under the spectral-norm physical dependence conditions (and by using Lemma
\ref{lemma:covariancephysicaldependence}) we have
\begin{eqnarray*}
&&\|C_{t,\tau}^{(N)} - C_{t-\tau}(t/N)\|_2 \\ 
&\leq& K\left(|t-\tau|N^{-1}\gu{t-\tau}^{-\kappa} + 2\cdot
       N^{-1}\gu{t-\tau}^{-\kappa}+
\frac{|t-\tau|}{N^{2}}\gu{t-\tau}^{-\kappa}+N^{-2}\gu{t-\tau}^{-\kappa}\right)
  \\
&\leq&  K\left(2\cdot N^{-1}\gu{t-\tau}^{-\kappa+1} + 3\cdot
       N^{-1}\gu{t-\tau}^{-\kappa}\right) \\
&\leq& 5K \cdot N^{-1}\gu{t-\tau}^{-\kappa+1},
\end{eqnarray*}
where the last line is due to $|t-\tau|\leq N$. The bounds for the two cases, 
$|t-\tau|\leq N$ and $|t-\tau|>N$ show that Assumption
\ref{assum:LS}(iii) equation (\ref{eq:LScovariance2}) holds. 

Finally, under condition (B) and by applying Theorem
\ref{lemma:spectral} we have that Assumption \ref{assum:LS}(i) holds.  
\end{proof}

In order to study the properties of the stochastic recurrence equation
defined in Example \ref{exam:sre} we state a general result for the
time series $\{Y_{t}\}$ where 
\begin{eqnarray}
\label{eq:Ytdefin}
Y_{t} = G_{t}(\mathcal{F}_{t}) = \sum_{s=0}^{\infty}B_{t-s+1,t}^{}(\mathcal{F}_{t-s+1,t}^{})b_{t-s}(\varepsilon_{t-s})
\end{eqnarray}
with $\mathcal{F}_{t} = (\varepsilon_{t},\varepsilon_{t-1},\ldots)$,
$\mathcal{F}_{t-s,t}^{} =
(\varepsilon_{t},\varepsilon_{t-1},\ldots,\varepsilon_{t-s})$ and
$\{\varepsilon_{t}\}$ are iid random variables. 
Using $Y_{t}$ we define the following coupled process 
\begin{eqnarray}
\label{eq:Ytdefin2}
Y_{t|\{t-j\}} &=&  G_{t}(\mathcal{F}_{t|\{t-j\}})  = 
\sum_{s=0}^{t-j+2}B_{t-s+1,t}^{}(\mathcal{F}_{t-s+1,t}^{})b_{t-s}(\varepsilon_{t-s})\nonumber\\
&&\quad  +B_{t-j+1,t}^{}(\mathcal{F}_{t-s,t}^{})b_{t-j}(\widetilde{\varepsilon}_{t-j})+
\sum_{s=0}^{t-j+2}B_{t-s+1,t}^{}(\mathcal{F}_{t-s+1,t}^{})b_{t-s}(\varepsilon_{t-s}),
\end{eqnarray}
where $\mathcal{F}_{t|\{t-j\}} =
(\varepsilon_{t},\varepsilon_{t-1},\ldots,\varepsilon_{t-j+1},\widetilde{\varepsilon}_{t-j},\varepsilon_{t-j-1},\ldots)$ and if $s<j$
$\mathcal{F}_{t-s,t|\{t-j\}}^{} = \mathcal{F}_{t-s,t}^{}$  else
$\mathcal{F}_{t-s,t|\{t-j\}}^{}=
(\varepsilon_{t},\varepsilon_{t-1},\ldots,\widetilde{\varepsilon}_{t-j},\varepsilon_{t-j-1}\ldots,\varepsilon_{t-s})$
and $\widetilde{\varepsilon}_{t-j}$ is independent of
$\varepsilon_{t-j}$. 

To prove the following result we use that if $B_{1}$ and $B_{2}$ are conformable
independent random matrices then
\begin{eqnarray}
\label{eq:BBbop}
\|\Ex(B_{1}B_{2}B_{2}^{\top}B_{1}^{\top})\|_{2} \leq \|\Ex[B_{1}B_{1}^{\top}]\|_{2}\|\Ex[B_{2}B_{2}^{\top}]\|_{2}.
\end{eqnarray}
Taking this further, if $B_{1},\ldots,B_{K}$ are independent
conformable random matrices then 
\begin{eqnarray}
\label{eq:BBbop2}
\|\Ex(B_{1}B_{2}\ldots
  B_{K}B_{K}^{\top}\ldots B_{2}^{\top}B_{1}^{\top})\|_{2} \leq
  \prod_{i=1}^{K}
\|\Ex[B_{i}B_{i}^{\top}]\|_{2}.
\end{eqnarray}
To simplify notation in the proofs below for the random vector or matrix $X$ we let $V(X) =
\Ex[XX^{\top}]$.

\begin{lemma}\label{lemma:markovbound}
Let $Y_{t}$ and $Y_{t|\{t-j\}}$ be defined as in (\ref{eq:Ytdefin}) and (\ref{eq:Ytdefin2})
respectively. Then we have 
\begin{eqnarray*}
\|V(Y_{t} - Y_{t|\{t-j\}})\|_{2} 
&\leq& 
 4(\sum_{s =j}^{\infty}\|V[b_{t-s}(\varepsilon_{t-s}) ]\|_{2}^{1/2}\|V[B_{t-s+1,t}^{}(\mathcal{F}_{t-s+1,t}^{}) 
]\|_{2}^{1/2})^{2}
\end{eqnarray*}
\end{lemma}
\begin{proof}
Using \eqref{eq:Ytdefin} and \eqref{eq:Ytdefin2} we have
\begin{eqnarray}
\label{eq:Ytdiff}
Y_{t} - Y_{t|\{t-j\}}  &=&
  B_{t-j+1,t}^{}(\mathcal{F}_{t-j+1,t}^{})[b_{t-j}(\varepsilon_{t-j})-b_{t-j}(\widetilde{\varepsilon}_{t-j})] \nonumber\\
&&  + 
\sum_{s=j+1}^{\infty}[B_{t-s+1,t}^{}(\mathcal{F}_{t-s+1,t}^{})-B_{t-s+1,t}^{}(\mathcal{F}_{t-s+1,t|\{t-j\}}^{})]b_{t-s}(\varepsilon_{t-s}).
\end{eqnarray}
Applying (\ref{eq:ExY}) to the above gives
\begin{eqnarray*}
&&\|V(Y_{t} - Y_{t|\{t-j\}})\|_2^{1/2}\leq
\|V(B_{t-j+1,t}^{}(\mathcal{F}_{t-j+1,t}^{})[b_{t-j}(\varepsilon_{t-j})-b_{t-j}(\widetilde{\varepsilon}_{t-j})])\|_{2}^{1/2} \nonumber\\
&&  + 
\sum_{s=j+1}^{\infty}\|V([B_{t-s+1,t}^{}(\mathcal{F}_{t-s+1,t}^{})-B_{t-s+1,t}^{}(\mathcal{F}_{t-s+1,t|\{t-j\}}^{})]b_{t-s}(\varepsilon_{t-s}))\|_{2}^{1/2}.
\end{eqnarray*}
By applying (\ref{eq:BBbop}) to each of the terms above we have
\begin{eqnarray}
&&\|V(Y_{t} - Y_{t|\{t-j\}})\|_2^{1/2}\leq 
\|V(b_{t-j}(\varepsilon_{t-j})-b_{t-j}(\widetilde{\varepsilon}_{t-j}))\|_{2}^{1/2} 
\|V[B_{t-j+1,t}^{}(\mathcal{F}_{t-j+1,t})]\|_{2}^{1/2}
\nonumber\\
&&  + 
\sum_{s=j+1}^{\infty}
\|V[B_{t-s+1,t}^{}(\mathcal{F}_{t-s+1,t}^{})-B_{t-s+1,t}^{}(\mathcal{F}_{t-s+1,t|\{t-j\}}^{})]\|_{2}^{1/2}
\|V(b_{t-s}(\varepsilon_{t-s}))\|_{2}^{1/2}. \label{eq:VYdiff}
\end{eqnarray}
We now bound the terms inside the above sum. By using  Lemma
\ref{lemma:CauchySchwarzMV} we have 
\begin{eqnarray*}
&&\|V(b_{t-j}(\varepsilon_{t-j})-b_{t-j}(\widetilde{\varepsilon}_{t-j}))\|_{2}
  \leq
  \|V(b_{t-j}(\varepsilon_{t-j}))\|_{2}+\|V(b_{t-j}(\widetilde{\varepsilon}_{t-j}))\|_{2}\\
  &&+
2\|\Ex[b_{t-j}(\varepsilon_{t-j})
  b_{t-j}(\widetilde{\varepsilon}_{t-j})^{\top}]\|_{2}\\
 &\leq&\|V(b_{t-j}(\varepsilon_{t-j}))\|_{2}+\|V(b_{t-j}(\widetilde{\varepsilon}_{t-j}))\|_{2}+
2\|V(b_{t-j}(\varepsilon_{t-j}))\|_{2}^{1/2}\|
        V(b_{t-j}(\widetilde{\varepsilon}_{t-j}))]\|_{2}^{1/2}\\
&=&4 \|V(b_{t-j}(\varepsilon_{t-j}))\|_{2},
\end{eqnarray*}
where the last line follows from the fact that $V(b_{t-j}(\varepsilon_{t-j}))=V(b_{t-j}(\widetilde{\varepsilon}_{t-j}))$.
By a similar set of arguments we have 
\begin{eqnarray*}
\|V[B_{t-s+1,t}^{}(\mathcal{F}_{t-s+1,t}^{})-B_{t-s+1,t}^{}(\mathcal{F}_{t-s+1,t|\{t-j\}}^{})]\|_{2}
&\leq& 4\|V[B_{t-s+1,t}^{}(\mathcal{F}_{t-s+1,t}^{}) ]\|_{2}.
\end{eqnarray*}
Substituting these bounds into (\ref{eq:VYdiff}) gives the result. 
\end{proof}

We now apply the above result to the nonstationary model
\begin{eqnarray*}
X_{t,N} = A(t/N,\varepsilon_{t})X_{t-1,N} + b(t/N,\varepsilon_{t})
\end{eqnarray*}
and its stationary approximation 
\begin{eqnarray*}
X_{t}(u) = A(u,\varepsilon_{t})X_{t-1}(u) + b(u,\varepsilon_{t}).
\end{eqnarray*}
In the case $\sup_{u}\|\Ex[A(u,\varepsilon_{t})
A(u,\varepsilon_{t})^{\top}]\|_{2}<1$, then both $X_{t,N}$ and
$X_{t}(u)$ admit the causal solutions
\begin{eqnarray*}
X_{t,N} &=& g_{t,N}(\mathcal{F}_{t}) = 
  \sum_{s=0}^{\infty}\prod_{i=0}^{s-1}A((t-i)/N,\varepsilon_{t-i})
  b((t-s)/N,\varepsilon_{t-s}) \\
X_{t}(u) &=&  g(u,\mathcal{F}_{t}) =\sum_{s=0}^{\infty}\prod_{i=0}^{s-1}A(u,\varepsilon_{t-i})
  b(u,\varepsilon_{t-s}). 
\end{eqnarray*}
We will use Lemma \ref{lemma:markovbound} to prove the assertion in
Example \ref{exam:sre}. The above expansions for $X_{t,N}$ and
$X_{t}(u)$ allow us to 
apply Lemma \ref{lemma:markovbound} to obtain the physical dependence
bound. In the same spirit we require analogous expansions for
$e_{t,N}  = X_{t,N}-X_{t}(t/N)$ and
$\widetilde{X}_{t}^{v_1,v_2} = X_{t}(v_1)-X_{t}(v_1)$
\begin{eqnarray*}
e_{t,N} &=&  \sum_{s=0}^{\infty}\prod_{i=0}^{s-1}\left[A\left(\frac{(t-i)}{N},\varepsilon_{t-i}\right)
  b\left(\frac{(t-s)}{N},\varepsilon_{t-s}\right) -
  \prod_{i=0}^{s-1}A\left(\frac{t}{N},\varepsilon_{t-i}\right) b\left(\frac{t}{N},\varepsilon_{t-s}\right)\right] \\
 &=&
     \sum_{s=0}^{\infty}\sum_{k=0}^{s-1}\prod_{i=0}^{k-1}A\left(\frac{t-i}{n},\varepsilon_{t-i}\right)\left[A\left(\frac{t-k}{n},
     \varepsilon_{t-k}\right)-
A\left(\frac{t}{n},\varepsilon_{t-k}\right)\right]\\
&&\times \prod_{i=k+1}^{s-1}A\left(\frac{t}{n},\varepsilon_{t-i}\right)
b(\frac{t-s}{n},\varepsilon_{t-s})
     \\
&& +\sum_{s=0}^{\infty}\prod_{i=0}^{s-1}A(u,\varepsilon_{t-i})
   \left[b\left(\frac{t-s}{N},\varepsilon_{t-s}\right)-
   b\left(\frac{t}{n},\varepsilon_{t-s}\right)\right]. 
\end{eqnarray*}
and 
\begin{eqnarray*}
\widetilde{X}_{t}^{v_1,v_2} 
 &=&
     \sum_{s=0}^{\infty}\sum_{k=0}^{s-1}[\prod_{i=0}^{k-1}A(v_1,\varepsilon_{t-i})][A(v_1,\varepsilon_{t-k})-
A(v_2,\varepsilon_{t-k})][\prod_{i=k+1}^{s-1}A(v_2,\varepsilon_{t-i})]b(\frac{t-s}{n},\varepsilon_{t-s})
     +\\
&& \sum_{s=0}^{\infty}\prod_{i=0}^{s-1}A(v_2,\varepsilon_{t-i})  [b(v_1,\varepsilon_{t-s})-b(v_2,\varepsilon_{t-s})]. 
\end{eqnarray*}
Using the above expansion we can prove Example \ref{exam:sre}, which
is  given in the following lemma.

\begin{lemma}\label{lemma:physicaldependencemarkov}
Suppose $\sup_{u}\|\Ex[A(u,\varepsilon_{t})
A(u,\varepsilon_{t})^{\top}]\|_{2}<\rho<1$, $\sup_{u}\|\Ex[b(u,\varepsilon_{t})
b(u,\varepsilon_{t})^{\top}]\|_{2}<\infty$ and for all $v_1$ and $v_2$ 
$\|\Ex[(A(v_1,\varepsilon_{t})-A(v_2,\varepsilon_{t}))(A(v_1,\varepsilon_{t})-A(v_2,\varepsilon_{t}))^{\top}]\|_{2}\leq
K|v_1-v_2|$ and $\|\Ex[(b(v_1,\varepsilon_{t})-b(v_2,\varepsilon_{t}))(b(v_1,\varepsilon_{t})-b(v_2,\varepsilon_{t}))^{\top}]\|_{2}\leq
K|v_1-v_2|$. Then 
\begin{eqnarray*}
\sup_{N}\sup_{t} \| \var(X_{t,N} - X_{t,N|\{t-j\}})\|_2&\leq&
                                                                 K\rho^{j} \\
  \\
\sup_{u} \|\var(X_{t|\{t-j\}}(u) - X_{t|\{t-j\}}(u))\|_2&\leq&
                                                                   K\rho^{j}
  \\
\|\var(X_t^{v_1,v_2}-X_{t|\{t-j\}}^{v_1,v_2}])\|_2 &\leq& K|v_1-v_2|(\sum_{s=j}^{\infty}s\rho^{(s-1)/2})^{2} \\
\textrm{ and }\sup_{N}\sup_{t} \|\var(e_{t,N} - e_{t,N|\{t-j\}})\|_2 &\leq&
  KN^{-1}(\sum_{s=j}^{\infty}s^{3/2}\rho^{(s-1)/2})^{2}.
\end{eqnarray*}
\end{lemma}
\begin{proof}
To prove the result we obtain bounds for $\|V[b_{t-s}(\varepsilon_{t-s}) ]\|_{2}$ and $\|V[B_{t-s+1,t}^{}(\mathcal{F}_{t-s+1,t}^{}) 
]\|_{2}$ and apply Lemma \ref{lemma:markovbound}. By using
that $\{\varepsilon_{t}\}_{t}$ are iid random vectors and applying
(\ref{eq:BBbop2}) we have 
\begin{eqnarray*}
\left\|V\left(\prod_{i=0}^{s-1}A(u,\varepsilon_{t-i})\right)\right\|_{2}
  \leq \|V(A(u,\varepsilon_{0}))\|_{2}^{s}\leq \rho^{s},
\end{eqnarray*}
\begin{eqnarray*}
\left\|V\left(\prod_{i=0}^{s-1}A((t-i)/n,\varepsilon_{t-i})\right)\right\|_{2}
  \leq \prod_{i=0}^{s-1}\|V(A((t-i)/n,\varepsilon_{0}))\|_{2}
  \leq \rho^{s}.
\end{eqnarray*}
Further, by using (\ref{eq:ExY}) we have
\begin{eqnarray*}
&&\left\|V\left(\sum_{k=0}^{s-1}[\prod_{i=0}^{k-1}A(v_1,\varepsilon_{t-i})][A(v_1,\varepsilon_{t-k})-
A(v_2,\varepsilon_{t-k})][\prod_{i=k+1}^{s-1}A(v_2,\varepsilon_{t-i})]
  \right)\right\|_{2}\\
  &\leq& \left( \sum_{k=0}^{s-1}
\left\|V\left([\prod_{i=0}^{k-1}A(v_1,\varepsilon_{t-i})][A(v_1,\varepsilon_{t-k})-
A(v_2,\varepsilon_{t-k})][\prod_{i=k+1}^{s-1}A(v_2,\varepsilon_{t-i})] \right) \right\|_{2}^{1/2}
         \right)^{2}  \\
  &\leq& \left( \sum_{k=0}^{s-1}
\|V(A(v_1,\varepsilon_{t-i})\|_{2}^{k/2}\|V[A(v_1,\varepsilon_{t-k})-
A(v_2,\varepsilon_{t-k})\|_{2}^{1/2}\|V(A(v_2,\varepsilon_{t-i})])\|_{2}^{(s-k-1)/2}
         \right)^{2} \\
  &\leq& Ks^{2}\rho^{s-1}|v_1-v_2|
\end{eqnarray*}
and 
\begin{eqnarray*}
 \left\|V\left(
\sum_{k=0}^{s-1}[\prod_{i=0}^{k-1}A(\frac{t-i}{n},\varepsilon_{t-i})][A(\frac{t-k}{n},\varepsilon_{t-k})-
A(\frac{t}{n},\varepsilon_{t-k})][\prod_{i=k+1}^{s-1}A(\frac{t}{n},\varepsilon_{t-i})]\right)
     \right\|_{2}\leq \frac{Ks^{3}\rho^{s-1}}{n}.
\end{eqnarray*}
We use the above bounds obtain the result. 

To bound $\|\var(X_{t,N} - X_{t,N|\{t-j\}})\|_{2}$ we set
\begin{eqnarray*}
B_{t-s+1,t}^{}(\mathcal{F}_{t-s+1,t}^{}) =
  \prod_{i=0}^{s-1}A((t-i)/n,\varepsilon_{t-i}) \textrm{ and }
  b_{t-s}(\varepsilon_{t-s}) = b((t-s)/n,\varepsilon_{t-s}).
\end{eqnarray*}
Now by applying Lemma \ref{lemma:markovbound} we have 
\begin{eqnarray*}
\sup_{N}\sup_{t} \|\var[X_{t,N} - X_{t,N|\{t-j\}})\|_{2}&\leq&
                                                                 (K\sum_{s=j}\rho^{s/2})^{1/2}
  \leq K\rho^{j}.
\end{eqnarray*}
By using similar arguments we obtain the remaining bounds. 
\end{proof}

\subsection{Proof of results in Sections~\ref{sec:LS1} and \ref{sec:LS2}}\label{sec:proof3}

\begin{proof}[Proof of Theorem~\ref{thm.D.smooth}]
We begin with the proof of \eqref{lem.D.cont}. Note that $\BS C^{-1}=\BS D$. Using the classical matrix inverse
expansion we have 
\begin{eqnarray}
\label{eq:DuDv}
\BS D(u)-\BS D(v) &=& \BS C(u)^{-1}-\BS C(v)^{-1} = \BS C(u)^{-1}[\BS
                      C(v)-\BS C(u)]\BS C(v)^{-1} \nonumber\\
    &=& \BS D(u)[\BS C(v)-\BS C(u)]\BS D(v).
\end{eqnarray}
Thus by the Lipschitz continuity of $\BS C$ (see
Assumption~\ref{assum:LS}(iii)) and Theorem~\ref{thm.1},  we have 
\begin{align*}
    &\|D_{t-\tau}(u)-D_{t-\tau}(v)\|_2 =\sum_{s_1,s_2 \in \Z} (\BS D(v))_{t,s_1} (\BS C(u)-\BS C(v))_{s_1,s_2} (\BS D(u))_{s_2,\tau} \\
    &\leq K \mathcal{K}^2 \sum_{s_1,s_2
      \in \Z} \zeta(t-s_1)^{\kappa-1}
      \frac{|u-v|}{\gu{s_2-s_1}^\kappa}
      \zeta(s_2-\tau)^{\kappa-1}\\
&= K \mathcal{K}^2 \sum_{s_1,s_2 \in \Z} \zeta(s_1)^{\kappa-1}
  \frac{|u-v|}{\gu{s_2+\tau-t-s_1}^\kappa} \zeta(s_2)^{\kappa-1} \\
    &\leq 49K\mathcal{K}^2 |u-v| \zeta(\tau-t)^{\kappa-1},
\end{align*}
where the last inequality follows from Lemma~\ref{lem.SSR.F1} and $\mathcal{K}$ is finite constant, independent of $u,v,t,\tau$. 
This proves \eqref{lem.D.cont}.

To prove \eqref{lem.Dn.smooth}, we note that using the classical inverse matrix expansion (analogous to
(\ref{eq:DuDv})) we have
\begin{eqnarray*}
 \BS D^{(N)}-\BS D(t/N)
&=&  \BS D^{(N)}\left(\BS C(t/N)- \BS C^{(N)}\right)\BS D(t/N).
\end{eqnarray*}
Theorem~\ref{thm.1}
 gives bounds for the entries in $\BS D(t/N)$
and $\BS D^{(N)}$. On the other hand, Assumption~\ref{assum:LS} gives the bound
\begin{eqnarray*}
\|\left(\BS C(t/N)- \BS C^{(N)}\right)_{s_1,s_2}\|_{2} &\leq&
                                                              \|\left(\BS C(t/N)-
  \BS C(s_1/N)\right)_{s_1,s_2}\|_{2}   +\|\left(\BS C(s_1/N)-
  \BS C^{(N)}\right)_{s_1,s_2}\|_{2}  \\
 &\leq& K\left(\min\left(\frac{|t-s_1|}{N\gu{s_1-s_2}^{\kappa}},
 \frac{2}{\gu{s_1-s_2}^{\kappa}}\right) + \frac{1}{N\gu{s_1-s_2}^{\kappa-1}}\right).
\end{eqnarray*}
Substituting these bounds into $[\BS D^{(N)}\left(\BS C(t/N)-
  \BS C^{(N)}\right)\BS D(t/N)]_{t,\tau}$ gives 
\begin{eqnarray*}
  &&  \|(\BS D^{(N)}-\BS D(t/N))_{t,\tau}\|_2 \\
    &\leq&  
    K\mathcal{K}^2\!\! \sum_{s_1,s_2\in \mathbb{Z}}\!\! \zeta(t-s_1)^{\kappa-1}\! \left( \gu{s_1-s_2}^{-\kappa}\min(\frac{|t-s_1|}{N},2) + \frac{1}{N\gu{s_1-s_2}^{\kappa-1}}\right) 
    \zeta(\tau-s_2)^{\kappa-1}\nonumber\\
 &\leq& K\mathcal{K}^2  \min\bigg(\sum_{s_1,s_2\in \mathbb{Z}} \zeta(t-s_1)^{\kappa-2} 
    \times \frac{1}{N\gu{s_1-s_2}^{\kappa}} \times 
    \zeta(\tau-s_2)^{\kappa-1},\\
    &&2\sum_{s_1,s_2\in \mathbb{Z}} \zeta(t-s_1)^{\kappa-1} 
    \times \frac{1}{\gu{s_1-s_2}^{\kappa}} \times 
    \zeta(\tau-s_2)^{\kappa-1}\bigg) \\
&& + K\mathcal{K}^2\sum_{s_1,s_2\in \mathbb{Z}} \zeta(t-s_1)^{\kappa-1} 
    \times \frac{1}{N\gu{s_1-s_2}^{\kappa-1}}\times
    \zeta(\tau-s_2)^{\kappa-1}\\
    &\leq&  98K\mathcal{K}^2 \zeta(t-\tau)^{\kappa-2}\min(1/N,2\zeta(t-\tau))
    \label{eq.thm.proof.2},
\end{eqnarray*}
where the last bound follows from Lemma \ref{lem.SSR.F1}. This proves \eqref{lem.Dn.smooth}.
\end{proof}

\begin{proof}[Proof of  equation \eqref{eq.claim.diff}]
By using (\ref{eq:DuDv}) we have 
\begin{eqnarray*}
D_r(u) - D_{r}(v) = 
\sum_{s_1,s_2\in \mathbb{Z}}D_{s_1}(u)[C_{s_1}(u)-C_{s_2}(v)]
D_{s_2-r}(v).
\end{eqnarray*}
Let $h\in \mathbb{R}\backslash\{0\}$, and substitute 
$v=u+h$ and $u=u$ into the above to give 
\begin{eqnarray*}
[D_{r}(u) - D_{r}(u+h)]/h = 
\sum_{s_1,s_2\in \mathbb{Z}}D_{s_1}(u)
\frac{[C_{s_1}(u)-C_{s_2}(u+h)]}{h}
D_{s_2-r}(u+h).
\end{eqnarray*}
Taking the limit $h\rightarrow 0$
(and using dominated convergence to exchange limit and sum)
gives the entry-wise matrix derivative
\begin{eqnarray*}
 \frac{d D_r(u)}{du} 
    &=&  -\sum_{s_1,s_2}D_{s_1}(u)\frac{d C_{s_1-s_2}(u)}{du}D_{s_2-r}(u)
\end{eqnarray*}
and the bound
\begin{eqnarray*}
 \left\|\frac{d D_r(u)}{du}\right\|_2 
    &\leq&  \sum_{s_1,s_2} \|D_{s_1}(u)\|_2\|\frac{d C_{s_1-s_2}(u)}{du}\|_2\|D_{s_2-r}(u)]\|_2
    \leq \mathcal{K} \zeta(r)^{\kappa-1},
\end{eqnarray*}
where the last inequality follows 
from  Theorem \ref{thm.D.smooth}, the condition $\sup_{u} \| \frac{d C_r(u)}{du}\|_2 \leq K\zeta(r)^{\kappa-1}$ and
Lemma~\ref{lem.SSR.F1}.
\end{proof}

\begin{proof}[Proof of Corollary~\ref{cor.var}]
To prove the result we start with the inverse matrix 
$\BS D^{(N)} = (\BS C^{(n)})^{-1}$ which we show below has simple
easily derivable properties. We then apply
Theorem~\ref{lemma:spectral}, Lemma~\ref{lem.demko}, and
Theorem~\ref{thm.D.smooth} to obtain analogous properties on its
inverse $\BS C^{(N)} = (\BS D^{(n)})^{-1}$. 

 Define the matrix 
\begin{eqnarray*}
\widetilde{\Phi}_{j}(t/N) = 
\left\{
\begin{array}{cc}
I_{p} & j=0 \\
-\Phi_{j}(t/N) & 1\leq j \leq p \\
0 & \textrm{ otherwise }
\end{array}
\right..
\end{eqnarray*}
Using $\{\Phi_{j}(u)\}_{j}$ we define the stationary time 
$X_{t}(u)=\sum_{j=1}^{d}\Phi_j(u)X_{t-j}(u) + \Sigma(u)^{1/2}\varepsilon_{t}$. This has the
inverse (stationary) covariance 
$\BS D(u)=(D_{t-\tau}(u);t,\tau\in \mathbb{Z})$ where 
\begin{eqnarray}
\label{eq:DARu}
D_{t-\tau}(u)=\sum_{\ell=0}^{d}\widetilde{\Phi}_{\ell}(u)^{\top} \Sigma(u)^{-1} \widetilde{\Phi}_{(t-\tau)+\ell}(u).
\end{eqnarray}
The corresponding inverse spectral density is 
$f(\omega;u)^{-1} = \sum_{r\in \mathbb{Z}}D_{r}(u)\exp(ir\omega)$.
Under the stated conditions on the roots associated with $\{\Phi_{j}(u)\}_{r}$ we have that for some 
$\gamma_{1}$ and $\gamma_{2}$ that
$0 <\gamma_{1}\leq\inf_{u}\inf_{\omega}\lambda_{\inf}(f(\omega;u)^{-1}) \leq \sup_{u}\sup_{\omega}\lambda_{\sup}(f(\omega;u)^{-1})\leq \gamma_2 <\infty$ and thus the eigenvalues of 
$\BS D(u)$ are uniformly bounded away from $\gamma_{1}$ and $\gamma_{2}$. Let $\BS C(u) = \BS D(u)^{-1}
=(C_{t-\tau};t,\tau\in \mathbb{Z})$. 
Then by using 
Lemma \ref{lem.demko}  we have
\begin{eqnarray}
\label{eq:Cexponental}
\sup_{u}\|C_{r}(u)\|_{2}\leq \mathcal{K}\rho^{|r|}
\end{eqnarray}
for some $0<\rho <1$. Further, by using \eqref{lem.D.cont} (applied to exponential decay rather than polynomial decay) we have $\|C_{r}(u)-C_{r}(v)\|_2\leq \mathcal{K}\rho^{|r|}|u-v|$.

Using the Cholesky decomposition it can be shown that the inverse covariance is
$\BS D^{(N)} = (D_{t,\tau};t,\tau\in \mathbb{Z})$ where 
\begin{eqnarray}
\label{eq:DARt}
D^{(N)}_{t,\tau}=\sum_{\ell=0}^{d}
  \widetilde{\Phi}_{\ell}\left(\frac{t+\ell}{N}\right)^{\top}
  \Sigma\left(\frac{t+\ell}{N}\right)^{-1} 
\widetilde{\Phi}_{(t-\tau)+\ell}\left(\frac{t+\ell}{N}\right).
\end{eqnarray}
The Lipschitz conditions on  $\Phi_{j}(\cdot)$ 
together with (\ref{eq:DARu}) and (\ref{eq:DARt}) 
imply that  $D^{(N)}_{t,\tau}$ is
approximated by $D_{t-\tau}(t/N)$. I.e.
\begin{eqnarray*}
|D^{(N)}_{t,\tau} - D_{t-\tau}(t/N)\|_{2}
\leq \left\{
\begin{array}{cc}
\frac{\mathcal{K}}{N} & |t-\tau|\leq d \\
0 & |t-\tau| > d
\end{array}
\right..
\end{eqnarray*}
Now by using the above 
and Theorem~\ref{lemma:spectral} for large enough $N$ the conditions in 
Assumption \ref{assum:LS} hold (in terms of the inverse covariance). Therefore for sufficiently large $N$, the rate $\|C_{t,\tau}^{(N)}\|_{2} \leq \mathcal{K}\rho^{|t-\tau|}$ follows from Lemma \ref{lem.demko}. Further, the conditions in
Theorem~\ref{thm.D.smooth} hold and we have
\begin{eqnarray*}
\|C_{t,\tau}^{(N)} - C_{t-\tau}(t/N)\|\leq \mathcal{K}\frac{\rho^{|r|}}{N},
\end{eqnarray*}
which gives $\|C_{t,\tau}^{(N)}-C_{t-\tau}(t/N)\|_{2} \leq \mathcal{K}\frac{\rho^{|t-\tau|}}{N}$. 
Thus we have proved the result. 
\end{proof}

\begin{proof}[Proof of Theorem \ref{cor.Dn.smooth}]
The result uses the bounds $\|[\BS
C^{(N)}(-\infty,T)^{-1}]_{s_1,s_2}\leq \mathcal{K}\zeta(s_1-s_2)^{\kappa-1}$
and $\|[\BS
C(-\infty,T;u)^{-1}]_{s_1,s_2}\leq \mathcal{K}\zeta(s_1-s_2)^{\kappa-1}$. The assertion follows by the same steps as in the proof of 
Theorem~\ref{thm.D.smooth}.
\end{proof}

\begin{proof}[Proof of Theorem~\ref{lemma:smoothVAR}]
To prove the result we start with the following identities 
\begin{eqnarray}
\label{eq:Phidef}
\Phi_{T,j}^{(N)} &=&  -([\BS C^{(N)}(-\infty,T)^{-1}]_{T,T})^{-1}[\BS   C(-\infty,T)^{-1}]_{T,T-j} \\
\textrm{ and } \nonumber
\Phi_{j}(u) &=&  -([\BS
  C(-\infty,T;u)^{-1}]_{T,T})^{-1}{[\BS
  C(-\infty,T;u)^{-1}]_{T,T-j}}
  \end{eqnarray}
  where $\BS C^{(N)}(-\infty,T) = (C_{t,\tau}^{(N)};t,\tau \leq T)$
  and $\BS C(-\infty,T;u) = (C_{t,\tau}(u);t,\tau\leq T)$.
These identities  together with Theorem
\ref{cor.Dn.smooth} will be used to 
prove the result. 

We first obtain a bound for $\|\Sigma_{T}^{(N)} -  \Sigma(T/N)
\|_2$. We note that 
\begin{eqnarray*}
&& \Sigma_{T}^{(N)} -  \Sigma(T/N) = ([\BS
                                    C^{(N)}(-\infty,T)^{-1}]_{T,T})^{-1}
                                    -  ([\BS
                                    C(-\infty,T;T/N)^{-1}]_{T,T})^{-1}
  \\
&=&([\BS C(-\infty,T;T/N)^{-1}]_{T,T})^{-1} ([\BS
    C(-\infty,T;T/N)^{-1}]_{T,T} \\
    &&- [\BS C^{(N)}(-\infty,T)^{-1}]_{T,T})
 ([\BS C^{(N)}(-\infty,T)^{-1}]_{T,T})^{-1}.
\end{eqnarray*}
Thus by using Theorem \ref{cor.Dn.smooth} (with $t=T$ and $\tau=T$)  we have 
\begin{eqnarray}
 \|\Sigma_{t}^{(N)} -  \Sigma(t/N)\|_{2} 
&\leq &\|([\BS C(-\infty,T;T/N)^{-1}]_{T,T})^{-1}\|_2 \|([\BS C^{(N)}(-\infty,T)^{-1}]_{T,T})^{-1}\|_2 \nonumber\\
&&\times \|[\BS
    C(-\infty,T;T/N)^{-1} - \BS C^{(N)}(-\infty,T)^{-1}]_{T,T}\|_2\nonumber\\
    &
\leq& \mathcal{K}N^{-1}. \label{eq:sigma2}
\end{eqnarray}
This proves the first part of (i)

To prove the second part of (ii), we use 
 (\ref{eq:Phidef}) to give the decomposition
$\Phi_{t,j}^{(N)} - \Phi_{j}(t/N) = J_{1} + J_{2}$,
where 
\begin{eqnarray*}
J_1 &=& -\left[ ([\BS C^{(N)}(-\infty,T)^{-1}]_{T,T})^{-1} - [\BS
  C(-\infty,t;t/N)^{-1}]_{T,T})^{-1}\right]
  [\BS C^{(N)}(-\infty,T)^{-1}]_{T,T-j}, \\
J_2 &=& -  ([\BS
  C(-\infty,T;T/N)^{-1}]_{T,T})^{-1}
  \left[ [\BS C^{(N)}(-\infty,T)^{-1} - \BS C(-\infty,T;T/N)^{-1}]_{T,T-j}\right].
\end{eqnarray*}
First we bound $J_1$ this gives 
\begin{eqnarray*}
\|J_{1}\|_2 &\leq& \left\| ([\BS C^{(N)}(-\infty,T)^{-1}]_{T,T})^{-1} - ([\BS
  C(-\infty,T;T/N)^{-1}]_{T,T})^{-1}\right\|_{2}\\
  &&\times 
  \|[\BS C^{(N)}(-\infty,T)^{-1}]_{T,T-j}\|_{2} \\
 &\leq&
                   \mathcal{K}\frac{1}{N}\zeta(0)^{\kappa-1}\times \zeta(j)^{\kappa-1}.
\end{eqnarray*}
where we have used the bounds in Theorem \ref{thm.1}  and
(\ref{eq:sigma2}) in the above. 
Using a similar argument (and Theorem \ref{cor.Dn.smooth} (with $t=T$ and $\tau=T-j$) ) we have 
\begin{eqnarray*}
\|J_{2}\|_2 &\leq& \|([\BS
  C(-\infty,T;T/N)^{-1}]_{T,T})^{-1}\|_{2}
  \left\| [\BS C^{(N)}(-\infty,T)^{-1} - \BS
                   C(-\infty,T;T/N)^{-1}]_{T,T-j}\right\|_{2}\\
&\leq&  \mathcal{K} \zeta(j)^{\kappa-2} \min(2\zeta(j),1/N).
\end{eqnarray*}
Altogether this gives 
$
\|\Phi_{T,j}^{(N)} - \Phi_{j}(T/N)\|_{2} \leq \mathcal{K} \zeta(j)^{\kappa-2} \min(2\zeta(j),1/N).
$
Thus we have proved the second part of (i). The proof for
(ii) follows a similar method as given in the proof of Theorem~\ref{thm.D.smooth}, and we omit the details. 
\end{proof}

We now prove Theorem~\ref{thm.part}. 
To prove this result we will use the alternative representation of the 
covariance operator $\BS C^{(N)}$ defined in Remark
\ref{remark:alternativegroup}. With this in mind,  we define the sub-operators  
 $\BS C^{(e,f)}:\ell_{2}\rightarrow \ell_{2}$ which are infinite
dimensional matrices where 
$[\BS C^{(e,f)}]_{t,\tau} = \cov[X_{t,N}^{(e)},
X_{t,N}^{(f)}]$. Note that to reduce cumbersome notation, we have dropped the $N$
from the definition ${\BS C}^{(e,f)}$. We also define the
corresponding ``stationary'' matrix operators $\BS C^{(e,f)}(u):\ell_{2}\rightarrow \ell_{2}$, where 
 $[\BS C^{(e,f)}(u)]_{t,\tau} = \cov[X_{t}^{(e)}(u),
X_{t}^{(f)}(u)]$. This representation is instrumental in proving the result below.

\begin{proof}[Proof of Theorem~\ref{thm.part}]
We first prove (\ref{eq:partial1}) and (\ref{eq:partial2}).  We start by obtaining an expression for 
\begin{eqnarray*}
\var\left[X_{t}^{(c)|\shortminus\{a,b\}};t\in \mathbb{Z}, c\in
  \{1,2\}\right] &=& (\Delta_{t,\tau,N}^{\shortminus\{a,b\}};t,\tau\in
                     \mathbb{Z}) \\
\textrm{ and }
\var\left[X_{t}(u)^{(c)|\shortminus\{a,b\}};t\in \mathbb{Z}, c\in
  \{1,2\}\right] &=& (\Delta_{t-\tau}^{\shortminus\{a,b\}}(u);t,\tau\in \mathbb{Z}).
\end{eqnarray*}
To simplify notation, and
without loss of generality, we focus on the case
$a=1,b=2$. We will represent the above in terms of block matrices of $\BS C^{(N)}$ and
$\BS C(u)$. We define 
 $\BS A^{(1,2)}: \ell_{2,2} \to \ell_{2,2}$, $\BS B^{(1,2)}: \ell_{2,p-2} \to \ell_{2,2}$ and $\BS E^{(1,2)} : \ell_{2,p-2}\to \ell_{2,p-2}$ where 
\begin{eqnarray*}
\BS A^{(1,2)}&=&\begin{pmatrix}
\BS C^{(1,1)}& \BS C^{(1,2)} \\
\BS C^{(2,1)}& \BS C^{(2,2)} \\
\end{pmatrix},
\BS B^{(1,2)}=\begin{pmatrix}
\BS C^{(1,3)}& \dots & \BS C^{(1,p)} \\
\BS C^{(2,3)}& \dots & \BS C^{(2,p)} \\
\end{pmatrix} \\
\textrm{ and }
\BS E^{(1,2)} &=& 
(\BS C^{(e,f)};e,f\in \{3,\ldots,p\}).
\end{eqnarray*}
Analogously, we define $\BS A^{(1,2)}(u),\BS B^{(1,2)}(u),\BS E^{(1,2)}(u)$.  It is clear the operators 
$\BS A^{(1,2)}$, $\BS B^{(1,2)}$ and $\BS E^{(1,2)}$ are comprised of an infinite number of 
$2\times 2$, $2\times (p-2)$ and $(p-2)\times (p-2)$
matrices respectively. To denote these sub-matrices 
we use the following notation. Suppose $\BS H : \ell_{2,p_1} \to \ell_{2,p_2}$ for some $p_1,p_2$ then $[\BS H]_{t,\tau}:=(I_{p_1} \otimes e_t)^\top B^{(1,2)}(I_{p_2} \otimes e_\tau)$ refers to their $p_1\times p_2$-dimensional submatrices.

It is well known that the conditional covariance of $X_{t,N}^{(c)}$
and $X_{t}^{(c)}(u)$ can be represented as the  Schur complement
\begin{eqnarray*}
\var\left[X_{t,N}^{(c)|\shortminus\{1,2\}};t\in \mathbb{Z}, c\in \{1,2\}\right]
&=& \BS A^{(1,2)}-\BS B^{(1,2)} (\BS E^{(1,2)})^{-1}(\BS B^{(1,2)})^\top
\end{eqnarray*}
and 
\begin{eqnarray*}
\var\left[X_{t}(u)^{(c)|\shortminus\{1,2\}};t\in \mathbb{Z}, c\in \{1,2\}\right]
&=& 
\BS A^{(1,2)}(u)-\BS B^{(1,2)}(u) (\BS E^{(1,2)}(u))^{-1}(\BS B^{(1,2)}(u))^\top.
\end{eqnarray*}
Then, we have  
\begin{eqnarray}
\Delta_{t,\tau,N}^{\shortminus\{a,b\}} &=& [\BS A^{(1,2)}-\BS B^{(1,2)}
                                           (\BS E^{(1,2)})^{-1}(\BS B^{(1,2)})^\top]_{t,\tau} \nonumber\\
\textrm{ and }
\Delta_{t-\tau}^{\shortminus\{a,b\}}(u)
&=& [\BS A^{(1,2)}(u)-\BS B^{(1,2)}(u)
    (\BS E^{(1,2)}(u))^{-1}(\BS B^{(1,2)}(u))^\top]_{t,\tau}. \label{eq:Gammattau} 
\end{eqnarray}

We use the above representations to prove (\ref{eq:partial1}). Using
(\ref{eq:Gammattau}) we have 
\begin{eqnarray*}
\|\Delta_{t,\tau,N}^{\shortminus\{a,b\}} -
  \Delta_{t-\tau}^{\shortminus\{a,b\}}(t/N)\|_2 \leq J_1 + J_2 + J_3 + J_4
\end{eqnarray*}
where 
\begin{eqnarray*}
J_1 &=& \| (\BS A^{(1,2)}-\BS A^{(1,2)}(t/N))_{t,\tau}\|_2 \\
J_2 &=& \|(\BS B^{(1,2)} (\BS E^{(1,2)})^{-1}(\BS B^{(1,2)}-\BS B^{(1,2)}(t/N))^\top)_{t,\tau}\|_2\\
J_3&=& \|(\BS B^{(1,2)} ((\BS E^{(1,2)})^{-1}-(\BS E^{(1,2)}(t/N))^{-1})(\BS B^{(1,2)}(t/N))^\top)_{t,\tau}\|_2\\
J_4&=& \|((\BS B^{(1,2)}-\BS B^{(1,2)}(t/N)) (\BS E^{(1,2)}(t/N))^{-1}(\BS B^{(1,2)}(t/N))^\top)_{t,\tau}\|_2.
\end{eqnarray*}

Under Assumption~\ref{assum:LS} and by using Theorem \ref{thm.1} we
bound the terms above (the proof is in the spirit of the
proof of Theorem \ref{thm.D.smooth}).  Assumption~\ref{assum:LS}
directly implies 
\begin{eqnarray*}
J_{1} = \| (\BS A^{(1,2)}-\BS A^{(1,2)}(t/N))_{t,\tau}\|_2 \leq K\frac{1}{N\gu{t-\tau}^{\kappa-1}}.
\end{eqnarray*}
The  bounds for   $J_{2}, J_3$ and $J_{4}$ are more involved, however all three follow a similar strategy. 
We  focus on obtaining a bound for $J_{3}$. Using standard matrix
multiplication it can be seen that 
\begin{eqnarray}
\label{eq:J3}
&&J_{3} =  \|\sum_{s_1,s_2\in \mathbb{Z}}
  [\BS B^{(1,2)}]_{t,s_1}
          [(\BS E^{(1,2)})^{-1} - (\BS E^{(1,2)}(t/N))^{-1}]_{s_1,s_2}[\BS B^{(1,2)}(t/N))^\top]_{s_2,\tau}\|_{2} \nonumber\\
&\leq&\!\!\sum_{s_1,s_2\in \mathbb{Z}}
  \|[\BS B^{(1,2)}]_{t,s_1}\|_{2}\cdot\|[(\BS E^{(1,2)})^{-1} - (\BS E^{(1,2)}(t/N))^{-1}]_{s_1,s_2}\|_2
  \cdot\|(\BS B^{(1,2)}(t/N))^\top]_{s_2,\tau}\|_{2}\qquad
\end{eqnarray}
To bound $ \|[\BS B^{(1,2)}]_{t,s_1}\|_{2}$ and
$\|(\BS B^{(1,2)}(t/N))^\top]_{s_2,\tau}\|_{2}$ we simply use Assumption
\ref{assum:LS}, which immediately gives 
\begin{eqnarray}
\label{eq:B12bound}
\|[\BS B^{(1,2)}]_{t,s_1}\|_{2}
\leq K\gu{t-s_1}^{-\kappa} \textrm{ and }
\|(\BS B^{(1,2)}(t/N))^\top]_{s_2,\tau}\|_{2} \leq
K\gu{s_2-\tau}^{-\kappa}. 
\end{eqnarray}
The bound for 
$\|[(\BS E^{(1,2)})^{-1} -
(\BS E^{(1,2)}(t/N))^{-1}]_{s_1,s_2}\|_2$ needs a little more work. 
We first note that the covariance operator
$\BS E^{(1,2)}$ is a suboperator of
${\BS C}^{(N)}$, thus it satisfies
Assumption~\ref{assum:LS} where $\BS E^{(1,2)}(u)$ is its locally stationary
approximation. Therefore we can apply the 
results of 
Theorem \ref{thm.D.smooth} to $(\BS E^{(1,2)})^{-1}$ and this gives 
\begin{eqnarray}
\label{eq:E12bound1}
\|((\BS E^{(N),(1,2)})^{-1}-(\BS E^{(1,2)}(s_1/N))^{-1})_{s_1,s_2}\|_2\leq \mathcal{K} \zeta(s_1-s_2)^{\kappa-2}\min(1/N,2\zeta(s_1-s_2))
\end{eqnarray}
and 
\begin{eqnarray}
\label{eq:E12bound2}
\|((\BS E^{(1,2)}(s_1/N))^{-1} - (\BS E^{(1,2)}(t/N))^{-1})_{s_1,s_2}\|_2\leq \mathcal{K} |s_1-t|\zeta(s_1-s_2)^{\kappa-1}/N.
\end{eqnarray}
Substituting (\ref{eq:B12bound}), (\ref{eq:E12bound1}) and
(\ref{eq:E12bound2}) into (\ref{eq:J3}) we have 
\begin{eqnarray*}
   J_{3}  &\leq& \mathcal{K}K^2\sum_{s_1,s_2\in \Z} 
\frac{1}{\gu{t-s_1}^\kappa}\times \left(\zeta(s_1-s_2)^{\kappa-2}\min(1/N,2\zeta(s_1-s_2))\right. \\
&& +\left. |s_1-t|\zeta(s_1-s_2)^{\kappa-1}/N\right) \frac{1}{\gu{s_2-\tau}^\kappa}\\
    &\leq& 2\times(49)K^2 \mathcal{K} \zeta(t-\tau)^{\kappa-2}\min(1/N,\zeta(t-\tau))=:\mathcal{K} \zeta(t-\tau)^{\kappa-2}\min(1/N,\zeta(t-\tau)),
\end{eqnarray*}
where the last line follows from Lemma \ref{lem.SSR.F1}.

To bound $J_2$, we use 
 Theorem \ref{thm.1} to give,
$\|[(\BS E^{(1,2)})^{-1}]_{s_1,s_2}\|\leq
\mathcal{K}\zeta(s_1-s_2)^{(\kappa-1)}$.
This together with (\ref{eq:B12bound}), using the bounds stated in Assumption \ref{assum:LS}(iii) and
following the same proof as above we can show that 
\begin{eqnarray*}
J_{2} \leq \mathcal{K} \zeta(t-\tau)^{\kappa-1}/N \textrm{
             and } J_{4} \leq \mathcal{K} \zeta(t-\tau)^{\kappa-1}/N.
\end{eqnarray*}
Altogether the bounds for $J_{1},J_{2},J_{3}$ and $J_{4}$ prove 
\begin{eqnarray*}
\|\Delta_{t,\tau,N}^{\shortminus\{a,b\}} - \Delta_{t-\tau}^{\shortminus\{a,b\}}(t/N)\|_2 
\leq \mathcal{K} \zeta(t-\tau)^{\kappa-2}\min(1/N,\zeta(t-\tau))
\end{eqnarray*}
thus proving (\ref{eq:partial1}). The proof of (\ref{eq:partial2})
follows a similar technique. 

Finally, the proofs for (\ref{eq:partial3}) and (\ref{eq:partial4})
are the same as the proofs for (\ref{eq:partial1}) and
(\ref{eq:partial2}), thus we omit the details.
\end{proof}





\end{document}